\documentclass[11pt]{article}

\usepackage{bm}
\usepackage{amsmath, amssymb}

\makeatletter

\@addtoreset{equation}{section}
\def\section{\@startsection{section}{1}{0pt}{-3.5ex plus -1ex minus
 -.2ex}{2.3ex plus .2ex}{\large\bf}}
\def\subsection{\@startsection{subsection}{2}{\z@}{-3.25ex plus -1ex minus
 -.2ex}{1.5ex plus .2ex}{\normalsize\bf}}
\makeatother

\newtheorem{remark}{Remark}[section]
\newtheorem{definition}{Definition}[section]
\newtheorem{theorem}{Theorem}[section]
\newtheorem{corollary}{Corollary}[section]
\newtheorem{proposition}{Proposition}[section]
\newtheorem{lemma}{Lemma}[section]
\newtheorem{assumption}{Assumption}[section]
\newtheorem{condition}{Condition}[section]

\newcommand{\delfrac}[2]{{\displaystyle\frac{\partial #1}{\partial #2}}}
\newcommand{\iprod}{\, \mbox{\raise.3ex\hbox{\tiny $\bullet$}}\, }
\title{Probabilistic representation of weak solutions to a parabolic 
boundary value problem on a non-smooth domain}

\author{Masaaki Tsuchiya\footnote{e-mail: mtsuchiya02@yahoo.co.jp}
\hspace{1.8cm} 
(Kanazawa University
\footnote{Kanazawa University, Kanazawa 920-1192, Japan}\ )
\and 
Appendix B:\ Hajime Kawakami\footnote{e-mail: kawakami@math.akita-u.ac.jp} 
\ (Akita University
\footnote{Mathematical Science Course, 
Akita University, Akita 010-8502, Japan}\ )}

\date{}

\begin{document}
\maketitle

\begin{abstract}
\noindent
The probabilistic representation of weak solutions to a parabolic boundary value problem is established in the following framework. 
The boundary value problem consists of a second order parabolic equation defined on a time-varying Lipschitz domain in a Euclidean space and of a mixed boundary condition composed of a Robin and the homogeneous Dirichlet conditions. 
It is  assumed that the time-varying domain is included in a fixed smooth domain and that a certain  part of the boundary of the time-varying domain is also included in the boundary of the fixed domain, namely the fixed boundary.  
The Robin condition is imposed on a part of the boundary included in the fixed one and the Dirichlet condition 
on the rest of the boundary.  Such a parabolic boundary value problem always has a unique weak solution for given data; however it does not possess a classical or strong solution in general, even in the case of constant coefficients. 
The stochastic solution to the boundary value problem is also considered.  
Then, by verifying the equality between both the solutions on the domain and its lateral boundary, the probabilistic representation for the weak solution is obtained. 
Thus, it is ensured that, for the weak solution, the stochastic solution provides a version which is continuous up to the lateral boundary of the domain except the border between the places imposed the Robin and Dirichlet conditions. 
As an application, it is shown the continuity property 
of a functional (cost function) related to an optimal stopping problem motivated by an inverse problem determining the shape of a domain.  
\medskip

\medskip

\medskip

\noindent
\textbf{AMS Subject Classification}:
\begin{tabular}[t]{l} 
primary: 60G46, 60J60;  \\
secondary: 35D30, 35K20, 35R37, 60J55 
\end{tabular} 

\noindent
\textbf{Keywords}: 
\begin{tabular}[t]{l} 
parabolic boundary value problem, weak solution, 
stochastic solution, 
\\
oblique reflecting diffusion process, local time on the boundary, 
\\
shape identification inverse problem  
\end{tabular} 
\end{abstract}

\newpage

\numberwithin{equation}{section} 

\section{Introduction}
\label{sect.intro}

 The probabilistic representation of solutions to second order parabolic equations is a useful tool to analyze these solutions (e.g. boundary sensitive analysis of solutions in \cite{CGK06}, numerical analysis of solutions in \cite{CPS98}) and it has several applications (e.g. see 
\cite{Fre85, Roz05, Zha12, Kaw15} and references therein). 
For classical solutions to elliptic and parabolic equations, such probabilistic representation has been studied extensively in \cite{Fre85}.  
Parabolic equations encountered  in applications have often no solutions with sufficient regularity. 
Several authors have studied the probabilistic representation of  weak solutions to parabolic equations on a whole Euclidean space via backward stochastic differential equations  (e.g. \cite{Roz05, Zha12}).

This paper concerns with the probabilistic representation of weak solutions to a parabolic equation with a mixed boundary condition on a time-varying Lipschitz domain in a Euclidean space: the boundary condition is composed of a Robin and the homogeneous Dirichlet conditions.  
Such a parabolic equation is treated in \cite{Kaw10} to study an inverse problem determining the shape of a time-varying domain; in general, it does not possess a classical or strong solution, even if the equation has constant coefficients. 

The probabilistic representation is considered in the following framework:    
The time-varying domain is included in a fixed smooth domain and a certain  part of the boundary of the time-varying domain is also included in the boundary of the fixed domain, namely the fixed boundary. 
The Robin condition is imposed on a part of the boundary included in the fixed boundary and the homogeneous Dirichlet condition is imposed on the rest of the boundary; we call the place setting   the Robin condition (resp. Dirichlet condition) the Robin part (resp. Dirichlet part). 
Accordingly the Dirichlet part may be varied with time and both of the parts may be adjoining. 
As a basic stochastic process for the representation, we take the corresponding oblique reflecting diffusion process on the closure of the fixed domain to the parabolic equation. 

Then the probabilistic representation of a weak solution is given by the \textit{stochastic solution}, 
which is the expectation of the quantity obtained by applying It\^o's formula formally to the weak solution and the diffusion process.  
The equality between both the solutions within the domain is verified based on an appropriate approximation of the weak solution; the boundedness and regularity of the weak solution and a Poincar\'e type inequality with weight (see Lemma 4.1 of \cite{Kaw10}) play key roles.
To show the equality on the lateral boundary, we need to examine the continuity up to the boundary for the stochastic solution, since the boundary values of weak solutions are given by their traces on the boundary. 

The continuity of the stochastic solution is examined based on the coupled martingale formulation for the oblique reflecting diffusion process (see \cite{Gra88, Tsu94, Kaw00} for such martingale formulations), because we have to use the continuity property of the local time with respect to random parameter and to treat the weak convergence property of functionals which contain the local time. 
Then we ensure that the stochastic solution is continuous up to the lateral boundary except the border between the Robin and Dirichlet parts:  
it is derived from the continuity property of the local time and the entrance time to the Dirichlet part with respect to random parameter, and it is further done from a certain estimate for the distribution of the entrance time in the case the process starts from a point near the Dirichlet part.

The paper is organized as follows. In Section 2, we provide a necessary setting and assumptions: they are concerned with the considered domain and the parabolic boundary value problem.  The main result (Theorem \ref{main_result} and Corollary \ref{main_result_cor}) is stated with proof in Section 3 after describing the coupled martingale problem for the oblique reflecting diffusion process. In Section 4, as a simple application of the main result, we show the continuity of a functional whose argument is a domain belonging to an admissible class. In Appendix A, we provide several technical details of the proof for  the main result:  the boundedness and regularity of weak solutions; no attainability of the process to closed null sets of the boundary; the continuity property of the entrance time to the Dirichlet part; the estimate for the distribution of the entrance time; a variant of the mapping theorem in weak convergence of probability measures. 
Finally, in Appendix B, based on a Bayesian framework,                       we give a brief exposition on the shape identification inverse problem mentioned in Section 4.  

\section{Preliminaries and a parabolic boundary value problem}
\label{sect.frame_assump}

We start with introducing necessary notations for describing the parabolic boundary value problem. 
It is treated in the backward form adapted to the probabilistic consideration. 
In what follows, we treat the time-varying domain in $\bm{R}^{n}$ as a non-cylindrical domain in the time-space $\bm{R}^{1+n} = \bm{R}_t \times \bm{R}_x^n$.  
For a subset $G$ of $\bm{R}^{1+n},$ 
denote the $t$--section of $G$ by  
\[
G(t) : = G \cap \left(\{t\} \times \bm{R}^n\right). 
\] 
If necessary, 
we identify $G(t)$ with the set $\left\{x: (t,x) \in G \right\}$ in $\bm{R}^n.$
For a bounded open subset $G$ of $\bm{R}^{1+n},$ we consider its parabolic boundary, lateral boundary and ceiling (the time reversed notion of the bottom in the forward form) in the backward form and denote by 
$\partial_P G,$ $\partial_L G$ and $\partial_C G,$ respectively, which are subsets of the boundary $\partial G$ in $\bm{R}^{1+n}.$    
Here, for the precise definition (in the forward form), we refer to page 1787 in \cite{BHL97}.
For $0 \leq a < b$ and $E \subset \bm{R}^n,$ denote by $E_{a,b} := (a, b) \times E$  the cylindrical set determined by $a , b$ and $E;$ in particular, we put $E_{T}:= E_{0,T}.$ 

Now we denote the usual Sobolev space on an open set $E$ in $\bm{R}^{n}$ with a nonnegative integral order $r$ by $H^{r}(E)$; that is, 
\[
H^{r}(E):= \{f \in L^{2}(E): \partial^{\alpha}_{x}f \in L^{2}(E) ~\mbox{ for } ~|\alpha| \leq r\}, 
\]
where $\partial_{x}^{\alpha}f$ indicates the weak derivative of $f$ with respect to multi-index $\alpha.$
Moreover, we recall Sobolev spaces and another function space on a domain $G$ in the time-space 
$\bm{R}^{1+n}:$ 
For nonnegative integers $r$ and $s$, we set 
\[
H^{r,s}(G):= \{g \in L^{2}(G): \partial_{t}^{\alpha}g \in L^{2}(G) ~\mbox{ for } ~0 \leq \alpha \leq r; ~
\partial_{x}^{\beta}g \in L^{2}(G) ~\mbox{ for } ~|\beta| \leq s\},
\]
and, when $G \subset [a, b] \times \bm{R}^{n}$, we set  
\[
V^{0,1}(G):= \{g \in H^{0,1}(G): \bar{g} \in C([a, b]; L^{2}(\bm{R}^{n}))\}, 
\]
where let $\bar{g}(t, x):= g(t, x)1_{G(t)}(x).$ 
These function spaces are equipped with the following norms, respectively:
\begin{align*}
\|f\|_{H^{r}(E)} &= \left\{\sum_{|\alpha|\leq r}\|\partial^{\alpha}_{x}f\|_{L^{2}(E)}^{2}\right\}^{1/2}, \\
\|g\|_{H^{r,s}(G)} &= \left\{\sum_{0\leq\alpha\leq r}\|\partial_{t}^{\alpha}g\|_{L^{2}(G)}^{2} + 
\sum_{|\beta|\leq s}\|\partial_{x}^{\beta}g\|_{L^{2}(G)}^{2}\right\}^{1/2}, \\
\|g\|_{V^{0,1}(G)} &= \max_{a\leq t\leq b}\|g(t, \cdot)\|_{L^{2}(G(t))} + \|\partial_{x}g\|_{L^{2}(G)}.
\end{align*}

For a Lipschitz domain $E$ in $\bm{R}^{n},$ denote by $\gamma_{\partial E}$ the boundary trace operator from $H^{1}(E)$ into $L^{2}(\partial E);$ it can be given by the following pointwise limit (see \cite{Eva92}, p. 133): for $v \in H^{1}(E)$ 
\begin{equation}
\gamma_{\partial E}v(x) = \lim_{r \to 0}\, -\!\!\!\!\!\!\int_{B_{r}(x) \cap E} v(y) dy  ~~~~(\mathcal{H}^{n-1} \mbox{--a.e.}~x \in \partial E),
\label{trace}
\end{equation}
where $B_{r}(x) = \{y \in \bm{R}^{n}: |y - x| < r\}$ ~($r > 0$) and $\mathcal{H}^{n-1}$ stands for the $(n-1)$--dimensional Hausdorff measure (then $\left.\mathcal{H}^{n-1}\right|_{\partial E}$ is the surface measure on $\partial E$). 
More generally, if $\Delta$ is an open subset of $\partial E,$ $V$ is a neighborhood of $\Delta$ in $\overline{E}$ and $v \in H^{1}(V \cap E)$, then the limit in (\ref{trace}) by replacing $E$ with $V \cap E$ exists for $\mathcal{H}^{n-1}$--a.e. $x \in \Delta.$ 
Hence, if necessary, we use the notion of boundary trace in this extended sense; the limit is also  denoted by $\gamma_{\Delta}v(x)$.   

We need the notion of regularized distance for the Euclidean distance. 
Following Theorem 2 on page 171 of \cite{Ste70},  we summarize the fundamental properties. 
For a closed set $F$ in a Euclidean space, there exists a function 
$\mathit{\Delta}_{F}$, the regularized distance for the Euclidean distance $\rho_{F}: = d(\cdot, F),$ such that (i) there are positive constants $c_{1}$, $c_{2}$ satisfying $c_{1}\rho_{F} \leq \mathit{\Delta}_{F} \leq c_{2}\rho_{F}$ for every closed set $F$, where $c_{1}$ is taken as an absolute constant and $c_{2}$ is taken as a constant depending only on and increasing with the dimension of the Euclidean space; (ii) $\mathit{\Delta}_{F}$ is a $C^{\infty}$ function in $F^{c}$ and for any multi-index $\alpha$ there is a positive constant $B_{\alpha}$ satisfying $\displaystyle{\left|\partial^{\alpha}\mathit{\Delta}_{F}\right| \leq B_{\alpha}\rho_{F}^{1-|\alpha|}}$ for every closed set $F.$
For a subset $A$ of the Euclidean space, let $\tilde{\rho}_{A}:= c_{1}^{-1}\mathit{\Delta}_{\bar{A}}.$  
Then $\rho_{A} \leq \tilde{\rho}_{A} \leq c^{\ast}\rho_{A}$ for every subset $A$ with $c^{\ast}:= c_{1}^{-1}c_{2} > 1$; so that,  in what follows, we use $\tilde{\rho}_{A}$ instead of $\mathit{\Delta}_{\bar{A}}$ in  $\bm{R}^{n}$ and $\bm{R}^{1+n}$, where $c_{2}$ is taken as the constant in the time-space. 

This paper concerns with a domain which varies with time of $(0, T)$ in a fixed bounded Lipschitz domain $\Omega$ in $\bm{R}^n.$  
Therefore  the time-varying domain, say $D,$ is described as a domain of $\bm{R}^{1+n}$ 
in the cylindrical domain $\Omega_T.$   
Moreover we take a Lipschitz dissection of the boundary 
$\Gamma := \partial \Omega$ of $\Omega;$ $\Gamma = 
\Gamma' \cup \Pi \cup \Gamma^{\prime\prime}.$   
Here by a Lipschitz dissection of 
$\Gamma$ we mean that $\Gamma'$ and $\Gamma^{\prime\prime}$ are disjoint open subsets of $\Gamma$ and $\Pi$ is an $(n-2)$--dimensional closed Lipschitz surface in $\Gamma$ (see \cite{McL00} for the precise definition). 
In the following, we further need the notion of Lipschitz lateral surface with extreme time edges in the time-space (see Condition 3.1 in \cite{Kaw10} for the definition).   

We first impose the following condition on the considered domain $D.$  
{\condition
\label{cond_D_1}
\begin{enumerate}
\item[(i)]  
$D$ is a bounded Lipschitz domain in $\bm{R}^{1+n};$
\item[(ii)]  
for each $t \in [0, T],$ $D(t)$ is a non-empty Lipschitz domain in $\bm{R}^n$ and $D(t) \subset \Omega,$ where 
$D(0) := \left(\overline{D}(0)\right)^{\circ}$ (the interior of $\overline{D}(0)$) and
$D(T) := \left(\overline{D}(T)\right)^{\circ};$
\item[(iii)]   
$\partial D \setminus \partial_{P}D = D(0)$ and $\partial_C D = D(T);$
\item[(iv)]   
\label{cond_D_1.v}
the lateral boundary $\partial_{L}D$  of $D$ is a Lipschitz lateral surface with extreme time edges and includes the set 
$[0, T] \times(\Gamma' \cup \Pi),$ and further 
there is  an open subset $U$ in $\Gamma$ 
satisfying $\Pi \subset U,$ $[0, T] \times U \subset \partial_LD$ and 
$\left\{\partial_L D \setminus 
\left([0, T] \times(\Gamma' \cup \Pi)\right)\right\} 
\cap ([0, T] \times U) \subset [0, T] \times \Gamma^{\prime\prime}.$ 
\end{enumerate}
}
\vspace{0.3cm}

Therefore the time-varying portion of the lateral boundary of $D$ is
included in the set $\Sigma := \partial_L D \setminus 
\left([0, T] \times(\Gamma' \cup \Pi)\right).$ 
In the paper, 
we mainly treat the case where $\Sigma$ has no connected component included in $[0, T] \times \Omega$ and $\Pi \neq \emptyset$, although we can treat some of general types of division of 
$\partial_{L}D$ into two parts: one includes the time-varying portion and the other does not. 
The case treated here is a typical one where each of both the parts has a component adjoining each other; such a case needs a complicated treatment to obtain the result. 

Let $A \equiv A(t,x),$ 
$\bm{a} \equiv \bm{a}(t,x),$ 
$\bm{b} \equiv \bm{b}(t,x)$ 
and $a \equiv a(t,x)$ be an $n \times n$ real matrix-valued 
function,  
$n$-dimensional real vector-valued functions and a real-valued 
function 
defined on $\overline{\Omega_{T}}$ respectively, and $\sigma(t, x)$ a real-valued function defined on 
$\overline{\Gamma'_{T}}$.  
Furthermore suppose that $A(t,x)$ is symmetric and positive definite.  
For such $A$, $\bm{a},$ $\bm{b}$ 
and $a$, define a second-order differential operator $\mathcal{L} \equiv \mathcal{L}_{x}(t)$ on $\Omega_{T}$: 
\[
\mathcal{L}u(t,x) :=   
\nabla_x \iprod (A(t,x) \nabla_x u(t,x) + 
\bm{a}(t,x) u(t,x)) 
- \bm{b}(t,x) \iprod \nabla_x u(t,x)
- a(t,x)u(t,x). 
\]

As mentioned in the beginning, we treat the parabolic equation in the backward form and hence we consider  
weak solutions to the terminal-boundary value problem for the parabolic equation. 
Define the parabolic operator $\mathcal{P} \equiv \mathcal{P}_{x}(t)$ on $D$ in the backward form,  
the conormal derivative $\partial/\partial \mathcal{N} \equiv \partial/\partial \mathcal{N}_{x}(t)$ 
relative to $\mathcal{L}$ and the boundary operator $\mathcal{B} \equiv \mathcal{B}_{x}(t)$ on  $\Gamma'_{T}$ by 
\begin{align*}
\mathcal{P}u(t,x)
& := \delfrac{u}{t}(t,x)
+ \mathcal{L}u(t,x) ~~\mbox{for}~ (t, x) \in D, 
\\
\delfrac{u}{\mathcal{N}}(t,x) 
& :=    
- \left[A(t,x) \nabla_x u(t,x) + \bm{a}(t,x) u(t,x)
\right] \iprod \bm{n}(t,x) ~~\mbox{for}~ (t, x) \in \Gamma'_{T}, 
\nonumber
\\
\mathcal{B}u(t,x) 
& := \delfrac{u}{\mathcal{N}}(t,x) + \sigma(t,x)u(t,x) ~~\mbox{for}~ (t, x) \in \Gamma'_{T}, 
\nonumber
\end{align*}
where $\bm{n} = \bm{n}(t,x)$ is the inward unit normal vector at 
$x\in \Gamma' \subset \partial (D(t)),$ the boundary of $D(t)$ in $\{t\} \times \bm{R}^n.$  
We note that $(\partial_{L}D)(t) = \partial(D(t))$ under Condition \ref{cond_D_1}. 

Now consider the following terminal-boundary value problem on the 
domain $D$:
\begin{equation}
\begin{cases}
\mathcal{P}u(t,x)= - \nabla_x \iprod \bm{f}(t,x) + f(t, x)  
& \text{if} ~ (t,x) \in D
\\
\mathcal{B}u(t,x) = - \psi(t,x)
& \text{if} ~ (t,x) \in \Gamma'_T
\\
u(t,x)=0 & \text{if} ~ (t,x) \in \Sigma
\\
u(T,x)=h(x) & \text{if} ~ x \in D(T).
\end{cases}
\label{eq_1}
\end{equation}
In what follows, assume that $h$ is extended onto $\Omega$ with value zero outside $D(T),$ and $\bm{f}, ~f$ are extended onto $\Omega_{T}$ with value zero outside $D.$ 
On the problem (\ref{eq_1})
we impose the following assumption. 
{\assumption
\label{assump_2}
The coefficients and source terms
of the terminal-boundary value problem 
(\ref{eq_1}) satisfy the following conditions.
\begin{enumerate}
\item[(i)] 
\label{item.assump_2.i}
There exists a constant $\nu>0$ such that 
for every $(t,x)\in \overline{\Omega_T}$ and   
$\xi \in \bm{R}^n$
\[ 
A(t,x)\xi\iprod\xi\ge\nu|\xi|^2. 
\] 
\item[(ii)] 
\label{item.assump_2.ii}
The coefficients of the parabolic operator $\mathcal{P}$ and $\sigma$ 
are bounded measurable on $\overline{\Omega_T}$ and $\overline{\Gamma'_T},$ 
respectively.   
\item[(iii)] 
Suppose that $h \in L^2(\Omega),$  
$f \in L^2(\Omega_T),$ $\bm{f} = (f_{1}, \ldots, f_{n}) \in L^2(\Omega_T)^n,$ 
$\psi \in L^2(\Gamma'_T),$
where $L^{2}(\Omega_{T})$ (resp. $L^{2}(\Gamma'_{T})$) stands for the $L^{2}$ space with respect to the measure $dt \times dx$ (resp. $dt \times S(dx)$, here $S(dx)$ is the surface measure on $\Gamma$). 
Moreover, 
$\gamma_{\Gamma'}\bm{f} =  (\gamma_{\Gamma'} f_{1}, \ldots, \gamma_{\Gamma'} f_{n})$ is defined and it belongs to $L^2(\Gamma'_T)^n,$ where $\gamma_{\Gamma'} f_{j}(t, x):= (\gamma_{\Gamma'}f(t, \cdot))(x)$ for $(t, x) \in \Gamma'_{T}$ ~$(j = 1, 2, \ldots, n).$ 
\end{enumerate}
}
\medskip

The notion of weak solutions to (\ref{eq_1}) 
is defined as follows.
{\definition
\label{definition.subsect.general.parabolic.weak.1}
A function $u(t,x) \in V^{0,1}(D)$ 
is called a \textit{weak solution} in $V^{0,1}(D)$ to the 
terminal-boundary value problem 
(\ref{eq_1}) 
if it satisfies the following conditions:  
\begin{enumerate}
\item[(i)] 
$\left.\,u \right|_{\Sigma}=0,$ that is, 
$\gamma_{\partial(D(t))}u(t,\cdot) = 0$ 
on ${\Sigma}(t)$ 
for almost every $t\in (0,T)$;  
\item[(ii)] 
\label{item.definition.subsect.general.parabolic.weak.1.ii}
for every $\eta \in H^{1,1}(D)$ 
with $\left.\,\eta \right|_{\Sigma}=0$ and $\eta(0, \cdot) = 0$ 
\begin{multline}
\int_D u\, \partial_{t}\eta\, dt dx - \int_{D(T)} h(x)\eta(T,x)\, dx 
\\
+\int_D\left\{\left(A\nabla_xu  + \bm{a}u \right) 
\iprod \nabla_x \eta
+ \left(\bm{b} \iprod \nabla_xu + au \right)\eta \right\}\, dt dx 
\\ 
+ \int_D\left(\bm{f} \iprod \nabla_x \eta + f \eta \right)\, dt dx
+ \int_{\Gamma'_T}\left(\psi + \gamma_{\Gamma'}\bm{f}\iprod \bm{n} \right) \gamma_{\Gamma'}\eta\, dt S(dx) 
\\
+ \int_{\Gamma'_T}\sigma \gamma_{\Gamma'} u\cdot \gamma_{\Gamma'} \eta\, dt S(dx)
= 0. 
\label{eq.weak.form.I}
\end{multline}
\end{enumerate}
}

By the results in \cite{Kaw10}, we see that, under Assumption \ref{assump_2}, 
there exists a unique weak solution in $V^{0,1}(D)$ 
to the problem (\ref{eq_1}).  
Moreover, in the proof of the main theorem, the boundedness of the weak solution plays a key role. 
The boundedness result is shown in Appendix \ref{subsect.appendix.boundedness}.
There we use the following norms: for a measurable function $\varphi(t, x)$ on $\Omega_{T}$ or 
on $\Xi_{ T}$ (determined by a measurable subset $\Xi$ of $\Gamma$) and for $q, ~r \geq 1,$ let
\begin{align*}
\|\varphi\|_{q,r;\Omega_{T}}&:= \left\{\int_{0}^{T}\left(\int_{\Omega}\left|\varphi(t, x)\right|^{q} dx\right)^{\frac{r}{q}} dt \right\}^{\frac{1}{r}}, \\
\|\varphi\|_{q,r;\Xi_{T}}&:= \left\{\int_{0}^{T}\left(\int_{\Xi}\left|\varphi(t, x)\right|^{q} S(dx)\right)^{\frac{r}{q}}dt \right\}^{\frac{1}{r}}.
\end{align*}

\section{Main result}
\label{sect.main_result}

\subsection{Coupled martingale problem} 
\label{subsect.coup.mart}

To state and prove the main result, we use  the coupled martingale formulation for the oblique reflecting diffusion process associated with the parabolic and boundary operators. 
For this purpose, we need a further assumption on the domain $\Omega,$ the set $\Sigma,$ and the coefficients and the source terms in (\ref{eq_1}). 
{\assumption
\label{assump_3}
\begin{enumerate}
\item[(i)]
$\Omega$ is a bounded $C^{2,\alpha}$ domain in $\bm{R}^n$ for an $\alpha \in (0, 1].$
\item[(ii)]
Let $\Sigma_{1}:= \Sigma \cap ([0, T] \times \Gamma)$ and $\Sigma_{2}:= \Sigma \cap ([0, T] \times \Omega).$ For $j = 1, 2$ denote by $\delta \Sigma_{j}$,  
$\stackrel{\circ}{\Sigma}_{j}$ and 
$\overline{\Sigma}_{j}$ the boundary of $\Sigma_{j}$, the interior of 
$\Sigma_{j}$ and the closure of $\Sigma_{j}$ in $\partial_{L}D,$ respectively, and  
let $\Lambda = \delta \Sigma_{2}.$ 
Suppose that   
\begin{align*}
\delta \Sigma_{1} = \Lambda\cup([0, T] \times \Pi), ~\overline{\Sigma}_{1} = \Lambda\cup \stackrel{\circ}{\Sigma}_{1} \cup([0, T] \times \Pi), \\
d(\Lambda, [0, T] \times \Pi) > 0, ~ S(\Lambda(t)) = 0 ~\mbox{for a.e.}~t \in [0, T], 
\end{align*}
where $d$ indicates the Euclidean distance. 
\item[(iii)]
$A$ and $\bm{a}$ are Lipschitz continuous in $t \in [0, T]$ and their Lipschitz constants are uniformly bounded in $x \in \overline{\Omega}$; furthermore, $A$ and $\bm{a}$ are differentiable in $x \in \overline{\Omega}$ and their derivatives are $(\alpha/2, \alpha)$--H\"older continuous on $\overline{\Omega_{T}}$. In addition, $\bm{b}$ and $a$ are also 
$(\alpha/2, \alpha)$--H\"older continuous on $\overline{\Omega_{T}}$, and $\sigma$ is Lipschitz continuous on 
$\overline{\Gamma^{\prime}_{T}}$.  
\item[(iv)]
Suppose that $\bm{f} = \bm{0}$, $f$ is continuous on $\overline{\Omega_{T}}$, 
$\psi$ is Lipschitz continuous on $\overline{\Gamma^{\prime}_{T}}$ and $h$ is Lipschitz continuous on $\overline{\Omega}$. 
\end{enumerate} 
}
\medskip
\medskip

In what follows, assume that such coefficients and source terms defined 
on $[0, \infty) \times \overline{\Omega}$ in such a way as $A(t, x) = A(T, x),$ $\bm{a}(t, x) = \bm{a}(T, x)$, and so on for $t \geq T$ and $x \in \overline{\Omega}.$ 
According to the extension of the coefficients and the source terms, we also extend the domain $D$ and a subset $\Xi$ of the lateral boundary $\partial_{L}D$ to $D_{(\infty)}$ and $\Xi_{(\infty)}$ along the time axis as follows:
\[
D_{(\infty)}: = D \cup ([T, \infty) \times D(T)), ~~~ \Xi_{(\infty)}: = \Xi \cup ([T, \infty) \times \Xi(T)).
\]
In addition, if necessary, $\sigma$ and $\psi$ are appropriately extended onto $\Gamma_{T}$ with the same property and further extended onto $\Gamma_{(\infty)}$ in the same way as above.  
Under Assumption \ref{assump_3}, the operator $\mathcal{L}$ and the boundary operator 
$\mathcal{B}$ are rewritable in the non-divergence and in the oblique derivative forms, respectively:               
\[
\mathcal{L} = {\rm{Tr}}(A(t, x)\nabla_{x}^{2}) + \bm{c}(t, x)\iprod\nabla_{x} + c(t, x); ~ 
\mathcal{B} = \bm{\beta}(t, x)\iprod\nabla_{x} + \gamma(t, x).
\]
Further set 
\[
\mathcal{L}_{0} \equiv \mathcal{L}_{0;x}(t):= \mathcal{L}_{x}(t) - c(t, x); ~\mathcal{B}_{0} \equiv \mathcal{B}_{0;x}(t):= \mathcal{B}_{x}(t) - \gamma(t, x); ~\mathcal{P}_{0}:= \frac{\partial}{\partial t} + \mathcal{L}_{0}.
\]
Then we consider the coupled martingale problem associated with $(\mathcal{P}_{0}, \mathcal{B}_{0})$ as in \cite{Tsu94, Kaw00}.  Let $\bm{W}:= C([0, \infty) \to \overline{\Omega}),$  $\bm{V}:= C([0, \infty) \to [0, \infty)) \cap \{increasing\, functions\}$ and 
$\bm{U}:= \bm{W} \times \bm{V}.$ 
We assume that $\bm{W}$, $\bm{V}$ and $\bm{U}$ each are equipped with the locally uniform convergence topology. 
Denote generic elements of $\bm{U}$, $\bm{W}$ and $\bm{V}$ by $\bm{\omega}$, $\bm{w}$ and 
$\bm{v}$, respectively; that is, $\bm{\omega} = (\bm{w}, \bm{v}).$ 
Then put 
\[
X(t, \bm{\omega}) \equiv X(t, \bm{w}):= \bm{w}(t), ~ L(t, \bm{\omega}) \equiv L(t, \bm{v}):= \bm{v}(t)
\]
and define the $\sigma$ fields $\mathcal{U}_{t}^{s}~~(0 \leq s \leq t < \infty)$ and $\mathcal{U}$ generated by $(X(r, \cdot), L(r, \cdot))$ $(r \in [s, t])$ and by $(X(r, \cdot), L(r, \cdot))$ $(r \in [0, \infty))$, respectively. 
We now introduce two spaces of test functions for the coupled martingale problem: $C^{1,2}_{b}([0, \infty) \times \overline{\Omega})$ stands for the space of functions $f$ on $[0, \infty) \times \overline{\Omega}$ having bounded continuous derivatives $\partial^{\alpha}_{t}\partial^{\beta}_{x}$ with $2\alpha + |\beta| \leq 2$ and $C_{0}(\Omega)$ 
for the space of continuous functions on $\Omega$ with $\mathrm{spt} f \subset \Omega$. 
By \cite{Tsu94, Kaw00}, under Assumptions \ref{assump_2} and \ref{assump_3}, for each 
$(s, x) \in [0, \infty) \times \overline{\Omega},$ the coupled martingale problem associated with $(\mathcal{P}_{0}, \mathcal{B}_{0})$ has a unique solution starting from $(s, x)$, say $P_{s,x};$ 
that is, $P_{s,x}$ is a unique probability measure on $(\bm{U},~ \mathcal{U})$ satisfying the following conditions:
\begin{enumerate}
\item[(i)] ~ 
$\displaystyle{
P_{s,x}(X(r) = x ~\mbox{ and } ~L(r) = 0 ~\mbox{for all}~0 \leq r \leq s) = 1
}$;
\item[(ii)] ~ 
$\displaystyle{
P_{s,x}\left(L(t) = \int_{s}^{t}1_{\Gamma}(X(r))L(dr) ~\mbox{for all}~ t \geq s\right) = 1
}$; 
\item[(iii)]
for every $g \in C^{1,2}_{b}([0, \infty) \times \overline{\Omega}),$ the process $\{M_{g}(t); t\geq s\}$ defined by 
\[
\hspace*{-2.5cm}~~M_{g}(t):= g(t, X(t)) - g(s, X(s)) - \int_{s}^{t}\mathcal{P}_{0}g(r, X(r))dr 
- \int_{s}^{t}\mathcal{B}_{0}g(r, X(r))L(dr) 
\]
is a martingale on the filtered probability space $(\bm{U}, \mathcal{U}, P_{s,x}; \mathcal{U}_{t}^{s}).$ 
\end{enumerate}

We know that, under the condition (i), the condition (iii) is equivalent to the fact: for every $g \in C^{1,2}_{b}([0, \infty) \times \overline{\Omega}),$ the process $\{M_{g}(t); t\geq s\}$ is a 
martingale on the filtered probability space $(\bm{U}, \mathcal{U}, P_{s,x}; \mathcal{U}_{t}^{0})$; that is, the filtration $\{\mathcal{U}_{t}^{s}\}$ may be replaced with $\{\mathcal{U}_{t}^{0}\}.$   
Moreover we note that the condition (ii) is equivalent to the condition  
$\displaystyle{
P_{s,x}\left(\int_{s}^{t}1_{\Omega}(X(r))L(dr) = 0 ~\mbox{for all}~ t \geq s\right) = 1}$, that is, 
for every $h \in C_{0}(\Omega)$, $\displaystyle{
P_{s,x}\left(\int_{s}^{t} h(X(r))L(dr) = 0 ~\mbox{for all}~ t \geq s\right) = 1}$. This is also equivalent to the condition that, for every $h \in C_{0}(\Omega)$, the process 
$\displaystyle{\left\{\tilde{M}_{h}(t):= \int_{s}^{t}h(X(r))L(dr); t \geq s\right\}}$
 is a 
$\{\mathcal{U}_{t}^{s}\}$--martingale, because $\{\tilde{M}_{h}(t); t \geq s\}$  is a 
$\{\mathcal{U}_{t}^{s}\}$--adapted continuous bounded variation process with $\tilde{M}_{h}(s) = 0$. 
We also observe that, for $s \geq 0$, if we set $L_{s}(t):= L(t) - L(t \wedge s)$ $(t \geq 0)$, then $L_{s}(dr) = L(dr)$ for $dr \subset [s, \infty)$.  
Accordingly we have another equivalent formulation of the coupled martingale problem as in the following remark. 
{\remark
\label{remark_1}
A probability measure $P$ on $(\bm{U},~ \mathcal{U})$ is called a solution to the coupled martingale problem associated with $(\mathcal{P}_{0}, \mathcal{B}_{0})$ starting from $(s, x) \in [0, \infty) \times \overline{\Omega},$ if it satisfies the following conditions: 
\begin{enumerate}
\item[(i)]  $\displaystyle{
P(X(r) = x ~\mbox{ and } ~L(r) = 0 ~\mbox{for all}~0 \leq r \leq s) = 1}$;
\item[(ii)]
for every $g \in C^{1,2}_{b}([0, \infty) \times \overline{\Omega})$ and $h \in C_{0}(\Omega),$ the process $\{(M_{g}(t), \tilde{M}_{h}(t)); t\geq s\}$ is a two-dimensional martingale on the filtered probability space $(\bm{U}, \mathcal{U}, P; \mathcal{U}_{t}^{0}).$ 
\end{enumerate}
(There are various equivalent types of the condition in (ii) just above as in Theorem 4.2.1 of \cite{Str79}.) 
By this unified martingale formulation, the fundamental facts on the martingale problem described in \S 6.1 and \S 6.2 of \cite{Str79} can be directly applied to the coupled martingale problem (see also \cite{Eth86}, Chap. \!4 for a general treatment of martingale problems and the Markov property of their solutions). 
}
\medskip

 The coupled martingale formulation for diffusion processes with reflecting barrier is crucial for our treatments, since the continuity property in $(t, \bm{\omega})$ of the canonical process $L(t, \bm{\omega})$ for the local time on the boundary is essentially used.
By virtue of the fundamental facts derived from the unified martingale formulation,   
noted in Remark \ref{remark_1}, the uniqueness of solutions implies that the family $\{P_{s,x}; (s, x) \in [0, \infty) \times \overline{\Omega}\}$ is strong Markov in the sense as in Theorem 6.2.2 of \cite{Str79}.

We further need the uniform exponential integrability of $L(T)$: 
\begin{equation}
\sup_{(s,x)\in [0,T]\times\overline{\Omega}}E_{s,x}\left[e^{\lambda L(T)}\right] < \infty ~\mbox{ for every } \lambda > 0
\label{unif.expo}
\end{equation}
(see Proof of Theorem 4.3 in \cite{Kaw00} and also Proposition 3.5 in \cite{CGK06}).

Moreover, under Assumptions \ref{assump_2} and \ref{assump_3}, the transition probability $P_{s,x}(X(t) \in dy)$ has a transition density $p(s, x; t, y)$ with respect to the $n$--dimensional Lebesgue measure $\mathcal{L}^{n}(dy) = dy$, 
which is a fundamental solution to the parabolic equation $\mathcal{P}_{0}u = 0$ on 
$\overline{\Omega}$ with the boundary condition $\mathcal{B}_{0}u = 0$. 
Then, the transition density and its first order derivatives in $x$ have the usual upper Gaussian bound: for every $T > 0,$ there exist positive constants $K$ and $C$ such that 
\begin{equation}
\left|\partial_{x}^{\alpha}p(s, x; t, y)\right| \leq K(t - s)^{-\frac{n+|\alpha|}{2}}\exp\left\{- C \frac
{|x - y|^{2}}{t - s}\right\}
\label{G.bound}
\end{equation}
for every $0 \leq s < t \leq T,$ $x, y \in \overline{\Omega}$ and $|\alpha| \leq 1$  (see \cite{Gar84, Kaw00}). 

\subsection{Statement and proof of the main result} 
\label{subsect.stat.main_result}

For a Borel set $G$ in $[0, \infty) \times \overline{\Omega}$ and $s \geq 0,$ we set 
\begin{align*}
\sigma_{s}(G)&: = \inf\{t \geq s: (t, X(t)) \in G\} =  \inf\{t \geq s: X(t) \in G(t)\}, \\
\sigma_{s+}(G)&: = \inf\{t > s: (t, X(t)) \in G\} =  \inf\{t > s: X(t) \in G(t)\}; 
\end{align*}
in particular, in the case $G = [0, \infty) \times E,$ $\sigma_{s}(G) =  \inf\{t \geq s: X(t) \in E\} =:\sigma_{s}(E)$ and $\sigma_{s+}(G) =  \inf\{t > s: X(t) \in E\} =:\sigma_{s+}(E).$ 
Whenever $(s, X(s)) \notin G,$ it holds that $\sigma_{s}(G) = \sigma_{s+}(G)$ and that 
$\sigma_{s}(G)(\bm{\omega}) = \sigma_{0}(G)(\bm{\omega})$ and $\sigma_{s+}(G)(\bm{\omega}) = \sigma_{0+}(G)(\bm{\omega})$ for $\bm{\omega}$ with $X(r, \bm{\omega}) = X(s, \bm{\omega})$ ($0 \leq r \leq s$).  
Moreover, for a closed set $G,$ $\sigma_{s}(G)$ is a stopping time relative to the filtration $\{\mathcal{U}^{s}_{t}; t \geq s\}$ and, for an open set $G,$ $\sigma_{s+}(G)$ is a stopping time relative to the filtration $\{\mathcal{U}^{s}_{t+}; t \geq s\}.$

Let 
\[
Z_{s}(t) := \int_{s}^{t}c(r, X(r))dr + \int_{s}^{t}\gamma(r, X(r))L(dr) \quad\mbox{for}~ 0 \leq s \leq t. 
\]
Then, under Condition \ref{cond_D_1}, Assumptions \ref{assump_2} and \ref{assump_3},  we have the following probabilistic representation of the weak solution $u = u^{D}$ to the problem (\ref{eq_1}). 

{\theorem
\label{main_result}
The equality 
\begin{align}
u^{D}(s, x) = &- E_{s,x}\left[\int_{s}^{\sigma_{s}(\overline{\Sigma}_{(\infty)})\wedge T}\exp Z_{s}(t)f(t, X(t))dt\right] \nonumber \\ 
& - E_{s,x}\left[\int_{s}^{\sigma_{s}(\overline{\Sigma}_{(\infty)})\wedge T}\exp Z_{s}(t)\psi(t, X(t))L(dt)\right]  \nonumber \\ 
& + E_{s,x}\left[\exp Z_{s}(T)h(X(T)); \sigma_{s}(\overline{\Sigma}_{(\infty)}) > T\right] 
\label{repre}
\end{align}
holds in the following each case: 
\begin{enumerate}
\item[(i)]  $ds \times dx$--a.e. $(s, x) \in D;$
\item[(ii)]  $ds \times S_{s}(dx)$--a.e. $(s, x) \in \partial_{L}D$, where $S_{s}(dx)$ indicates the surface measure on $(\partial_{L}D)(s)$. 
\end{enumerate}
Furthermore, the right hand side of (\ref{repre}) is continuous on $D \cup \left(\partial_{L}D \setminus ([0, T] \times \Pi)\right)$. 
}
\medskip

\noindent
\textbf{Proof.}
We first verify the equality (\ref{repre}) in the case (i). 
Throughout the proof, we extend $u^{D}$ onto $\Omega_{T}$ with value zero outside $D$ and use the same symbol for the extension.
 Choose a $C^{\infty}$ increasing function $\tilde{\lambda}$ on $[0, \infty)$ such that 
$\tilde{\lambda}(\xi) = 0$ for $0 \leq \xi \leq 1,$ 
$\tilde{\lambda}^{'}(\xi) > 0$ for $1 < \xi < 2$ and $\tilde{\lambda}(\xi) = 1$ for $\xi \geq 2.$
In the following, we take $\varepsilon > 0$ and $\delta > 0$ satisfying the conditions: 
\begin{enumerate}
\item[(1)] 
The set $\{y \in \overline{\Omega}: \rho_{\Pi}(y) \leq 2\varepsilon\}$ is contained in a tubular  neighborhood of $\Gamma$ in $\overline{\Omega}$ (see \cite{Gil01}, Chap. 14, Appendix for such a neighborhood).  
\item[(2)]
$[0, T] \times \{y \in \overline{\Omega}: \rho_{\Pi}(y) \leq 2\varepsilon\} \subset D \cup \partial_{L}D$. 
\item[(3)] 
$\delta < \frac{1}{2} d(\Sigma_{2(\infty)}, [0, \infty) \times (\Pi \cup \Gamma'))$.
\end{enumerate}
For such $\varepsilon$ and $\delta$, let 
\begin{align*}
\tilde{\mu}_{\delta}(t, y)&:= \tilde{\lambda}(\delta^{-1}\tilde{\rho}_{\Sigma_{2(\infty)}}(t, y)), ~~ u^{D}_{\centerdot\, \delta}(t, y) := \tilde{\mu}_{\delta}(t, y)u^{D}(t, y), \\        
\chi_{\varepsilon}(y) &:=  \tilde{\lambda}(\varepsilon^{-1}\tilde{\rho}_{\Pi}(y)), ~~ 
u^{D}_{\varepsilon \,\centerdot}(t, y) := \chi_{\varepsilon}(y) u^{D}(t, y), \\
u^{D}_{\varepsilon,\delta}(t, y)&:= \chi_{\varepsilon}(y) u^{D}_{\centerdot\, \delta}(t, y) = 
\tilde{\mu}_{\delta}(t, y)u^{D}_{\varepsilon \,\centerdot}(t, y)
= \chi_{\varepsilon}(y)\tilde{\mu}_{\delta}(t, y)u^{D}(t, y).
\end{align*}
Using the argument on page 27 in \cite{Kaw10}, we know that  
the family $\{u^{D}(t, \cdot)^{2}; 0\leq t \leq T\}$ is uniformly integrable and that 
for a positive constant  $M,$ $\sup_{t \geq 0}\mathcal{L}^{n}(\Sigma_{2(\infty)}^{[2\delta]}(t)) \leq M\delta$ for sufficiently small $\delta > 0,$ where 
$\Sigma_{2(\infty)}^{[2\delta]} := \{(t, y) \in D_{(\infty)}: \rho_{\Sigma_{2(\infty)}}(t, y) < 2 \delta\}.$ 
Therefore 
\begin{equation}
\max_{0\leq t \leq T}\|u^{D}(t, \cdot) - u^{D}_{\varepsilon,\delta}(t, \cdot)\|_{L^{2}(\Omega)} \longrightarrow 0 ~~\mbox{ as }~\delta \downarrow 0~\mbox{ and } \varepsilon \downarrow 0.  
\label{approx.sol.1}
\end{equation}
Moreover, we see that  
\begin{equation}
u^{D}_{\varepsilon,\delta} \in H^{1,2}(\Omega_{T});
\label{reg.sol}
\end{equation}
which is verified in Appendix \ref{subsect.appendix.regularity}. 
Hence, applying the smoothing procedure in Theorem 3 on page 127 in \cite{Eva92} separately to the time variable and  the space variable, 
we can take an approximating sequence 
$\left\{u^{D}_{\varepsilon,\delta;m}\right\}_{m=1}^{\infty}$ to the function $u^{D}_{\varepsilon,\delta}$ such that 
\begin{align*}
&u^{D}_{\varepsilon,\delta;m} \in H^{1,2}(\Omega_{T})\cap C^{\infty}\left(\overline{\Omega_{T}}\right)~~  (m = 1, 2, \ldots), \\   
&u^{D}_{\varepsilon,\delta;m} \to u^{D}_{\varepsilon,\delta} \mbox{ in } H^{1,2}(\Omega_{T})  \mbox{ as } m \to \infty.
\end{align*}
We recall briefly the smoothing procedure for functions to use again it later.  
It is based on a partition of unity for 
$\overline{\Omega_{T}}$ subordinate to a finite open covering of $\overline{\Omega_{T}}$ and on mollifying the localized ones of the functions on each open set consisting of the covering; in particular, in each open set including some part of the boundary, the mollifying is combined with a certain inward parallel displacement.    
Moreover, by Theorem \ref{boundedness}, we see that 
\[
\|u^{D}\|_{\infty;D}:= \mathrm{ess\, sup}_{(t,y)\in D}|u^{D}(t, y)| < \infty; 
\]
then, from the smoothing procedure for $u^{D}_{\varepsilon,\delta},$ it follows that 
\[
\|u^{D}_{\varepsilon,\delta;m}\|_{\infty;D} \leq \|u^{D}\|_{\infty;D} 
\]
for $\varepsilon > 0, \delta > 0$ and $m = 1, 2, \ldots$ . 
For $\eta > 0,$ let 
\begin{align*}
\Sigma_{2(\infty)}^{\eta} &:= \{(t, y) \in D_{(\infty)}: \tilde{\rho}_{\Sigma_{2(\infty)}}(t, y) < \eta\}, \\
\Pi^{\eta} &:= \{y \in \Omega: \tilde{\rho}_{\Pi}(y) < \eta\}. 
\end{align*}
Then, for $s \in (0, T)$ and sufficiently small $\varepsilon$ and $\delta$ 
mentioned above, define
\begin{align*}
&\hat{\sigma}_{\delta}:= \sigma_{s}(\overline{\Sigma}_{(\infty)}) \wedge \sigma_{s}\left(\overline{\Sigma_{2(\infty)}^{2\delta}}\right), \\
&\tilde{\sigma}_{\varepsilon}:= \sigma_{s}([0, \infty) \times \overline{\Pi^{2\varepsilon}}) = \sigma_{s}(\overline{\Pi^{2\varepsilon}})
\end{align*}
and, for $0 < s < s^{\prime} < T,$ set  
$
\hat{\sigma}_{\delta}\wedge T \wedge\tilde{\sigma}_{\varepsilon}\vee s^{\prime}:= (\hat{\sigma}_{\delta}\wedge T \wedge\tilde{\sigma}_{\varepsilon})\vee s^{\prime}.
$
It\^{o}'s formula implies 
\begin{align*}
&\exp Z_{s}(\hat{\sigma}_{\delta}\wedge T\wedge\tilde{\sigma}_{\varepsilon}\vee s^{\prime})\,u^{D}_{\varepsilon,\delta;m}(\hat{\sigma}_{\delta}\wedge T\wedge\tilde{\sigma}_{\varepsilon}\vee s^{\prime}, X(\hat{\sigma}_{\delta}\wedge T\wedge\tilde{\sigma}_{\varepsilon}\vee s^{\prime})) \\
&\hspace*{8cm}- \exp Z_{s}(s^{\prime})\,u^{D}_{\varepsilon,\delta;m}(s^{\prime}, X(s^{\prime})) \\
&= \int_{s^{\prime}}^{\hat{\sigma}_{\delta}\wedge T\wedge\tilde{\sigma}_{\varepsilon}\vee s^{\prime}}\exp Z_{s}(t)\, \mathcal{P}u^{D}_{\varepsilon,\delta;m}(t, X(t))dt  \\
&+ \int_{s^{\prime}}^{\hat{\sigma}_{\delta}\wedge T\wedge\tilde{\sigma}_{\varepsilon}\vee s^{\prime}}\exp Z_{s}(t)\, \mathcal{B}u^{D}_{\varepsilon,\delta;m}(t, X(t))L(dt) + \mbox{ a martingale difference}.
\end{align*} 
Now let 
\[
\tilde{D}_{\varepsilon,\delta}:= \{(t, y) \in D: \tilde{\rho}_{\Sigma_{2(\infty)}}(t, y) > 2\delta, ~\tilde{\rho}_{\Pi}(y) > 2\varepsilon\}. 
\]
For $(s, x) \in D,$ we further take $\varepsilon$ and $\delta$ so small that $(s, x) \in 
\tilde{D}_{\varepsilon,\delta}.$ 
Then 
\begin{align}
&E_{s,x}\left[\exp Z_{s}(s^{\prime})\,u^{D}_{\varepsilon,\delta;m}(s^{\prime}, X(s^{\prime})); s^{\prime} < \hat{\sigma}_{\delta}\wedge T\wedge\tilde{\sigma}_{\varepsilon}\right] \nonumber \\
&= - E_{s,x}\left[\int_{s^{\prime}}^{\hat{\sigma}_{\delta}\wedge T\wedge\tilde{\sigma}_{\varepsilon}\vee s^{\prime}}\exp Z_{s}(t)\, \mathcal{P}u^{D}_{\varepsilon,\delta;m}(t, X(t))dt\right]  \nonumber \\
&- E_{s,x}\left[\int_{s^{\prime}}^{\hat{\sigma}_{\delta}\wedge T\wedge\tilde{\sigma}_{\varepsilon}\vee s^{\prime}}\exp Z_{s}(t)\, \mathcal{B}u^{D}_{\varepsilon,\delta;m}(t, X(t))L(dt)\right]  \nonumber \\
&+ E_{s,x}\left[\exp Z_{s}(\hat{\sigma}_{\delta}\wedge T\wedge\tilde{\sigma}_{\varepsilon}\vee s^{\prime})\,u^{D}_{\varepsilon,\delta;m}(\hat{\sigma}_{\delta}\wedge T\wedge\tilde{\sigma}_{\varepsilon}\vee s^{\prime}, X(\hat{\sigma}_{\delta}\wedge T\wedge\tilde{\sigma}_{\varepsilon}\vee s^{\prime})) \right]. 
\label{eq-expec-1}
\end{align}
From now on, denote by $I_{0}(s, x; s^{\prime}, \varepsilon,\delta, m)$ the left hand side and by $I_{j}(s, x; s^{\prime}, \varepsilon,\delta, m)$ the $j$ th term of the right hand side $(j = 1, 2, 3)$ of (\ref{eq-expec-1}).
For showing the assertion (i) of the theorem, it is enough to verify the following equalities: 
\begin{equation}
\widetilde{\lim_{s^{\prime}\downarrow s}}\lim_{\varepsilon\downarrow 0}\lim_{\delta\downarrow 0}\widetilde{\lim_{m\to\infty}} I_{0}(s, x; s^{\prime}, \varepsilon,\delta, m) = u^{D}(s, x) 
\label{fact_0}
\end{equation}
for each $s \in (0, T)$, a.e. $x \in D(s)$,   
and
\begin{align}
&\lim_{s^{\prime}\downarrow s}\lim_{\varepsilon\downarrow 0}\lim_{\delta\downarrow 0}\lim_{m\to\infty} I_{1}(s, x; s^{\prime}, \varepsilon,\delta, m) = - E_{s,x}\left[\int_{s}^{\sigma_{s}(\overline{\Sigma}_{(\infty)})\wedge T}\exp Z_{s}(t)\, f(t, X(t))dt\right],  \label{fact_1}  \\
&\lim_{s^{\prime}\downarrow s}\lim_{\varepsilon\downarrow 0}\lim_{\delta\downarrow 0}\lim_{m\to\infty} I_{2}(s, x; s^{\prime}, \varepsilon,\delta, m) = E_{s,x}\left[\int_{s}^{\sigma_{s}(\overline{\Sigma}_{(\infty)})\wedge T}\exp Z_{s}(t)\, \psi(t, X(t))L(dt)\right], 
\label{fact_2}  \\
&\lim_{s^{\prime}\downarrow s}\lim_{\varepsilon\downarrow 0}\lim_{\delta\downarrow 0}\lim_{m\to\infty} I_{3}(s, x; s^{\prime}, \varepsilon,\delta, m) = E_{s,x}\left[\exp Z_{s}(T)\, h(X(T)); \sigma_{s}(\overline{\Sigma}_{(\infty)}) > T\right],
\label{fact_3}
\end{align}
for $(s, x) \in D.$ 
Here ``$\displaystyle{\widetilde{\lim_{m\to\infty}}}$'' means that there is a subsequence 
$\{m_{k}\}_{k=1}^{\infty}$ (independent of $s$, $x$; $s'$, $\varepsilon$ and 
$\delta$) of $\{m\}_{m=1}^{\infty}$ and the limit takes along the subsequence, and \\
`` $\displaystyle{\widetilde{\lim_{s'\downarrow s}}}$ '' means that for each $s \in (0, T)$ there is a subsequence 
$\{s_{\ell}\}_{\ell=1}^{\infty}$ with $s_{\ell}\downarrow s$ ($\ell \to \infty$) and the limit takes along the subsequence. 

In the following, for a given $\varepsilon > 0$, we take $\delta > 0$ as $\delta < \frac{1}{2} d(\Sigma_{2(\infty)}, [0, \infty) \times (\Pi^{2\varepsilon} \cup \Gamma'))$. \\
Verifying (\ref{fact_0}): First divide $I_{0}(s, x; s^{\prime}, \varepsilon,\delta, m) - u^{D}(s, x)$ as follows:
\begin{align*}
&I_{0}(s, x; s^{\prime}, \varepsilon,\delta, m) - u^{D}(s, x) \\
&= \left\{
E_{s,x}\left[\exp Z_{s}(s^{\prime})\,u^{D}_{\varepsilon,\delta;m}(s^{\prime}, X(s^{\prime})); s^{\prime} < \hat{\sigma}_{\delta}\wedge T\wedge\tilde{\sigma}_{\varepsilon}\right] - 
E_{s,x}\left[u^{D}_{\varepsilon,\delta;m}(s', X(s'))\right] 
\right\}
\\
&+ \left\{E_{s,x}\left[u^{D}_{\varepsilon,\delta;m}(s', X(s'))\right] - E_{s,x}\left[u^{D}_{\varepsilon,\delta;m}(s, X(s'))\right] 
\right\}
\\
&+ \left\{
E_{s,x}\left[u^{D}_{\varepsilon,\delta;m}(s, X(s'))\right] - u^{D}_{\varepsilon,\delta;m}(s, x) 
\right\}
\\
&+ \left\{u^{D}_{\varepsilon,\delta;m}(s, x) - u^{D}(s, x)\right\}.
\end{align*}
In the equality above, 
we denote by $I_{0,j}(s, x; s^{\prime}, \varepsilon,\delta, m)$ 
the $j$ th difference of the right hand side $(j = 1, 2, 3)$ 
and by $I_{0,4}(s, x; \varepsilon,\delta;m)$ 
the fourth difference of the right hand side.
Since $\|u^{D}_{\varepsilon,\delta;m}\|_{\infty;D} \leq \|u^{D}\|_{\infty;D},$ 
we have 
\[
\left|I_{0,1}(s, x; s', \varepsilon,\delta, m)\right| \leq \|u^{D}\|_{\infty;D}E_{s,x}\left[ \left|\exp Z_{s}(s')1_{\{s'<\hat{\sigma}_{\delta}\wedge T\wedge\tilde{\sigma}_{\varepsilon}\}} - 1\right|\right]
\]
for every $(s, x) \in D$ and $m \in \bm{N}.$ 
It holds that $P_{s,x}(\tilde{\sigma}_{\varepsilon} \uparrow  \sigma_{s}(\Pi) \mbox{ as } \varepsilon \downarrow 0) = 1$ by the continuity of the sample paths $X(\cdot)$ and 
$P_{s,x}(\sigma_{s}(\Pi) = \infty) = 1$ (see Appendix \ref{subsect.appendix.proof.hitting.time.1}). 
The continuity of the sample paths $X(\cdot)$ implies also $P_{s,x}(\hat{\sigma}_{\delta} \uparrow \sigma_{s}(\overline{\Sigma}_{(\infty)}) ~\mbox{ as } \delta \downarrow 0) = 1.$  
Therefore 
\begin{align*}
\lim_{s'\downarrow s}\lim_{\varepsilon\downarrow 0}\lim_{\delta\downarrow  0}\overline{\lim_{m\to\infty}}\left|I_{0,1}(s, x; s',\varepsilon,\delta, m)\right| = 0~~\mbox{ for } 
(s, x) \in D.
\end{align*}
Next examine the terms $I_{0,2}$ and $I_{0,3}$. 
Since $u^{D}_{\varepsilon,\delta;m} \to u^{D}_{\varepsilon,\delta}$ in $H^{1,2}(\Omega_{T})$ as $m \to \infty$, using the upper Gaussian bound (\ref{G.bound}) for the transition density, we have 
\begin{align*}
&\overline{\lim_{\delta\downarrow 0}}~ \overline{\lim_{m\to\infty}} I_{0,2}(s, x; s',\varepsilon,\delta, m)^{2} \\
&\leq \overline{\lim_{\delta\downarrow 0}}\,K(s' - s)^{-\frac{n}{2}}\int_{\Omega}\left\{u^{D}_{\varepsilon,\delta}(s', y) - u^{D}_{\varepsilon,\delta}(s, y)\right\}^{2}\exp\left\{- C\frac{|y - x|^{2}}{s' - s}\right\} dy  \\
&\leq K \int_{\mathbb R^{n}}\left\{u^{D}(s', x + (s'-s)^{1/2}z) - u^{D}(s, x + (s'-s)^{1/2}z)\right\}^{2}{\rm e}^{-C|z|^{2}} dz.
\end{align*}
Here $u^{D}(t, \cdot)$ is extended onto $\bm{R}^{n}$ with value zero outside $\Omega.$ 
Let 
\[
J_{0,2}(s, x; s'):= \overline{\lim_{\varepsilon\downarrow 0}}\, \overline{\lim_{\delta\downarrow 0}}\, \overline{\lim_{m\to \infty}} I_{0,2}(s, x; s', \varepsilon, \delta, m)^{2}.
\]
Noting that $u^{D} \in C([0, T]; L^{2}(\Omega)),$ we see that 
\begin{equation}
\lim_{s'\downarrow s} \int_{\Omega}J_{0,2}(s, x; s') dx = 0.
\label{J_{0,2}}
\end{equation}
For the term $I_{0,3}$,  
letting 
\[
J_{0,3} \equiv J_{0,3}(s, x; s'):= \overline{\lim_{\varepsilon\downarrow 0}}\, \overline{\lim_{\delta\downarrow 0}}\, \overline{\lim_{k\to \infty}} I_{0,3}(s, x; s', \varepsilon, \delta, m_{k})^{2}, 
\]
we get for each $s \in (0, T)$
\begin{equation}
\lim_{s'\downarrow s}\int_{\Omega}J_{0,3}(s, x; s') dx = 0
\label{J_{0,3}}
\end{equation}
as in the calculation on $I_{0,2}.$
Thus, combining (\ref{J_{0,2}}) with (\ref{J_{0,3}}), for each $s \in (0, T)$, we can select a subsequence $\{s_{\ell}\}_{\ell=1}^{\infty}$ with $s_{\ell} \downarrow s$ such that 
\[
\lim_{\ell\to \infty}J_{0,2}(s, x; s_{\ell}) = 0 ~~\mbox{ and } \lim_{\ell\to \infty}J_{0,3}(s, x; s_{\ell}) = 0
\]
for $x \in \Omega \setminus \mathcal{N}'(s)$ with $\mathcal{L}^{n}(\mathcal{N}'(s)) = 0.$ \\
We proceed to examine the term $I_{0,4}.$ 
For $v \in C([0, T]; L^{2}(\Omega)),$ let $\{v_{m}\}_{m=1}^{\infty}$ be its approximating sequence given by the smoothing procedure described above. 
Then 
\[
\sup_{m}\max_{0\leq t \leq T}\|v_{m}(t, \cdot)\|_{L^{2}(\Omega)} 
\leq B \max_{0\leq t \leq T}\|v(t, \cdot)\|_{L^{2}(\Omega)}
\]
for a constant $B$.  
In particular, if $\{u^{D}_{;m}\}_{m=1}^{\infty}$ means the approximating sequence for $u^{D}$,  
\begin{align}
\sup_{m}\max_{0\leq t \leq T}\|u^{D}_{;m}(t, \cdot) - 
u^{D}_{\varepsilon,\delta;m}(t, \cdot)\|_{L^{2}(\Omega)} \leq B \max_{0\leq t \leq T}\|u^{D}(t, \cdot) - u^{D}_{\varepsilon,\delta}(t, \cdot)\|_{L^{2}(\Omega)}. 
\label{approx.sol.2}
\end{align}
Using (\ref{approx.sol.1}),  (\ref{approx.sol.2}) and fact that $u^{D}_{\varepsilon,\delta;m} \longrightarrow u^{D}_{\varepsilon,\delta}$ in $H^{1,2}(\Omega_{T})$ as $m \to \infty$ (hence in $V^{0,1}(\Omega_{T})$ as $m \to \infty$), we ensure that 
\begin{equation*}
\max_{0\leq t \leq T}\|u^{D}_{;m}(t, \cdot) - u^{D}(t, \cdot)\|_{L^{2}(\Omega)} \longrightarrow 0 ~~\mbox{ as } ~m \to \infty.
\end{equation*}
This implies that 
\begin{align}
\lim_{m\to \infty}\max_{0\leq s \leq T}\mathcal{L}^{n}(\{x \in \Omega: \left|u^{D}_{;m}(s, x) - u^{D}(s, x)\right| > \eta\}) = 0 ~~\mbox{ for each } \eta > 0; 
\end{align}
accordingly, we can select a subsequence $\{m_{k}\}_{k=1}^{\infty}$ of $\{m\}_{m=1}^{\infty}$ so that 
\[
u^{D}_{;m_{k}}(s, x) \longrightarrow u^{D}(s, x) ~\mbox{ as } k \to \infty ~~
\]
for each $s \in (0, T)$ and $x \in \Omega \setminus \mathcal{N}(s)$ with $\mathcal{L}^{n}(\mathcal{N}(s)) = 0.$ 
Now we show that 
\begin{equation}
\lim_{\varepsilon\downarrow 0} \lim_{\delta\downarrow 0} \lim_{k\to \infty} \left|u^{D}_{\varepsilon, \delta, m_{k}}(s, x) - u^{D}(s, x)\right| = 0  
\label{a.e.conv}
\end{equation}
for each $s \in (0, T)$ and $x \in \Omega \setminus (\Sigma_{2}(s) \cup \mathcal{N}(s)).$ 
Let $\Upsilon:= \Sigma_{2(\infty)} \cup \Pi$ and for $\eta > 0$ set 
\[
\Upsilon^{(\eta)}:= \{(t, y) \in \Omega_{T}: \tilde{\rho}_{\Sigma_{2(\infty)}}(t, y) < \eta, ~
\tilde{\rho}_{\Pi}(y) < \eta\}.
\]
By the procedure of constructing $u^{D}_{\varepsilon,\delta;m}$ and $u^{D}_{;m},$ we see that 
for any $\eta > 0$ there is a positive integer $m_{0}$ such that 
\begin{align*}
u^{D}_{\varepsilon,\delta;m}(s, x) = u^{D}_{;m}(s, x) 
\end{align*}
for every $(s, x) \in \Omega_{T} \setminus \Upsilon^{(3\eta)}$, $m \geq m_{0}$ and $0 < \varepsilon,\delta < \eta$, 
because $u^{D}_{\varepsilon,\delta}(s, x) = u^{D}(s, x)$ for every $(s, x) \in \Omega_{T} \setminus \Upsilon^{(2\eta)}$ and $0 < \varepsilon,\delta < \eta.$
These imply the relationship (\ref{a.e.conv}). 
Instead of $\Sigma_{2}(s) \cup \mathcal{N}(s)$ we write it $\mathcal{N}(s)$ simply because  
$\mathcal{L}^{n}(\Sigma_{2}(s)) = 0.$ 
Therefore we obtain 
\begin{align}
\lim_{\varepsilon\downarrow 0} \lim_{\delta\downarrow 0}\lim_{k\to \infty}\left|I_{0,4}(s, x; \varepsilon, \delta, m_{k})\right| = 0~~\mbox{ for each }~s \in (0, T) \mbox{ and }~ x \in \Omega \setminus \mathcal{N}(s).
\end{align}
Thus we have checked the relationship (\ref{fact_0}). 

Next we will verify the equality (\ref{fact_3}). 
For a given $(s, x) \in D, $ without generality, we further take $\varepsilon$ and $\delta$ so small  that 
\begin{equation*}
(s, x) \in \left\{(t, y) \in D: \tilde{\rho}_{\Sigma_{(\infty)}}(t, y) > 4c^{\ast}\delta, ~ \tilde{\rho}_{\Pi}(y) > 2\varepsilon\right\}.
\end{equation*}
We divide $ I_{3}(s, x; s^{\prime}, \varepsilon, \delta, m)$ into the four terms: 
\begin{align*}
I_{3}(s, x; s^{\prime}, \varepsilon, \delta, m)  
&= E_{s,x}\left[\exp Z_{s}(s^{\prime})\,u^{D}_{\varepsilon,\delta;m}(s^{\prime}, X(s^{\prime})); 
\hat{\sigma}_{\delta}\wedge T< s^{\prime}, \hat{\sigma}_{\delta}\wedge T < \tilde{\sigma}_{\varepsilon}\right] \\ 
&+ E_{s,x}\left[\exp Z_{s}(s^{\prime})\,u^{D}_{\varepsilon,\delta;m}(s^{\prime}, X(s^{\prime})); 
\tilde{\sigma}_{\varepsilon} < s^{\prime}, \tilde{\sigma}_{\varepsilon} \leq \hat{\sigma}_{\delta}\wedge T\right] \\ 
&+ E_{s,x}\left[\exp Z_{s}(\hat{\sigma}_{\delta}\wedge T)\,u^{D}_{\varepsilon,\delta;m}(\hat{\sigma}_{\delta}\wedge T, X(\hat{\sigma}_{\delta}\wedge T)); s^{\prime} \leq \hat{\sigma}_{\delta}\wedge T < \tilde{\sigma}_{\varepsilon}\right] \\ 
&+ E_{s,x}\left[\exp Z_{s}(\tilde{\sigma}_{\varepsilon})\,u^{D}_{\varepsilon,\delta;m}(\tilde{\sigma}_{\varepsilon}, X(\tilde{\sigma}_{\varepsilon})); s^{\prime} \leq \tilde{\sigma}_{\varepsilon} \leq \hat{\sigma}_{\delta}\wedge T\right]
\end{align*}
and denote by $I_{3j}(s, x; s^{\prime}, \varepsilon, \delta, m)$ the $j$ th divided term of $I_{3}(s, x; s^{\prime}, \varepsilon, \delta, m)$ \\ $(j = 1, 2, 3, 4).$ 
Then it is easily seen that for $j = 1, 2, 4$
\[
\overline{\lim_{s'\downarrow s}}~\overline{\lim_{\varepsilon\downarrow 0}}~\overline{\lim_{\delta\downarrow 0}}~ \overline{\lim_{m\to\infty}}|I_{3j}(s, x; s^{\prime}, \varepsilon, \delta, m)| = 0.
\]
Furthermore divide $I_{33}(s, x; s^{\prime}, \varepsilon, \delta, m)$ as 
\[
I_{33}(s, x; s^{\prime}, \varepsilon, \delta, m) = I_{331}(s, x; s^{\prime}, \delta, \varepsilon, m) + I_{332}(s, x; \varepsilon, \delta, m),
\]
where for each $s' \in (s, T)$ 
\begin{align*}
I_{331}(s, x; s^{\prime}, \varepsilon, \delta, m) &:= E_{s,x}\left[\exp Z_{s}(\hat{\sigma}_{\delta})\,u^{D}_{\varepsilon,\delta;m}(\hat{\sigma}_{\delta}, X(\hat{\sigma}_{\delta})); s^{\prime} \leq \hat{\sigma}_{\delta} < \tilde{\sigma}_{\epsilon} \wedge T\right],  \\
I_{332}(s, x; \varepsilon, \delta, m) &:= E_{s,x}\left[\exp Z_{s}(T)\,u^{D}_{\varepsilon,\delta;m}(T, X(T)); T \leq \hat{\sigma}_{\delta}, ~T < \tilde{\sigma}_{\epsilon}\right].
\end{align*}
Note that 
\begin{align}
&I_{331}(s, x; s^{\prime}, \varepsilon, \delta, m)^{2} \nonumber \\
&\leq {\rm e}^{2\|c\|_{\infty}T}E_{s,x}\left[{\rm e}^{2\|\gamma\|_{\infty}L(T)}\right] 
E_{s,x}\left[\left(u^{D}_{\varepsilon,\delta;m}\right)^{2}(\hat{\sigma}_{\delta}, X(\hat{\sigma}_{\delta})); s^{\prime} \leq \hat{\sigma}_{\delta} < T \wedge \tilde{\sigma}_{\epsilon} \right].
\label{estim.331}
\end{align}
In the sequel, for a given $\varepsilon > 0,$ we write $\varepsilon^{\ast} := \varepsilon/c^{\ast}$, and further take $\delta > 0$ so small that 
\[
0 < \delta < \delta^{\ast}_{\varepsilon} : = \frac{1}{4c^{\ast}}
\left\{d(\Sigma_{(\infty)}, [0, \infty) \times \Gamma'^{[-\varepsilon^{\ast}/2]}) \wedge \frac{\varepsilon^{\ast}}{4}\right\}, 
\]
where $\Gamma'^{[-\eta]} := \Gamma' \setminus \Pi^{[\eta]}$ $(\eta > 0)$, and 
$\Pi^{[\eta]} = \{y \in \overline{\Omega}: \rho_{\Pi}(y) < \eta\}$. 
For such $\delta$, set 
\[
\bar{\mu}_{2c^{\ast}\delta}(t, y):= \tilde{\lambda}((2c^{\ast}\delta)^{-1}\tilde{\rho}_{\Sigma_{(\infty)}}(t, y)).
\]
Then $\bar{\mu}_{2c^{\ast}\delta}(t, y) = 1$ if $\tilde{\rho}_{\Sigma_{(\infty)}}(t, y) \geq 4c^{\ast}\delta$, so that the equality also holds on a neighborhood of $[0, \infty) \times \Gamma'^{(-\varepsilon)}$, where $\Gamma'^{(-\eta)} := \Gamma' \setminus \Pi^{\eta}$ 
$(\eta > 0)$.  
Therefore if we set 
\[
\tilde{u}^{D}_{\varepsilon,\delta;m}(t, y): = \left(1 - \bar{\mu}_{2c^{\ast}\delta}(t, y)\right)u^{D}_{\varepsilon,\delta;m}(t, y),
\]then  
\begin{align}
P_{s,x}\left(\tilde{u}^{D}_{\varepsilon,\delta;m}(\hat{\sigma}_{\delta}, X(\hat{\sigma}_{\delta})) = 
u^{D}_{\varepsilon,\delta;m}(\hat{\sigma}_{\delta}, X(\hat{\sigma}_{\delta})); \hat{\sigma}_{\delta} < \infty\right) = P_{s,x}\left(\hat{\sigma}_{\delta} < \infty\right), 
\label{331.key.1}
\end{align}
because 
\[
 P_{s,x}\left(\tilde{\rho}_{\Sigma_{(\infty)}}(\hat{\sigma}_{\delta}, X(\hat{\sigma}_{\delta})) \leq 2c^{\ast}\delta;  \hat{\sigma}_{\delta} < \infty\right) = P_{s,x}\left(\hat{\sigma}_{\delta} < \infty\right).
\]
In addition, for $\eta > 0,$ let 
\begin{align*}
\Sigma_{(\infty)}^{\eta}: = \left\{(t, y) \in (D \cup \partial_{L}D)_{(\infty)}: \tilde{\rho}_{\Sigma_{(\infty)}}(t, y) < \eta\right\} ~ ~\mbox{and} ~~ 
\check{\sigma}_{\delta}:= \sigma_{s}(\overline{\Sigma_{(\infty)}^{4c^{\ast}\delta}}).
\end{align*}
By definition, 
$P_{s,x}(\check{\sigma}_{\delta} \leq \hat{\sigma}_{\delta}) = 1.$ Moreover  
\begin{equation}
P_{s,x}\left(\tilde{u}^{D}_{\varepsilon,\delta;m}(\check{\sigma}_{\delta}, X(\check{\sigma}_{\delta})) = 0; 
\check{\sigma}_{\delta} < \infty\right) = P_{s,x}\left(\check{\sigma}_{\delta} < \infty\right),
\label{331.key.2}
\end{equation}
because 
$P_{s,x}\left(\tilde{\rho}_{\Sigma_{(\infty)}}(\check{\sigma}_{\delta}, X(\check{\sigma}_{\delta})) = 4c^{\ast}\delta; \check{\sigma}_{\delta} < \infty\right) = P_{s,x}\left(\check{\sigma}_{\delta} < \infty\right).$
Hence, by It\^{o}'s formula, it holds 
\begin{align}
&E_{s,x}\left[\left(\tilde{u}^{D}_{\varepsilon,\delta;m}\right)^{2}(\hat{\sigma}_{\delta}\wedge T\wedge\tilde{\sigma}_{\varepsilon}\vee s^{\prime}, X(\hat{\sigma}_{\delta}\wedge T\wedge\tilde{\sigma}_{\varepsilon}\vee s^{\prime}))\right]  \nonumber \\
&\hspace*{2.5cm}- E_{s,x}\left[\left(\tilde{u}^{D}_{\varepsilon,\delta;m}\right)^{2}(\check{\sigma}_{\delta}\wedge T\wedge\tilde{\sigma}_{\varepsilon}\vee s^{\prime}, X(\check{\sigma}_{\delta}\wedge T\wedge\tilde{\sigma}_{\varepsilon}\vee s^{\prime}))\right] \nonumber \\
&= E_{s,x}\left[\int_{\check{\sigma}_{\delta}\wedge T\wedge\tilde{\sigma}_{\varepsilon}\vee s^{\prime}}^{\hat{\sigma}_{\delta}\wedge T\wedge\tilde{\sigma}_{\varepsilon}\vee s^{\prime}}\mathcal{P}_{0}\left(\tilde{u}^{D}_{\varepsilon,\delta;m}\right)^{2}(t, X(t))dt\right] \nonumber \\
&\hspace*{2.5cm}+ E_{s,x}\left[\int_{\check{\sigma}_{\delta}\wedge T\wedge\tilde{\sigma}_{\varepsilon}\vee s^{\prime}}^{\hat{\sigma}_{\delta}\wedge T\wedge\tilde{\sigma}_{\varepsilon}\vee s^{\prime}}\mathcal{B}_{0}\left(\tilde{u}^{D}_{\varepsilon,\delta;m}\right)^{2}(t, X(t))L(dt)\right].
\label{Ito.formula.square}
\end{align}
Here, by the equality $\tilde{u}^{D}_{\varepsilon,\delta;m}(t, y) = 0$ on a neighborhood of 
$(\Gamma'^{(-\varepsilon)})_{T},$ we see that 
\begin{align*}
&E_{s,x}\left[\int_{\check{\sigma}_{\delta}\wedge T\wedge\tilde{\sigma}_{\varepsilon}\vee s^{\prime}}^{\hat{\sigma}_{\delta}\wedge T\wedge\tilde{\sigma}_{\varepsilon}\vee s^{\prime}}\mathcal{B}_{0}\left(\tilde{u}^{D}_{\varepsilon,\delta;m}\right)^{2}(t, X(t))L(dt)\right] \\
&= E_{s,x}\left[\int_{\check{\sigma}_{\delta}\wedge T\wedge\tilde{\sigma}_{\varepsilon}\vee s^{\prime}}^{(\hat{\sigma}_{\delta}\wedge T\wedge\tilde{\sigma}_{\varepsilon}\vee s^{\prime})-}\mathcal{B}_{0}\left(\tilde{u}^{D}_{\varepsilon,\delta;m}\right)^{2}(t, X(t))L(dt)\right] \\
&= E_{s,x}\left[\int_{\check{\sigma}_{\delta}\wedge T\wedge\tilde{\sigma}_{\varepsilon}\vee s^{\prime}}^{(\hat{\sigma}_{\delta}\wedge T\wedge\tilde{\sigma}_{\varepsilon}\vee s^{\prime})-}1_{(\Gamma'^{(-\varepsilon)})_{T}}(t, X(t)) \mathcal{B}_{0}\left(\tilde{u}^{D}_{\varepsilon,\delta;m}\right)^{2}(t, X(t))L(dt)\right] \\
&= 0.
\end{align*}

On the other hand, letting $\tilde{D}^{\ast}_{\varepsilon,\delta}:= \{(t, y) \in D: \tilde{\rho}_{\Sigma_{(\infty)}}(t, y) < 4c^{\ast}\delta, ~\tilde{\rho}_{\Pi}(y) > 2\varepsilon\}$ \\ $(\subset D \setminus (\overline{\Pi^{2\varepsilon}})_{T}),$ we have 
\begin{align*}
&\left|E_{s,x}\left[\int_{\check{\sigma}_{\delta}\wedge T\wedge\tilde{\sigma}_{\varepsilon}\vee s^{\prime}}^{\hat{\sigma}_{\delta}\wedge T\wedge\tilde{\sigma}_{\varepsilon}\vee s^{\prime}}\mathcal{P}_{0}\left(\tilde{u}^{D}_{\varepsilon,\delta;m}\right)^{2}(t, X(t))dt\right]        \right|  \\ 
&\leq E_{s,x}\left[\int_{\check{\sigma}_{\delta}\wedge T\wedge\tilde{\sigma}_{\varepsilon}\vee s^{\prime}}^{(\hat{\sigma}_{\delta}\wedge T\wedge\tilde{\sigma}_{\varepsilon}\vee s^{\prime})-}\left|\mathcal{P}_{0}\left(\tilde{u}^{D}_{\varepsilon,\delta;m}\right)^{2}(t, X(t))\right|dt\right]  \\
&= E_{s,x}\left[\int_{\check{\sigma}_{\delta}\wedge T\wedge\tilde{\sigma}_{\varepsilon}\vee s'}^{(\hat{\sigma}_{\delta}\wedge T\wedge\tilde{\sigma}_{\varepsilon}\vee s')-} 1_{\tilde{D}^{\ast}_{\varepsilon,\delta}}(t, X(t))\left|\mathcal{P}_{0}\left(\tilde{u}^{D}_{\varepsilon,\delta;m}\right)^{2}(t, X(t))\right|dt\right] \\
&\leq E_{s,x}\left[\int_{s'}^{T} 1_{\tilde{D}^{\ast}_{\varepsilon,\delta}}(t, X(t))\left|\mathcal{P}_{0}\left(\tilde{u}^{D}_{\varepsilon,\delta;m}\right)^{2}(t, X(t))\right|dt\right] \\
&= \int_{s'}^{T}dt \int_{\tilde{D}^{\ast}_{\varepsilon,\delta}(t)} \left|\mathcal{P}_{0}\left(\tilde{u}^{D}_{\varepsilon,\delta;m}\right)^{2}(t, y)\right| p(s, x; t, y) dy \\
&\leq K(s' - s)^{-\frac{n}{2}} \int_{s'}^{T}dt \int_{\tilde{D}^{\ast}_{\varepsilon,\delta}(t)} \left|\mathcal{P}_{0}\left(\tilde{u}^{D}_{\varepsilon,\delta;m}\right)^{2}(t, y)\right| dy. 
\end{align*}
Since, for each $\varepsilon$ and $\delta,$ $\tilde{u}^{D}_{\varepsilon,\delta;m}$ converges to 
$\tilde{u}^{D}_{\varepsilon,\delta}:= (1 - \bar{\mu}_{2c^{\ast}\delta})u^{D}_{\varepsilon,\delta}$ in $H^{1,2}(\Omega_{T})$ as $m \to \infty$ (which is derived from the property (\ref{reg.sol}) for $u^{D}_{\varepsilon,\delta}$), the last term of the right hand side of the inequality just above converges to 
\begin{equation}
K(s' - s)^{-\frac{n}{2}} \int_{s'}^{T}dt \int_{\tilde{D}^{\ast}_{\varepsilon,\delta}(t)} \left|\mathcal{P}_{0}\left(\tilde{u}^{D}_{\varepsilon,\delta}\right)^{2}(t, y)\right| dy 
\label{pre-integral}
\end{equation}
as $m \to \infty.$

In the following, we will show that the integral of (\ref{pre-integral}) converges to zero as $\delta \downarrow 0$ for a given $\varepsilon$. 
Let $\nu_{\delta}:= (1 - \bar{\mu}_{2c^{\ast}\delta})\tilde{\mu}_{\delta}.$ 
Then $\mbox{spt}\, \nu_{\delta} \cap D \subset \{(t, y) \in D: c^{\ast}\delta \leq \tilde{\rho}_{\Sigma_{(\infty)}}(t, y) \leq 4c^{\ast}\delta\}$ and 
$\tilde{u}^{D}_{\varepsilon,\delta} =\nu_{\delta}u_{\varepsilon \, \centerdot}^{D}$ on $D$. 
Therefore 

\begin{align}
\mathcal{P}_{0}\left(\tilde{u}^{D}_{\varepsilon,\delta}\right)^{2} &= 2\nu_{\delta}^{2}u^{D}_{\varepsilon \, \centerdot}\left\{\frac{\partial u^{D}_{\varepsilon \, \centerdot}}{\partial t} + \mbox{Tr}(A\nabla^{2}u^{D}_{\varepsilon \, \centerdot}) + \bm{c}\iprod\nabla u^{D}_{\varepsilon \, \centerdot}\right\} \nonumber  \\ 
&+ 2\Biggl\{\nu_{\delta}(u^{D}_{\varepsilon \, \centerdot})^{2}\frac{\partial\nu_{\delta}}{\partial t} + 
(u^{D}_{\varepsilon \, \centerdot})^{2}(A\nabla\nu_{\delta})\iprod\nabla\nu_{\delta} + \nu_{\delta}(u^{D}_{\varepsilon \, \centerdot})^{2}\mbox{Tr}(A\nabla^{2}\nu_{\delta})   \nonumber \\ 
&+ 4\nu_{\delta}u^{D}_{\varepsilon \, \centerdot}(A\nabla\nu_{\delta})\iprod\nabla u^{D}_{\varepsilon \, \centerdot} + \nu_{\delta}^{2}(A\nabla u^{D}_{\varepsilon \, \centerdot})\iprod\nabla u^{D}_{\varepsilon \, \centerdot} + 
\nu_{\delta}(u^{D}_{\varepsilon \, \centerdot})^{2}\bm{c}\iprod\nabla\nu_{\delta}\Biggr\}
\label{strong_form}
\end{align}
on $D$  in the strong form. 
Letting $f_{\varepsilon}$ and $g_{\varepsilon}$ the first and second terms of the right hand side of (\ref{strong_form}), we note 
\[
f_{\varepsilon} = 2\nu_{\delta}^{2}u^{D}(f - cu^{D})~~\mbox{ on } D \setminus (\overline{\Pi^{2\varepsilon}})_{T},
\] since $u^{D}_{\varepsilon \, \centerdot} = u^{D}$ on $D \setminus (\overline{\Pi^{2\varepsilon}})_{T}$. 
Here, for variables $\xi$ and $\eta$, the notations $\xi \lessapprox \eta$ and $\xi \approx \eta$ designate the relationships $\xi \leq c\eta$ and  $c_{1}\eta \leq \xi \leq c_{2}\eta$ with some positive constants $c$ and $c_{1}, ~ c_{2}$, respectively. 
Now we take a number in $(0, \delta^{\ast}_{\varepsilon})$, say $\delta_{0}$, and put $D^{\ast}_{\delta} := \{(t, y) \in D: \rho_{\Sigma_{(\infty)}}(t, y) < 4c^{\ast} \delta\}$ for $0 < \delta \leq \delta_{0}$.
Note that 
\begin{align}
&0 \leq \nu_{\delta} \leq 1,~~ \nu_{\delta} \longrightarrow 0 ~\mbox{ as } \delta \downarrow 0, 
\label{fact-1} \\
&\left|\frac{\partial\nu_{\delta}}{\partial t}(t, y)\right| \approx 
\left|\frac{\partial\nu_{\delta}}{\partial y_{i}}(t, y)\right| \approx \frac{1}{\delta} \approx \frac{1}{\tilde{\rho}_{\Sigma_{(\infty)}}(t, y)} ~~\mbox{ on }~ \mathrm{spt}\, \nu_{\delta} \cap D, \label{fact-2} \\
&\left|\frac{\partial^{2}\nu_{\delta}}{\partial y_{i}\partial y_{j}}(t, y)\right| \approx  
\frac{1}{\delta^{2}} \approx \frac{1}{\tilde{\rho}_{\Sigma_{(\infty)}}(t, y)^{2}} ~~\mbox{ on }~ \mathrm{spt}\, \nu_{\delta} \cap D, \label{fact-3} \\
&\rho_{\Sigma_{(\infty)}}(t, y) \approx \rho_{\Sigma_{(\infty)}^{[-\varepsilon^{\ast}/2]}}(t, y)  \quad \mbox{on} ~ D^{\ast}_{\delta_{0}} \setminus \bigl(\Pi^{[\varepsilon^{\ast}]}\bigr)_{T}, 
\label{fact-4}
\end{align}
where $\Sigma^{[-\varepsilon^{\ast}/2]} := \Sigma \setminus \bigl(\Pi^{[\varepsilon^{\ast}/2]}\bigr)_{T}$. 
Therefore, by (\ref{fact-1}), for $s' \in (s, T)$ and sufficiently small $\varepsilon > 0,$ 
\[
\int_{s^{\prime}}^{T}dt \int_{\tilde{D}^{\ast}_{\varepsilon,\delta}(t)}\left|f_{\varepsilon}(t, y)\right| dy \longrightarrow 0  ~\mbox{ as } \delta \downarrow 0. 
\]
Next examine the integral for $|g_{\varepsilon}(t, y)|$. 
Noting the facts (\ref{fact-3}), (\ref{fact-4}) and $u^{D}_{\varepsilon \, \centerdot} = 0$ on $\left(\Pi^{\varepsilon}\right)_{T}$ which includes $\bigl(\Pi^{[\varepsilon^{\ast}]}\bigr)_{T}$, and then applying Lemma 4.1 in \cite{Kaw10} by replacing $G$, $\Xi$ and $\Xi'$ with $D^{\ast}_{\delta_{0}}$, $\Sigma$ and $\Sigma^{[-\varepsilon^{\ast}/2]}$, respectively, we have 
$\displaystyle{\frac{u^{D}_{\varepsilon \, \centerdot}(t, y)}{d((t, y), \Sigma^{[-\varepsilon^{\ast}/2]})} \in L^{2}(D^{\ast}_{\delta_{0}})}$; hence for $0 < \delta < \delta_{0}$  
\begin{align*}
&\int_{s^{\prime}}^{T}dt \int_{\tilde{D}^{\ast}_{\varepsilon,\delta}(t)}\left|\frac{\partial^{2}\nu_{\delta}}{\partial y_{i}\partial y_{j}}(t, y)\right|(u^{D}_{\varepsilon \, \centerdot})^{2}(t, y) dy \\
&\leq \int_{s^{\prime}}^{T}dt \int_{\tilde{D}^{\ast}_{\delta}(t)}\left|\frac{\partial^{2}\nu_{\delta}}{\partial y_{i}\partial y_{j}}(t, y)\right|(u^{D}_{\varepsilon \, \centerdot})^{2}(t, y) dy \\
&\lessapprox \int_{s^{\prime}}^{T}dt \int_{D^{\ast}_{\delta}(t)} \frac{u^{D}_{\varepsilon \, \centerdot}(t, y)^{2}}{d((t, y), \Sigma^{[-\varepsilon^{\ast}/2]})^{2}} dy 
\longrightarrow 0  ~\mbox{ as } \delta \downarrow 0.
\end{align*}
Similar integrals for the other terms consisting of the function $g_{\varepsilon}$ also converge to zero as $\delta \downarrow 0;$ hence for $s' \in (s, T)$ and sufficiently small $\varepsilon > 0$ 
\[
\int_{s^{\prime}}^{T}dt \int_{\tilde{D}^{\ast}_{\varepsilon,\delta}(t)}|g_{\varepsilon}(t, y)| dy 
\longrightarrow 0  ~\mbox{ as } \delta \downarrow 0. 
\]
These ensure that 
\[
E_{s,x}\left[\int_{s^{\prime}}^{T}1_{\tilde{D}^{\ast}_{\varepsilon,\delta}}(t, X(t)) \left|\mathcal{P}_{0}\left(\tilde{u}^{D}_{\varepsilon,\delta}\right)^{2}(t, X(t))\right|dt\right] \longrightarrow 0 ~\mbox{ as } \delta \downarrow 0;
\]
as a result, for each $(s, x) \in D,$ $s' \in (s, T)$ and sufficiently small $\varepsilon > 0,$ 
\[
\lim_{\delta\downarrow 0}\overline{\lim_{m\to \infty}}\left|E_{s,x}\left[\int_{\check{\sigma}_{\delta}\wedge T\wedge\tilde{\sigma}_{\varepsilon}\vee s^{\prime}}^{\hat{\sigma}_{\delta}\wedge T\wedge\tilde{\sigma}_{\varepsilon}\vee s^{\prime}}\mathcal{P}_{0}\left(\tilde{u}^{D}_{\varepsilon,\delta;m}\right)^{2}(t, X(t))dt\right]\right| = 0.
\]
That is, 
\begin{align}
&\lim_{\delta\downarrow 0}\overline{\lim_{m\to \infty}}\left|E_{s,x}\left[\left(\tilde{u}^{D}_{\varepsilon,\delta;m}\right)^{2}(\hat{\sigma}_{\delta}\wedge T\wedge\tilde{\sigma}_{\varepsilon}\vee s^{\prime}, X(\hat{\sigma}_{\delta}\wedge T\wedge\tilde{\sigma}_{\varepsilon}\vee s^{\prime}))\right] \right. \nonumber \\
&- \left.E_{s,x}\left[\left(\tilde{u}^{D}_{\varepsilon,\delta;m}\right)^{2}(\check{\sigma}_{\delta}\wedge T\wedge\tilde{\sigma}_{\varepsilon}\vee s^{\prime}, X(\check{\sigma}_{\delta}\wedge T\wedge\tilde{\sigma}_{\varepsilon}\vee s^{\prime}))\right]\right| = 0  
\label{Ito.formula.left}
\end{align}
for each $(s, x) \in D,$ $s' \in (s, T)$ and sufficiently small $\varepsilon > 0.$ 
Now we divide the first and second terms of the left hand side of (\ref{Ito.formula.square}) as follows: 
\begin{align}
&E_{s,x}\left[\left(\tilde{u}^{D}_{\varepsilon,\delta;m}\right)^{2}(\hat{\sigma}_{\delta}\wedge T\wedge\tilde{\sigma}_{\varepsilon}\vee s^{\prime}, X(\hat{\sigma}_{\delta}\wedge T\wedge\tilde{\sigma}_{\varepsilon}\vee s^{\prime}))\right] \nonumber \\
&= E_{s,x}\left[\left(\tilde{u}^{D}_{\varepsilon,\delta;m}\right)^{2}(\hat{\sigma}_{\delta}, X(\hat{\sigma}_{\delta})); s' \leq \hat{\sigma}_{\delta} < T \wedge \tilde{\sigma}_{\varepsilon}         \right] \nonumber \\ 
&+ E_{s,x}\left[\left(\tilde{u}^{D}_{\varepsilon,\delta;m}\right)^{2}(T, X(T)); T \leq \hat{\sigma}_{\delta} \wedge \tilde{\sigma}_{\varepsilon}\right] \nonumber \\ 
&+ E_{s,x}\left[\left(\tilde{u}^{D}_{\varepsilon,\delta;m}\right)^{2}(\tilde{\sigma}_{\varepsilon}, X(\tilde{\sigma}_{\varepsilon})); s' \leq  \tilde{\sigma}_{\varepsilon} \leq T \wedge \hat{\sigma}_{\delta}\right] \nonumber \\ 
&+ E_{s,x}\left[\left(\tilde{u}^{D}_{\varepsilon,\delta;m}\right)^{2}(s', X(s')); \hat{\sigma}_{\delta} \wedge \tilde{\sigma}_{\varepsilon} < s'\right]
\label{decomp.1}
\end{align}
and 
\begin{align}
&E_{s,x}\left[\left(\tilde{u}^{D}_{\varepsilon,\delta;m}\right)^{2}(\check{\sigma}_{\delta}\wedge T\wedge\tilde{\sigma}_{\varepsilon}\vee s^{\prime}, X(\check{\sigma}_{\delta}\wedge T\wedge\tilde{\sigma}_{\varepsilon}\vee s^{\prime}))\right] \nonumber \\
&= E_{s,x}\left[\left(\tilde{u}^{D}_{\varepsilon,\delta;m}\right)^{2}(\check{\sigma}_{\delta}, X(\check{\sigma}_{\delta})); s' \leq \check{\sigma}_{\delta} < T \wedge \tilde{\sigma}_{\varepsilon}\right] \nonumber \\ 
&+ E_{s,x}\left[\left(\tilde{u}^{D}_{\varepsilon,\delta;m}\right)^{2}(T, X(T)); T \leq \check{\sigma}_{\delta} \wedge \tilde{\sigma}_{\varepsilon}\right] \nonumber \\ 
&+ E_{s,x}\left[\left(\tilde{u}^{D}_{\varepsilon,\delta;m}\right)^{2}(\tilde{\sigma}_{\varepsilon}, X(\tilde{\sigma}_{\varepsilon})); s' \leq  \tilde{\sigma}_{\varepsilon} \leq T \wedge \check{\sigma}_{\delta}\right] \nonumber \\ 
&+ E_{s,x}\left[\left(\tilde{u}^{D}_{\varepsilon,\delta;m}\right)^{2}(s', X(s')); \check{\sigma}_{\delta} \wedge \tilde{\sigma}_{\varepsilon} < s'\right]. 
\label{decomp.2}
\end{align}
Here we indicate by $\tilde{I}$ and $\tilde{I}_{k}$ the left hand side and the $k$ th term of the right hand side of (\ref{decomp.1}) ($k = 1, 2, 3, 4$)  , and by $\tilde{J}$ and $\tilde{J}_{k}$ the left hand side and the $k$ th term of the right hand side of (\ref{decomp.2}) ($k = 1, 2, 3, 4$), respectively. 
Then by (\ref{331.key.1}) and (\ref{331.key.2}), 
\begin{align}
&\tilde{I}_{1} = E_{s,x}\left[\left(u^{D}_{\varepsilon,\delta;m}\right)^{2}(\hat{\sigma}_{\delta}, X(\hat{\sigma}_{\delta})); s' \leq \hat{\sigma}_{\delta} < T \wedge \tilde{\sigma}_{\varepsilon}\right],  \nonumber \\
&\tilde{J}_{1} = 0.
\label{key-equality}
\end{align}
These equalities in (\ref{key-equality}) imply
\begin{align}
0 &\leq E_{s,x}\left[\left(u^{D}_{\varepsilon,\delta;m}\right)^{2}(\hat{\sigma}_{\delta}, X(\hat{\sigma}_{\delta})); s' \leq \hat{\sigma}_{\delta} < T \wedge \tilde{\sigma}_{\varepsilon}\right] \nonumber \\ 
&= (\tilde{I} - \tilde{J}) - (\tilde{I}_{2} - \tilde{J}_{2}) - (\tilde{I}_{3} - \tilde{J}_{3}) - \tilde{I}_{4} + \tilde{J}_{4} \nonumber \\ 
&\leq |\tilde{I} - \tilde{J}| + \tilde{J}_{4}, 
\label{decomp.3}
\end{align}
because $\tilde{I}_{k} - \tilde{J}_{k} \geq 0$ $(k = 2, 3)$ and $\tilde{I}_{4} \geq 0.$ 
By (\ref{Ito.formula.left}), we have   
\begin{align*}
\lim_{\delta\downarrow 0}\overline{\lim_{m\to\infty}} |\tilde{I} - \tilde{J}| = 0. 
\end{align*}
Moreover, the inequality  $0 \leq \tilde{J}_{4} \leq \|u^{D}\|_{\infty;D}^{2}\left\{P_{s,x}(\check{\sigma}_{\delta} < s') +  P_{s,x}(\tilde{\sigma}_{\varepsilon} < s')\right\}$ yields 
\[
\lim_{s'\downarrow s}\, \overline{\lim_{\varepsilon \downarrow 0}} \,\overline{\lim_{\delta\downarrow 0}}\,\overline{\lim_{m\to \infty}}\tilde{J}_{4} = 0. 
\]
Hence, by (\ref{estim.331}) and (\ref{decomp.3}), it holds  
\[
\lim_{s'\downarrow s}\, \overline{\lim_{\varepsilon \downarrow 0}} \,
\overline{\lim_{\delta\downarrow 0}}\,\overline{\lim_{m\to \infty}}|I_{331}(s, x; s', \varepsilon, \delta, m)|^{2} = 0
\]
for each $(s, x) \in D.$ 
Noting the facts that $P_{s,x}(\tilde{\sigma}_{\varepsilon} \uparrow \infty ~\mbox{ as } \varepsilon\downarrow 0) = 1$ and $P_{s,x}(\sigma_{s}(\overline{\Sigma}_{(\infty)}) = T) = 0$ (which is derived from the facts $\{\sigma_{s}(\overline{\Sigma}_{(\infty)}) = T\} \subset \{X(T) \in \overline{\Sigma}(T)\}$ and $\mathcal{L}^{n}(\overline{\Sigma}(T)) = 0$), we also have 
\begin{align*}
\lim_{\varepsilon \downarrow 0} \,\overline{\lim_{\delta\downarrow 0}}\,\overline{\lim_{m\to \infty}}\left|I_{332}(s, x; \delta, \varepsilon, m) - 
E_{s,x}\left[\exp Z_{s}(T) h(X(T)); \sigma_{s}(\overline{\Sigma}_{(\infty)}) > T\right] \right| = 0
\end{align*}
for each $(s, x) \in D.$ 
Therefore  (\ref{fact_3}) is verified. 

The equalities (\ref{fact_1}) and (\ref{fact_2}) are verified more easily by using the upper Gaussian bound of the transition density and the regularity of the weak solution.  
Thus the proof in the case (i) is complete.

Next we proceed to verify the equality (\ref{repre}) in the case (ii). 
In the first stage, we show that the right hand side of (\ref{repre}), say $v^{D}$, is continuous up to the Robin part. 
Note that for $s \in [0, T]$ and $x \in D(s) \cup \Gamma^{\prime}$ 
\begin{align}
P_{s,x}(\sigma_{s}(\overline{\Sigma}_{(\infty)}) = \sigma_{0}(\overline{\Sigma}_{(\infty)})) = 1, ~~~P_{s,x}(\sigma_{s}(\delta\Sigma_{1(\infty)}) = \sigma_{0}(\delta\Sigma_{1(\infty)})) = 1;
\label{s.to.0}
\end{align}
hence 
\begin{align}
 v^{D}(s, x) = &-\exp\left(- \int_{0}^{s}c(r, x)dr\right)\Biggl\{E_{s,x}\left[\int_{0}^{\sigma_{0}(\overline{\Sigma}_{(\infty)})\wedge T}\exp Z_{0}(t)f(t, X(t))dt\right] \nonumber\\
&\hspace*{6cm}- \int_{0}^{s}\exp\left(- \int_{0}^{t}c(r, x)dr\right)\, f(t, x)dt\Biggr\} \nonumber\\
&- \exp\left(- \int_{0}^{s}c(r, x)dr\right)E_{s,x}\left[\int_{0}^{\sigma_{0}(\overline{\Sigma}_{(\infty)})\wedge T}\exp Z_{0}(t)\,\psi(t, X(t))L(dt)\right] \nonumber\\
&+ \exp\left(- \int_{0}^{s}c(r, x)dr\right)E_{s,x}\left[\exp Z_{0}(T)\, h(X(T))\, 1_{(T,\infty)}
(\sigma_{0}(\overline{\Sigma}_{(\infty)}))\right].
\label{s_to_0}
\end{align}
For each $m \in \bm{N},$ take an appropriate cut-off function $h_{m}$ for the indicator function $1_{[0,\infty)}$ (e.g. the cut-off function $h_{m}$ on page 275 in \cite{Tsu94}).  Then \\
$\displaystyle{h_{m}(L(T))\int_{0}^{t\wedge T}g(r, X(r))L(dr)}$ is continuous in $(t, \bm{\omega}) \in [0, \infty) \times \bm{U}$ for every bounded continuous function $g$ on $[0, \infty) \times \overline{\Omega}.$
On the other hand, in Appendix \ref{subsect.continuity.hitting.time} (see Theorem \ref{continuity.final}), we will show that there exists a measurable subset 
$\widetilde{\bm{U}}$ of $\bm{U}$ such that $P_{s,x}(\widetilde{\bm{U}}) = 1$ for every $s \in [0, T], ~ x \in D(s) \cup \Gamma^{\prime}$ and the restriction of $\sigma_{0}(\overline{\Sigma}_{(\infty)})$ into $\widetilde{\bm{U}}$, say  $\left.\sigma_{0}(\overline{\Sigma}_{(\infty)})\right|_{\widetilde{\bm{U}}},$ is continuous. 
Therefore the function 
\[
h_{m}(L(T))\int_{0}^{\sigma_{0}(\overline{\Sigma}_{(\infty)})\wedge T}g(r, X(r))L(dr)
\]
of $\bm{\omega}$ is bounded on $\bm{U}$ and its restriction into $\widetilde{\bm{U}}$ is continuous. 
Combining the continuity property of $\{P_{s,x}\}$ on weak convergence and with a mapping theorem in weak convergence given in Appendix \ref{subsect.mapping.theorem} (see Theorem \ref{mapping.theorem}), we ensure that 
\[
E_{s,x}\left[h_{m}(L(T))\int_{0}^{\sigma_{0}(\overline{\Sigma}_{(\infty)})\wedge T}g(r, X(r))L(dr)\right]
\]
is continuous in $(s, x)$ with $s \in [0, T], x \in D(s) \cup \Gamma^{\prime}.$ 
Noting (\ref{unif.expo}), the uniform exponential integrability of $L(T)$, we see that $v^{D}$ is continuous in $(s, x)$ with $s \in [0, T], x \in D(s) \cup \Gamma^{\prime},$ and also the first term of the right hand side of (\ref{s_to_0}) is continuous in $(s, x)$ with $s \in [0, T], x \in D(s) \cup \Gamma^{\prime}.$  
By virtue of the fact $P_{s,x}(\sigma_{0}(\overline{\Sigma}_{(\infty)}) = T) = 0$ with the semicontinuity of functions $1_{(T,\infty)}(\xi)$ and $1_{[T,\infty)}(\xi)$, it follows that $1_{(T,\infty)}(\sigma_{0}(\overline{\Sigma}_{(\infty)}))$ is continuous on $\widetilde{\bm{U}} \cap \{\sigma_{0}(\overline{\Sigma}_{(\infty)}) \neq T\}.$ 
Hence the third term of the right hand side of (\ref{s_to_0}) is continuous in $(s, x)$ with $s \in [0, T], x \in D(s) \cup \Gamma^{\prime}.$  
As a result, the trace of $v^{D}$ into the set $\Gamma^{\prime}_{T}$ coincides with its ordinary boundary value almost everywhere with respect to the measure $ds \times S(dx).$ 

Finally we examine the continuity of $v^{D}$ on the Dirichlet part $\Sigma$. Since $v^{D} = 0$ on 
$\Sigma$, it is enough to see that for every $(s_{0}, x_{0}) \in \Sigma$ 
\begin{equation}
\lim_{D \ni (s, x) \to (s_{0}, x_{0})}v^{D}(s, x) = 0. 
\label{cont.stochas.sol}
\end{equation}
In the case of $(s_{0}, x_{0}) \in \Sigma_{2}$, by using the same way as in the proof of Theorem 13.1 in \cite{Wen81}, we have for any $\eta > 0$ 
\begin{equation}
P_{s,x}(\sigma_{s}(\overline{\Sigma}_{(\infty)}) \geq s + \eta) \rightarrow 0 ~~\mbox{as}~(s, x) \to (s_{0}, x_{0}),
\label{cont.Dirichlet_{2}}
\end{equation}
because 
$
P_{s,x}\left(\sigma_{s}(\overline{\Sigma}_{2(\infty)}) = \sigma_{s+}(((0, \infty) \times \Omega)\setminus \overline{D_{(\infty)}})\right) =1
$ 
for $s \geq 0$ and $x \in D(s) \cup \Gamma^{\prime}$;  
it is shown in (2--2) of (2) in Appendix \ref{subsect.continuity.hitting.time}. 
The relationship (\ref{cont.Dirichlet_{2}}) implies the continuity property (\ref{cont.stochas.sol}) of $v^{D}$ at $(s_{0}, x_{0}).$ The continuity property is verified as in the more complicated case of $(s_{0}, x_{0}) \in \Sigma_{1}$; so it is omitted. 
Next consider the case of $(s_{0}, x_{0}) \in \Sigma_{1}.$ 
Take an arbitrarily fixed positive number $\eta.$ 
Then we see that  
for any $\epsilon > 0$ there exists $\delta \in (0, \eta)$ such that 
\begin{equation} 
P_{s,x}(\sigma_{s}(\overline{\Sigma}_{(\infty)}) \geq s + 2\eta) < 2\epsilon + P_{s_{0},x_{0}}(\sigma_{s_{0}}(\overline{\Gamma'}) < s_{0} + 4\eta)
\label{cont.Dirichlet_{1}}
\end{equation} 
for every $(s, x) \in D$ with $d((s, x), (s_{0},x_{0})) < \delta;$ which will be verified in Appendix \ref{subsect.continuity.Dirichlet.part_{1}}.
Then the continuity property (\ref{cont.stochas.sol}) is proved as follows. 
From now on, suppose that $(s, x) \in D$ with $d((s, x), (s_{0},x_{0})) < \delta.$
We examine the second term of $v^{D}:$ 

\begin{align*}
&\left|E_{s,x}\left[\int_{s}^{\sigma_{s}(\overline{\Sigma}_{(\infty)})\wedge T}\exp Z_{s}(t) \psi(t, X(t)) L(dt)\right]\right| \\
&\leq \|\psi\|_{\infty} e^{\|c\|_{\infty}T}\left\{E_{s,x}\left[\int_{s}^{\sigma_{s}(\overline{\Sigma}_{(\infty)})\wedge T}e^{\|\gamma\|_{\infty}L(T)} L(dt); \sigma_{s}(\overline{\Sigma}_{(\infty)})\wedge T < s + 2\eta\right] \right. \\
&\hspace*{1.5cm}+ \left.E_{s,x}\left[\int_{s}^{\sigma_{s}(\overline{\Sigma}_{(\infty)})\wedge T}e^{\|\gamma\|_{\infty}L(T)} L(dt); \sigma_{s}(\overline{\Sigma}_{(\infty)})\wedge T \geq s + 2\eta\right]\right\}. 
\end{align*} 
Denote by I and II the first and second terms in the braces of the right hand side of the inequality just above, respectively. Then we have 
\begin{align*}
\text{I} &\leq E_{s,x}\left[e^{\|\gamma\|_{\infty}L(T)}L(s + 2\eta)\right]  \\
&\leq \sqrt{\sup_{(s,x) \in D}E_{s,x}\left[e^{2(\|\gamma\|_{\infty}+1)L(T)}\right]} \sqrt{E_{s,x}\left[e^{- 2L(T)}L(s_{0} + 3\eta)^{2}\right]}.  
\end{align*} 
Since 
\[
E_{s,x}\left[e^{- 2L(T)}L(s_{0} + 3\eta)^{2}\right] \rightarrow 
E_{s_{0},x_{0}}\left[e^{- 2L(T)}L(s_{0} + 3\eta)^{2}\right]
\]
as $(s, x) \to (s_{0},x_{0})$ and $E_{s_{0},x_{0}}\left[e^{- 2L(T)}L(s_{0} + 3\eta)^{2}\right] \to 0$ as $\eta \to 0,$ we have $\text{I} \to 0$ as $(s, x) \to (s_{0},x_{0}).$
On the other hand,
\begin{align*}
\text{II} &\leq  E_{s,x}\left[e^{\|\gamma\|_{\infty}L(T)}L(T); \sigma_{s}(\overline{\Sigma}_{(\infty)}) \geq s + 2\eta \right] \\
&\leq e^{-1} \sqrt{\sup_{(s,x) \in D}E_{s,x}\left[e^{2(\|\gamma\|_{\infty}+1)L(T)}\right]}\sqrt{2\epsilon + P_{s_{0},x_{0}}\left(\sigma_{s_{0}}(\overline{\Gamma'}) < s_{0} + 4\eta \right)}.  
\end{align*} 
Therefore, letting $\epsilon \to 0$ and then $\eta \to 0,$ we see that 
\[
\displaystyle{\limsup_{D \ni (s,x) \to (s_{0},x_{0})} \text{II} = 0},
\]
since 
$P_{s_{0},x_{0}}\left(\sigma_{s_{0}}(\overline{\Gamma'}) > s_{0}\right) = 1.$ 
The continuity property of the other terms of $v^{D}$ is checked easily. 
Consequently, the last assertion of the theorem is verified. 
That is, the weak solution $u^{D}$ has a version which is continuous up to the lateral boundary 
$\partial_{L}D$ except 
$[0, T] \times \Pi$, which is given by the stochastic solution $v^{D}$. 
\hspace*{\fill}$\Box$

\medskip 

\medskip 

Furthermore, the continuity of the source term $f$ and the Lipschitz continuity of the source term $\psi$ and the terminal value $h$ in (iv) of Assumption \ref{assump_3} are relaxed as follows.
{\corollary
\label{main_result_cor}
Replace the conditions on $f, \psi, h$ in (iv) of Assumption \ref{assump_3} with the following ones:
\begin{enumerate}
\item[(a)]  $\displaystyle{f \in L^{p_{1}}(\Omega_{T})~~\mbox{for some}~ p_{1} > \frac{n + 2}{2}}$, 
\item[(b)]  $\displaystyle{\psi \in L^{p_{2}}(\Gamma'_{T})~~\mbox{for some}~ p_{2} > n + 1}$,
\item[(c)]  $\displaystyle{h \in L^{p_{3}}(\Omega)~~\mbox{for some}~ p_{3} \geq 2}$.
\end{enumerate} 
Then the conclusion of Theorem \ref{main_result} holds.
}
\medskip

\noindent
\textbf{Proof.}  
For given $p_{1}, p_{2}, p_{3}$, choose $q_{1}, q_{2}, q_{3}$ and $r_{1}, r_{2}, r_{3}$ in such a way: 
\begin{align*}
p_{1} > q_{1} > \frac{n + 2}{2},~ r_{1} = \frac{p_{1}}{q_{1}}; ~p_{2} > q_{2} > n + 1,~ r_{2} = \frac{p_{2}}{q_{2}}; ~p_{3} > q_{3} > 1,~ r_{3} = \frac{p_{3}}{q_{3}}.
\end{align*}
First assume that $f, \psi, h$ satisfy the condition (iv) of Assumption \ref{assump_3}. 
Then, using the upper Gaussian bound (\ref{G.bound}) and the estimate (1.2) in Lemma 7.1.1 in \cite{Str79} and noting Theorems 4.1 and 4.3 in \cite{Kaw00}, we have the following: 
there exist positive constants $\Lambda_{i} = \Lambda_{i}(K, C, n, q_{i}, r_{i}, T)$ $(i = 1, 2)$ and $\Lambda_{3} = \Lambda_{3}(K, C, n, q_{3}, r_{3})$ such that 
\begin{align*}
&(a') ~~ \sup_{(s,x) \in D \cup \partial_{L}D}\left|E_{s,x}\left[\int_{s}^{\sigma_{s}(\overline{\Sigma}_{(\infty)})\wedge T}\exp Z_{s}(t) f(t, X(t)) dt\right]\right| \\
&\hspace*{2cm}\leq \Lambda_{1}\sup_{(s,x) \in \overline{\Omega_{T}}} E_{s,x}\left[\exp\{r_{1}^{\ast}\|\gamma\|_{\infty}L(T)\}\right]^{1/r_{1}^{\ast}}\|f\|_{p_{1};\Omega_{T}}; \\
&(b') ~~ \sup_{(s,x) \in D \cup \partial_{L}D}\left|E_{s,x}\left[\int_{s}^{\sigma_{s}(\overline{\Sigma}_{(\infty)})\wedge T}\exp Z_{s}(t) \psi(t, X(t)) L(dt)\right]\right| \\
&\hspace*{2cm}\leq \Lambda_{2}\sup_{(s,x) \in \overline{\Omega_{T}}} E_{s,x}\left[\exp\{(r_{2}^{\ast}\|\gamma\|_{\infty} + 1)L(T)\}\right]^{1/r_{2}^{\ast}}\|\psi\|_{p_{2};\Gamma'_{T}}; \\
&(c') ~~\mbox{for each}~ s\in [0, T) \\
&\hspace*{1.3cm} \sup_{x \in \overline{D(s)}}\left|E_{s,x}\left[\exp Z_{s}(T) h( X(T)); \sigma_{s}(\overline{\Sigma}_{(\infty)}) > T\right]\right| \\
&\hspace*{2cm}\leq \Lambda_{3}\sup_{(s',x') \in \overline{\Omega_{T}}} E_{s',x'}\left[\exp\{r_{3}^{\ast}\|\gamma\|_{\infty}L(T)\}\right]^{1/r_{3}^{\ast}} (T - s)^{-(1/2)(n/r_{3})(1-1/q_{3}^{\ast})}
\|h\|_{p_{3};\Omega}, 
\end{align*}
where $q_{i}^{\ast}$ and $r_{i}^{\ast}$ $(i = 1, 2, 3)$ indicate the conjugate H\"older exponents of $q_{i}$ and $r_{i}$ $(i = 1, 2, 3)$, respectively.  

Now, for $f \in L^{p_{1}}(\Omega_{T})$, $\psi \in L^{p_{2}}(\Gamma'_{T})$ and $h \in L^{p_{3}}(\Omega)$, choose sequences $\{f_{m}\}$, $\{\psi_{m}\}$ and $\{h_{m}\}$ such that   each $f_{m}$, $\psi_{m}$, $h_{m}$ has the same condition for $f$, $\psi$, $h$ in the condition (iv) of Assumption \ref{assump_3}, respectively, and 
\[
f_{m} \to f ~~\mbox{in} ~ L^{p_{1}}(\Omega_{T}), ~
\psi_{m} \to \psi ~~\mbox{in}~ L^{p_{2}}(\Gamma'_{T}), ~
h_{m} \to h ~~\mbox{in}~ L^{p_{3}}(\Omega)  
\]
as $m \to \infty$.
Denote by $v^{D}$ and $v^{D}_{m}$ the right hand side of (\ref{repre}) corresponding to $f, \psi, h$ and $f_{m}, \psi_{m}, h_{m}$, respectively. 
Then $v^{D}_{m}$ converges to $v^{D}$ locally uniformly on $D \cup \partial_{L}D$ as $m \to \infty$. 
Therefore $v^{D}$ has the same continuity property as $v^{D}_{m}$. 
On the other hand, denote by $u^{D}$ and $u^{D}_{m}$ the left hand side of (\ref{repre}) corresponding to $f, \psi, h$ and $f_{m}, \psi_{m}, h_{m}$, respectively. 
Then, using the energy estimate for weak solutions (see Theorem 7.2 in \cite{Kaw10} and also Theorem 3.1 in \cite{KMT07}), we see that $u^{D}_{m} \to u^{D}$ in $V^{0,1}(D)$ as $m \to \infty$. 
Consequently, the equalities between $u^{D}$ and $v^{D}$ in both of the cases of Theorem \ref{main_result} hold and hence the assertion is verified. 
\hspace*{\fill}$\Box$  


{\remark
\label{remark_2}
Note that if the weak solution $u^{D}$ is continuous on $D$, then the equality (\ref{fact_0}) holds everywhere on $D$; hence the equality (\ref{repre}) in the case (i) holds everywhere on $D$ and the equality is verified in a more simple way.  
On the other hand, combining the result on the H\"older continuity of bounded weak solutions of parabolic equations on cylindrical domains with Dirichlet boundary conditions (e.g. Theorem 1.1 on page 419 of \cite{LSU68}) with a localization argument as in the proof of the existence of weak solutions in \cite{Kaw10}, we see that the weak solution $u^{D}$ has a version which is continuous within $D$ and up to the Dirichlet part (except the Robin part and the border between both parts).  
}
\section{A simple application}
\label{sect.appli}

Here we apply the probabilistic representation of weak solutions given in Theorem \ref{main_result} to show the continuity of  a functional concerned with an estimation problem treated in \cite{Kaw10} for the shape of a domain. 

Let $\mathcal{D}$ be the class of domains satisfying the properties described in Condition \ref{cond_D_1} and Assumption \ref{assump_3}. 
Hence the lateral boundary of each domain of $\mathcal{D}$ includes $\Gamma'_{T}$ commonly, but the rest of the lateral boundary may be varied each other and suppose that its shape is unknown. 
The operator $\mathcal{L}$ on $\Omega_{T}$ and $\mathcal{B}$ on $\Gamma'_{T}$ are given 
and the terminal-boundary problem (\ref{eq_1}) on each domain $D \in \mathcal{D}$ is considered with given source terms and terminal value. 

We want to estimate the shape of the unknown portion (i.e. the Dirichlet part $\Sigma$) via thermal data on a part $\Gamma^{\omega}$ of the accessible portion $\Gamma'$ of the boundary of each section $D(t)$ for $t \in [0, T]$; the thermal data is regarded as the boundary value of the weak solution of the initial-boundary value problem corresponding to the terminal-boundary value problem (\ref{eq_1}) on some domain of  $\mathcal{D}$ (see outline of the proof of Theorem \ref{boundedness} in Appendix on the correspondence between forward and backward problems). 
In the sequel, we treat the solution in the backward form as before, that is, the solution to the terminal-boundary value problem (\ref{eq_1}). 
Then consider the following functional $\mathcal{V}(D)$ for $D \in \mathcal{D}$ associated to the inverse problem:
\begin{equation*}
\mathcal{V}(D) := \int_{0}^{T} dt \int_{\Gamma^{\omega}}\left|u^{D}(t, x) - d(t, x)\right|^{2}S(dx),
\end{equation*}
where $d(t, x)$ $(t \in [0, T], x \in \Gamma^{\omega})$ is a given observed data, supposing that 
$d(t, x) = u^{D_{0}}(t, x)$ for some domain $D_{0} \in \mathcal{D}.$
We examine a certain continuity property of the functional $\mathcal{V}(D)$ related to the Hausdroff metric $d_{\mathcal{H}}$ on the space of non-empty compact sets in $\bm{R}^{1+n}$.  
To verify the result, we use the first exit time from some sets:  
For a Borel set $G$ of $\bm{R}^{1+n}$, let   
\[
\tau_{s}(G):= \inf\{t \geq s: (t, X(t)) \notin G\}; 
\]
then set $\tau_{s;T}(G):= \tau_{s}(G) \wedge T$ and also $\sigma_{s;T}(G):= \sigma_{s}(G) \wedge T.$

{\theorem
\label{continuity.func}
Let $D, D_{m} \in \mathcal{D}$ $(m = 1, 2, \dots)$ and $\Sigma, \, \Sigma_{m}$ their Dirichlet parts, respectively. 
Suppose that 
$d_{\mathcal{H}}(\overline{\Sigma_{m}}, \overline{\Sigma}) \to 0$ as $m \to \infty.$ 
Then $\mathcal{V}(D_{m}) \to \mathcal{V}(D)$ as $m \to \infty.$
}
\medskip

\noindent
\textbf{Proof.} 
As in the proof to the second case in Theorem \ref{main_result}, we denote by $v^{D_{m}}$ and $v^{D}$ the right hand side of (\ref{repre}) for $D_{m}$ and $D$, respectively. 
Since $v^{D_{m}}(s, x)$ $(m = 1, 2, \dots; (s, x) \in \Gamma'_{T})$ are uniformly bounded,  for verifying the assertion of the theorem, it is enough to show 
\begin{equation}
v^{D_{m}}(s, x) \longrightarrow v^{D}(s, x) ~~\mbox{as} ~ m \to \infty 
\label{conv.prop}
\end{equation}
for each $(s, x) \in \Gamma'_{T}.$ 
Put $\tilde{D}:= D \cup \Pi_{T} \cup \Gamma'_{T}$, and for $\varepsilon > 0$ set 
\begin{align*}
&\tilde{D}^{+\varepsilon}:= \{(t, y) \in \Omega_{T}: d((t, y), \tilde{D}) <  \varepsilon\} \cup \tilde{D},  
\\
&\tilde{D}^{-\varepsilon}:= \{(t, y) \in \overline{\Omega_{T}}: d((t, y), \overline{\Omega_{T}} \setminus \tilde{D}) \geq \varepsilon\}. 
\end{align*}
By assumption, for any $\varepsilon > 0,$ 
$
\overline{\Sigma_{m}} \subset \overline{\tilde{D}^{+\varepsilon}} \setminus \tilde{D}^{-\varepsilon} 
$
for sufficiently large $m \in \bm{N}.$ 
Take a sufficiently small $\varepsilon > 0.$ Then we see that, for $(s, x) \in (\Gamma'\setminus\Pi^{\varepsilon})_{T},$  
\[
P_{s,x}\left(\tau_{s;T}(\tilde{D}^{-\varepsilon}_{(\infty)}) \leq \sigma_{s;T}(\overline{\Sigma}_{m(\infty)})\leq \tau_{s;T}(\tilde{D}^{+\varepsilon}_{(\infty)})\right) = 1
\]
for sufficiently large $m.$ 
On the other hand, by using the continuity of the sample paths of $X(\cdot)$ and that each point of 
$\Sigma_{2}$ is regular for the set $\Omega_{T}\setminus D$ (see \ref{subsect.continuity.hitting.time} in details), it holds 
\[
P_{s,x}\left(\lim_{\varepsilon\downarrow 0}\tau_{s;T}(\tilde{D}^{-\varepsilon}) = 
\lim_{\varepsilon\downarrow 0}\tau_{s;T}(\tilde{D}^{+\varepsilon}) = 
\sigma_{s;T}(\overline{\Sigma}_{(\infty)})\right) = 1
\]
for $(s, x) \in \Gamma'_{T}.$ 
Accordingly we have 
\[
P_{s,x}\left(\lim_{m \to \infty}\sigma_{s;T}(\overline{\Sigma}_{m(\infty)}) = 
\sigma_{s;T}(\overline{\Sigma}_{(\infty)}) \right) = 1
\]
for $(s, x) \in \Gamma'_{T};$ 
that is, the convergence property (\ref{conv.prop}) is verified. 
\hspace*{\fill}$\Box$ 

\medskip
Using the continuity of the functional $\mathcal{V}(D)$ and the result on the uniqueness in shape identification in \cite{Kaw10}, we can show that any minimizing sequence for the functional converges to a unique optimal domain in some cases, under certain a priori information. 
Furthermore, based on the functional $\mathcal{V}(D)$ and the probabilistic representation of solutions, the paper \cite{Kaw15} gives an algorithm for reconstruction of the position and shape of unknown cavities. 

{\remark
\label{remark_3}
A probabilistic cost function which plays the same role as $\mathcal{V}(D)$ is given by 
\[
E_{0,\nu}\left[\int_{0}^{T}t^{\gamma}1_{\Gamma^{\omega}}(X(t))\left|u^{D}(t, X(t)) - d(t, X(t))\right|^{2}L(dt)\right],
\]
where the symbol ``$E_{0,\nu}$'' designates the expectation with respect to the probability measure 
\[
\displaystyle{P_{0,\nu}( \cdot ):= \int_{\Gamma^{\omega}}P_{0,x}( \cdot )\nu(dx)}
\]
with a given initial distribution $\nu$ on $\Gamma^{\omega}$ and $\gamma \geq 0$; the initial distribution $\nu$ and the weight $t^{\gamma}$ may be chosen according to actual applications. 
Such a cost function driven by local times is used in \cite{Kot04} to provide an algorithm for an inverse problem determining coefficients of thermal conductivity.  
}

\section*{Acknowledgments}
The author is grateful to  Professors Hajime Kawakami and Kiyomasa Narita for their valuable comments and useful discussions with them.
Particularly, in Appendix B, Kawakami provides an exposition on an approach based on a Bayesian framework to the shape identification inverse problem mentioned in Section 4.  

\appendix
\section{Details of the proof}
\label{sect.appendix}
\subsection{Boundedness of weak solutions} 
\label{subsect.appendix.boundedness}

\noindent
The following boundedness result for weak solutions is proved by the method given by 
Lady\u{z}enskaja et al in \cite{LSU68}. 
They treat the case of a Dirichlet boundary condition and the case without boundary condition on a cylindrical domain in the time-space. 
However, their method is applicable to the case of a mixed boundary condition on a non-cylindrical domain under some additional considerations. 
{\theorem
\label{boundedness}
In addition to Definition \ref{cond_D_1} and Assumption \ref{assump_2}(i), suppose that the following integrability conditions for the coefficients and the source terms of the terminal-boundary problem (\ref{eq_1}) are fulfilled: 
\begin{enumerate}
\item[(i)]  
$\displaystyle{
\mu \equiv 
\left\||\bm{a}|^{2}\right\|_{q,r;D} \vee 
\left\||\bm{b}|^{2}\right\|_{q,r;D} \vee 
\left\||\bm{f}|^{2}\right\|_{q,r;D} \vee 
\|a\|_{q,r;D} \vee 
\|f\|_{q,r;D} 
<\infty}$ 

\medskip

for a pair $(q, r)$ satisfying the relationships 
\[
\frac{1}{r} + \frac{n}{2q} = 1 - \kappa, ~~ q \in \left[\frac{n}{2(1-\kappa)}, \infty\right], ~~
r \in \left[\frac{1}{1-\kappa}, \infty\right]
\]
for some $\kappa \in (0, 1)$, 
and further 
$\| \|A\|\|_{\infty; D} \vee \||\bm{a}|\|_{\infty; D} \vee 
\||\bm{b}| \|_{\infty; D} 
< \infty$  and $\|
|\bm{f}| \|_{2; D} < \infty$;  \\
\item[(ii)] \hspace*{1cm} 
$\displaystyle{
\tilde{\mu} \equiv 
\|\sigma\|_{\tilde{q},\tilde{r};\Gamma^{\prime}_{T}} \vee 
\|\psi\|_{\tilde{q},\tilde{r};\Gamma^{\prime}_{T}} \vee 
\||\gamma_{\Gamma'}\bm{f}|\|_{\tilde{q},\tilde{r};\Gamma^{\prime}_{T}} 
<\infty}$ 

\medskip

for a pair $(\tilde{q}, \tilde{r})$ satisfying the relationships 
\[
\frac{1}{\tilde{r}} + \frac{n-1}{2\tilde{q}} = \frac{1 - \tilde{\kappa}}{2}, ~~ \tilde{q} \in \left[\frac{n-1}{1-\tilde{\kappa}}, \infty\right], ~~
\tilde{r} \in \left[\frac{2}{1-\tilde{\kappa}}, \infty\right]
\]
for some $\tilde{\kappa} \in (0, 1)$. 
\end{enumerate}
Moreover assume that $\displaystyle{\mathrm{ess\, sup}_{x\in D(T)}h(x) \leq \check{k}}$ with 
$\check{k} \geq 0.$ 
Then, for any weak solution $u$ in $V^{0,1}(D)$ to the problem (\ref{eq_1}), there is a positive constant $K$ such that $\displaystyle{\mathrm{ess\, sup}_{(t,x)\in D}u(t, x) \leq K}.$ 
}

{\remark
\label{remark_A1} 
The additional conditions for $A, ~\bm{a}, ~\bm{b}$ and $\bm{f}$ in (i) of Theorem \ref{boundedness} are unnecessary when $D$ is a cylindrical domain, that is, $D = \Omega_{T}$ and hence $\Sigma = [0, T] \times \Gamma''$.
}

\medskip
\medskip

We begin with preparing key results corresponding to those in \cite{LSU68}. 
Let 
\begin{align*}
\stackrel{\circ}{V}\!{ }_{\Sigma}^{0,1}(D):= \{u \in V^{0,1}(D): \left.u\right|_{\Sigma} = 0\}, ~~
\stackrel{\circ}{H}\!{ }_{\Sigma}^{1,1}(D):= \{u \in H^{1,1}(D): \left.u\right|_{\Sigma} = 0\}. 
\end{align*}
For $u \in \stackrel{\circ}{V}\!{ }_{\Sigma}^{0,1}(D),$ denote by $\bar{u}$ the extension of $u$ by zero: 
\begin{align*}
\bar{u} = 
\begin{cases}
u & \text{on} ~ D,
\\
0 & \text{on} ~ D^{c}.
\end{cases}
\end{align*}
Then $\bar{u} \in V^{0,1}(\Omega_{T})$ with the properties:
\begin{align*}
&\|\bar{u}\|_{V^{0,1}(\Omega_{T})} = \|u\|_{V^{0,1}(D)}; ~~
\left.\bar{u}\right|_{\Gamma^{\prime}_{T}} = \left.u\right|_{\Gamma^{\prime}_{T}}, ~~
\left.\bar{u}\right|_{\Gamma^{\prime\prime}_{T}} = 0; \\
&\mathrm{ess\, sup}_{(t,x)\in \Omega_{T}}\bar{u}(t, x) = \mathrm{ess\, sup}_{(t,x)\in D}u(t, x)
\end{align*}
(see Proof of Lemma 4.1 and Lemma A.1 in \cite{Kaw10}; see also Example 10.2.1 in \cite{ABM06}).
Therefore it is enough to show that 
$\mathrm{ess\, sup}_{(t,x)\in \Omega_{T}}\bar{u}(t, x) \leq K$ 
for an arbitrarily given weak solution $u.$ 

First we recall the key estimates (3.3), (3.8) and (3.10), (3.11) in Chapter II of \cite{LSU68} for functions in $V^{0,1}(\Omega_{T})$ as the following lemmas. 
{\lemma
\label{key.lemma.1}
Take a pair $(q, r)$ with the following conditions:
\begin{align*}
&q \in \left[2, \frac{2n}{n - 2}\right], ~r \in [2, \infty) ~~\mbox{ if } n \geq 3, \\
&q \in [2, \infty), ~r \in (2, \infty] ~~\mbox{ if } n = 2
\end{align*}
and 
\[
\frac{1}{r} + \frac{n}{2q} = \frac{n}{4}.
\]
Then there is a positive constant $\beta$ such that 
\[
\|v\|_{q,r;\Omega_{T}} \leq \beta \|v\|_{V^{0,1}(\Omega_{T})}
\]
for every $v \in V^{0,1}(\Omega_{T}).$
}

{\lemma
\label{key.lemma.2}
Take a pair $(\tilde{q}, \tilde{r})$ with the following conditions:
\begin{align*}
&\tilde{q} \in \left[\frac{2(n - 1)}{n}, \frac{2(n - 1)}{n - 2}\right], ~\tilde{r} \in [2, \infty] ~~\mbox{ if } n \geq 3, \\
&\tilde{q} \in [1, \infty), ~\tilde{r} \in (2, \infty] ~~\mbox{ if } n = 2
\end{align*}
and 
\[
\frac{1}{\tilde{r}} + \frac{n - 1}{2\tilde{q}} = \frac{n}{4}.
\]
Then there is a positive constant $\tilde{\beta}$ such that 
\[
\|v\|_{\tilde{q},\tilde{r};\Gamma_{T}} \leq \tilde{\beta} \|v\|_{V^{0,1}(\Omega_{T})}
\]
for every $v \in V^{0,1}(\Omega_{T}).$ 
}


\medskip
\medskip

Given $v \in V^{0,1}(\Omega_{T}),$ for each $t \in (0, T),$ we can take the following version 
$\tilde{v}$ which is quasi-continuous with respect to the space variable: 
\[
\tilde{v}(t, x):= \lim_{r \to 0}\, -\!\!\!\!\!\!\int_{B_{r}(x) \cap \Omega} v(t, y) dy ~~ \mbox{ for } ~x \in \overline{\Omega}.
\]
In what follows, we always take the version for any element of $V^{0,1}(\Omega_{T}).$ 
Hence we can regard as  $v = \gamma_{\Gamma} v$ on the lateral boundary.

For $v \in V^{0,1}(\Omega_{T})$, $k \in \bm{R}$ and $t \in (0, T)$ we set 
\begin{align*}
v^{(k)} &:= (v - k) \vee 0, \\
A_{k}(t) &:= \{x \in \Omega: v(t, x) > k\}, \\ 
\tilde{A}_{k}(t) &:= \{x \in \Gamma: v(t, x) > k\}.
\end{align*}
Moreover, define 
\begin{align*}
&\mu(k) \equiv \mu_{T}(k):= \int_{0}^{T}\mathrm{mes}^{r/q}(A_{k}(t)) dt, \\
&\tilde{\mu}(k) \equiv \tilde{\mu}_{T}(k):= \int_{0}^{T}\widetilde{\mathrm{mes}}^{\tilde{r}/\tilde{q}}(\tilde{A}_{k}(t)) dt,
\end{align*}
where ``$\mathrm{mes}$'' and ``$\widetilde{\mathrm{mes}}$'' denote the $n$--dimensional Lebesgue measure and the surface measure on $\Gamma,$ respectively. 
Then we have the following result which is a miner modification of Theorem 6.1 in Chapter II of \cite{LSU68} (see also Remark 6.2 for the theorem), and is proved in the same way.

{\theorem
\label{key.theorem}
Suppose that for a given $v \in V^{0,1}(\Omega_{T})$ there exists $\hat{k} \geq 0$ satisfying the following condition: for a pair $(q, r)$ with the conditions in Lemma \ref{key.lemma.1} and for a pair $(\tilde{q}, \tilde{r})$ with the conditions in Lemma \ref{key.lemma.2}, it can be chosen positive constants $\varrho, ~\lambda, ~\tilde{\lambda}$ such as 
\begin{equation}
\|v^{(k)}\|_{V^{0,1}(\Omega_{T})} \leq \varrho k \left\{\mu^{(1+\lambda)/r}(k) + \tilde{\mu}^{(1+\tilde{\lambda})/\tilde{r}}(k)\right\}
\label{key-inequality}
\end{equation}
for every $k \geq \hat{k}.$
Then, for some constant $m = m_{T} > 1,$ it holds 
\[
\mathrm{ess\, sup}_{(t,x) \in \Omega_{T} } v(t, x) \leq 2m\hat{k};
\]
such a constant $m$ is given by 
\begin{align}
&m \equiv m_{T} := 1 + (\beta + \tilde{\beta})\varrho\, \zeta^{1/(\lambda\wedge\tilde{\lambda})} 2^{1/(\lambda\wedge\tilde{\lambda})^{2}} \nonumber \\
&\times \Bigl\{T^{(1+\lambda)/r}\mathrm{mes}^{(1+\lambda)/q}(\Omega) + T^{(1+\tilde{\lambda})/\tilde{r}}\widetilde{\mathrm{mes}}^{(1+\tilde{\lambda})/\tilde{q}}(\Gamma)\Bigr\},
\label{m.def}
\end{align}
where 
\begin{align*}
\zeta \equiv \zeta_{T}:= 
\begin{cases}
4(\beta + \tilde{\beta})\varrho 2^{\lambda}\Bigl\{T^{(\tilde{\lambda}-\lambda)/\tilde{r}}\widetilde{\mathrm{mes}}^{(\tilde{\lambda}-\lambda)/\tilde{q}}(\Gamma) + 1\Bigr\} 
& if ~ \lambda < \tilde{\lambda} \\ 
4(\beta + \tilde{\beta})\varrho
& if ~  \lambda = \tilde{\lambda}  \\
4(\beta + \tilde{\beta})\varrho 2^{\tilde{\lambda}}\Bigl\{T^{(\lambda-\tilde{\lambda})/r}\mathrm{mes}^{(\lambda-\tilde{\lambda})/q}(\Omega) + 1\Bigr\} 
& if ~ \lambda > \tilde{\lambda}.
\end{cases}
\end{align*}
}
\medskip 

\noindent
\textbf{Outline of the proof of Theorem \ref{boundedness}.}  

It is enough to show the boundedness of the extension $\bar{u}$ of the weak solution $u.$ 
For this purpose, we verify that $\bar{u}$ fulfills the condition on $v$ in Theorem \ref{key.theorem}. 
Replacing each test function $\eta \in \stackrel{\circ}{H}\!\!{ }_{\Sigma}^{1,1}(D)$ with 
$\tilde{\mu}_{\delta}\eta$ for $\delta \in (0, \frac{1}{2} d(\Sigma_{(2)}, \Gamma'_{T}))$ in the weak form (\ref{eq.weak.form.I}) for the weak solution $u$, we see that $u_{\delta} = \tilde{\mu}_{\delta}u$ becomes a weak solution to the terminal-boundary value problem (\ref{eq_1}) by replacing $\bm{f}$, $f$ and $h$ with $\bm{g},$  $g$ and $\tilde{h}$ respectively, 
where
\begin{align}
&\bm{g}:= \tilde{\mu}_{\delta}\bm{f} - u A\nabla_{x}\tilde{\mu}_{\delta}, \nonumber\\
&g:= \tilde{\mu}_{\delta}f + u\partial_{t}\tilde{\mu}_{\delta} + (A\nabla_{x}u + \bm{a}u - \bm{b}u +  \bm{f})\iprod\nabla_{x}\tilde{\mu}_{\delta}, \nonumber\\
&\tilde{h}:= \tilde{\mu}_{\delta}(T, \cdot)h. 
\label{source.term}
\end{align}
Moreover, since $\tilde{\mu}_{\delta}\eta \in \stackrel{\circ}{H}\!\!{ }_{\Sigma}^{1,1}(D)$ for 
$\eta \in \stackrel{\circ}{H}\!\!{ }_{\Gamma^{\prime\prime}_{T}}^{1,1}(\Omega_{T})$ (which is defined in the same way as $\stackrel{\circ}{H}\!\!{ }_{\Sigma}^{1,1}(D)$), we also see that $\bar{u}_{\delta} = \tilde{\mu}_{\delta}\bar{u}$ becomes a weak solution in 
$V^{0,1}(\Omega_{T})$ to the following terminal-boundary value problem on $\Omega_{T}:$ 
\begin{equation}
\begin{cases}
\mathcal{P}u(t,x)= -\nabla_x \iprod \bm{g}(t,x) + g(t, x)  
& \text{if} ~ (t,x) \in \Omega_{T}
\\
\mathcal{B}u(t,x) = - \psi(t,x)
& \text{if} ~ (t,x) \in \Gamma^{\prime}_T
\\
u(t,x)=0 & \text{if} ~ (t,x) \in \Gamma^{\prime\prime}_T
\\
u(T,x)= \tilde{h}(x) & \text{if} ~ x \in \Omega.
\end{cases}
\label{eq.A.1}
\end{equation}
Here $u$ in $\bm{g}$ and $g$ in (\ref{source.term}) is replaced  
with $\bar{u}.$  

In the following, to obtain the boundedness result, we treat the problem (\ref{eq.A.1}) and its weak form in the forward form. 
For this purpose, we consider the time reversion $t \longrightarrow T - t$ on the operators, solutions, etc., that is, we set 
$\mathcal{L}_{T}:= \mathcal{L}_{x}(T-t)$,  $\mathcal{B}_{T}:= \mathcal{B}_{x}(T-t)$, $u_{T}(t, x):= u(T-t, x)$, $A_{T}(t, x):= A(T-t, x)$, etc. 
Then $(\bar{u}_{\delta})_{T}\in V^{0,1}(\Omega_{T})$ is a weak solution to the following initial-boundary value problem on $\Omega_{T}:$ 
\begin{equation}
\begin{cases}
\frac{\partial v}{\partial t} - \mathcal{L}_{T}v(t,x)= \nabla_x \iprod \bm{g}_{T}(t,x) - g_{T}(t, x)  
& \text{if} ~ (t,x) \in \Omega_{T}
\\
\mathcal{B}_{T}v(t,x) = - \psi_{T}(t,x)
& \text{if} ~ (t,x) \in \Gamma^{\prime}_T
\\
v(t,x)=0 & \text{if} ~ (t,x) \in \Gamma^{\prime\prime}_T
\\
v(0,x)= \tilde{h}(x) & \text{if} ~ x \in \Omega.
\end{cases}
\label{eq.A.2}
\end{equation}
In the rest of the proof, for notational simplicity, we drop the subindex ``$T$'' and suppose $k \geq 0.$ 
Then we can take the Steklov average $\bigl(\bar{u}_{\delta}^{(k)}\bigr)_{h}$ of 
$\bar{u}_{\delta}^{(k)}:= (\bar{u}_{\delta})^{(k)}$ as a test function to the weak form for the weak solution $\bar{u}_{\delta}$ to the problem (\ref{eq.A.2}); 
so that using the argument on page 141 in \cite{LSU68} and then taking limit as $h \to 0,$ we get 
\begin{align*} 
&\left.
\frac{1}{2}\int_\Omega \left(\bar{u}_{\delta}^{(k)}(t,x)\right)^{2} dx
\right|_{t=t_1}^{t=t_2}
+
\int_{\Omega_{t_1,t_2}}
\Bigl\{\left(A\nabla_x\bar{u}_{\delta} 
+ \bm{a}\bar{u}_{\delta} + \bm{g}\right)\iprod \nabla_x\bar{u}_{\delta}^{(k)} \nonumber \\
\nonumber 
&\hspace*{5cm}+ \left(\bm{b} \iprod \nabla_x\bar{u}_{\delta} 
+ a \bar{u}_{\delta} + g \right)\bar{u}_{\delta}^{(k)} \Bigr\} dt dx 
\nonumber 
\\ 
&\hspace*{3cm}+ \int_{\Gamma'_{t_1,t_2}}
\left(\sigma \bar{u}_{\delta} + \psi + \gamma_{\Gamma'}\bm{f}\iprod\bm{n}\right)\gamma_{\Gamma'}\bar{u}_{\delta}^{(k)} dt S(dx) = 0
\end{align*}
for every $0 \leq t_{1} \leq t_{2} \leq T.$ 
Noting $\bar{u}_{\delta} \to \bar{u}$ in $V^{0,1}(\Omega_{T})$ as $\delta \to 0$ and using similar treatments to the second term of the right hand side of (\ref{strong_form})  (see also the proof of Lemma A1 in \cite{Kaw10}), we see that 
\begin{align*} 
&\int_{\Omega_{t_1,t_2}}(\bm{g} - \tilde{\mu}_{\delta}\bm{f})\iprod \nabla_{x}\bar{u}_{\delta}^{(k)} \, dt dx = - \int_{\Omega_{t_1,t_2}}\bar{u} A\nabla_{x}\tilde{\mu}_{\delta}\iprod \nabla_{x}\bar{u}_{\delta}^{(k)} \,  dt dx \to 0, \\ 
&\int_{\Omega_{t_1,t_2}}(g - \tilde{\mu}_{\delta} f)\bar{u}_{\delta}^{(k)} dt dx =  \int_{\Omega_{t_1,t_2}}\left\{\bar{u}\,\partial_{t}\tilde{\mu}_{\delta} + (A\nabla_{x}\bar{u} + \bm{a}\bar{u} - \bm{b}\bar{u} +  \bm{f})\iprod\nabla_{x}\tilde{\mu}_{\delta}\right\}\bar{u}_{\delta}^{(k)} dt dx  \to 0  
\end{align*} 
as $\delta \to 0;$ 
this yields the following key equality:
\begin{align} 
&\left.
\frac{1}{2}\int_\Omega \left(\bar{u}^{(k)}(t,x)\right)^{2} dx
\right|_{t=t_1}^{t=t_2}
+
\int_{\Omega_{t_1,t_2}}
\Bigl\{\left(A\nabla_x\bar{u} 
+ \bm{a}\bar{u} + \bm{f}\right)\iprod \nabla_x\bar{u}^{(k)} \nonumber \\
\nonumber 
&\hspace*{5cm}+ \left(\bm{b} \iprod \nabla_x\bar{u} 
+ a \bar{u} + f \right)\bar{u}^{(k)}\Bigr\} dt dx 
\nonumber 
\\ 
&\hspace*{3cm}+ \int_{\Gamma'_{t_1,t_2}}
\left(\sigma \bar{u} + \psi + \gamma_{\Gamma'}\bm{f}\iprod\bm{n}\right)\gamma_{\Gamma'}\bar{u}^{(k)} dt S(dx) = 0
\label{key.equality.2}
\end{align}
for every $0 \leq t_{1} \leq t_{2} \leq T.$
Let 
\begin{align*} 
&\mathcal{D}:= 2\left\{\frac{4}{\nu}|\bm{a}|^{2} + \frac{4}{\nu}|\bm{f}|^{2} + \frac{2}{\nu}|\bm{b}|^{2} + 2|a| + |f|\right\}, \\
&\mathcal{E}:= 2\left\{2|\sigma| + |\psi| +  |\gamma_{\Gamma'}\bm{f}|\right\}. 
\end{align*}
Then we have 
\begin{align*} 
\left(\frac{1}{2} \wedge \frac{\nu}{2}\right)\|\bar{u}^{(k)}\|_{V^{0,1}(\Omega_{t_{1}})} &\leq 
\|\mathcal{D}\|_{q,r;\Omega_{t_{1}}(k)} \|(\bar{u} - k)^{2} + k^{2}\|_{q',r';\Omega_{t_{1}}(k)} \\
& ~~~\qquad  + \|\mathcal{E}\|_{\tilde{q},\tilde{r};\Gamma'_{t_{1}}(k)} \|(\bar{u} - k)^{2} + k^{2}\|_{\tilde{q}',\tilde{r}';\Gamma'_{t_{1}}(k)}~ , 
\end{align*}
where 
\begin{align*} 
&q':= \frac{q}{q-1}, ~~r':= \frac{r}{r-1}, ~~\tilde{q}':= \frac{\tilde{q}}{\tilde{q}-1}, ~~\tilde{r}':= \frac{\tilde{r}}{\tilde{r}-1}, ~~ \\
&\Omega_{t_{1}}(k):= \{(t, x) \in \Omega_{t_{1}}: \bar{u}(t, x) > k\}, \\
&\Gamma'_{t_{1}}(k):= \{(t, x) \in \Gamma'_{t_{1}}: \bar{u}(t, x) > k\} ~(= \Gamma_{t_{1}}(k):= \{(t, x) \in \Gamma_{t_{1}}: \bar{u}(t, x) > k\}).
\end{align*}
Now letting $\hat{q}:= 2(1 + \vartheta)q'$ and $\hat{r}:= 2(1 + \vartheta)r'$ with 
$\displaystyle{\vartheta:= \frac{\kappa}{2n}},$ we find that the pair $(\hat{q}, \hat{r})$ fulfills the conditions on a pair 
$(q, r)$ in Lemma \ref{key.lemma.1}. 
Similarly, letting $\check{q}:= 2(1 + \tilde{\vartheta})\tilde{q}'$ and 
$\check{r}:= 2(1 + \tilde{\vartheta})\tilde{r}'$ with 
$\displaystyle{\tilde{\vartheta}:= \frac{\tilde{\kappa}}{n}},$ we also find that the pair 
$(\check{q}, \check{r})$ fulfills the conditions on a pair 
$(\tilde{q}, \tilde{r})$ in Lemma \ref{key.lemma.2}.  
Therefore, in the same way as in \cite{LSU68}, we see the following: taking $t_{1} > 0$ so small that 
\begin{align*} 
\beta^{2}\|\mathcal{D}\|_{q,r;\Omega_{t_{1}}}t_{1}^{2\vartheta/\hat{r}}\mathrm{mes}^{2\vartheta/\hat{q}}(\Omega) + \tilde{\beta}^{2}\|\mathcal{E}\|_{\tilde{q},\tilde{r};\Gamma'_{t_{1}}}t_{1}^{2\tilde{\vartheta}/\check{r}}\widetilde{\mathrm{mes}}^{2\tilde{\vartheta}/\check{q}}(\Gamma') \leq \frac{1}{2}\left(\frac{1}{2} \wedge \frac{\nu}{2}\right), 
\end{align*} 
we have 
\[
\|\bar{u}^{(k)}\|_{V^{0,1}(\Omega_{t_{1}})}^{2} \leq \varrho k \Bigl\{\mu_{t_{1}}^{(1+\vartheta)/\hat{r}}(k)  + \tilde{\mu}_{t_{1}}^{(1+\tilde{\vartheta})/\check{r}}(k)\Bigr\}~~\mbox{for}~k \geq \hat{k} :=\check{k} \vee 1 
\]
with 
\[
\varrho:= \left(\frac{1}{4} \wedge \frac{\nu}{4}\right)^{-1/2}\Bigl\{\Bigl(\|\mathcal{D}\|_{q,r;\Omega_{t_{1}}}\Bigr)^{1/2} + \left(\|\mathcal{E}\|_{\tilde{q},\tilde{r};\Gamma'_{t_{1}}}\right)^{1/2}\Bigr\}.
\]
This shows that Theorem \ref{key.theorem} holds for $\bar{u}$ by replacing $T$ with $t_{1}$,  
and hence  
\[
\mathrm{ess~sup}_{(t,x)\in \Omega_{t_{1}}}\bar{u}(t, x) \leq 2m_{t_{1}}\hat{k},
\]
where $m_{t_{1}}$ indicates the constant defined in (\ref{m.def}) by replacing $T$ with 
$t_{1}.$
In the case $t_{1} < T,$ we divide the interval $(0, T)$ into a finite number of intervals with length less than $t_{1}$, say $L$ the number of the intervals, and then, noting $m_{t_{1}} \leq m_{T}$,  we conclude that 
\[
\mathrm{ess~sup}_{(t,x)\in \Omega_{T}}\bar{u}(t, x) \leq (2m_{T})^{L}\hat{k}. 
\]
Consequently the theorem is proved.
\hspace*{\fill}
$\Box$ 

\subsection{Regularity of weak solutions} 
\label{subsect.appendix.regularity}

\noindent
As described in the outline of the proof of Theorem \ref{boundedness}, the extension 
$\bar{u}^{D}_{\delta}:= \overline{u^{D}_{\delta}}$ of $u^{D}_{\delta}$ for $u = u^{D}$ becomes a weak solution to the terminal-boundary value problem (\ref{eq.A.1}) on the cylindrical domain $\Omega_{T}.$ 
Moreover, since $\bar{u}^{D}A\nabla_{x}\tilde{\mu}_{\delta} = \bar{u}^{D}_{\delta/2}A\nabla_{x}\tilde{\mu}_{\delta} \in V^{0,1}(\Omega_{T})^{n}$ under Assumption \ref{assump_3} (iii) and since its boundary value $\left.\bar{u}^{D}A\nabla_{x}\tilde{\mu}_{\delta}\right|_{\Gamma_{T}} = \bm{0},$ 
$\bar{u}^{D}_{\delta}$ is a weak solution to the following terminal-boundary value problem on 
$\Omega_{T}:$ 
\begin{equation*}
\begin{cases}
\mathcal{P}u(t,x)= \bar{g}(t, x)  
& \text{if} ~ (t,x) \in \Omega_{T}
\\
\mathcal{B}u(t,x) = - \psi(t,x)
& \text{if} ~ (t,x) \in \Gamma^{\prime}_T
\\
u(t,x)=0 & \text{if} ~ (t,x) \in \Gamma^{\prime\prime}_T
\\
u(T,x)= \tilde{h}(x) & \text{if} ~ x \in \Omega
\end{cases}
\end{equation*}
with $\bar{g}(t, x):= \nabla_x \iprod(\bar{u}^{D}A\nabla_{x}\tilde{\mu}_{\delta})(t,x) + g(t, x),$ because we suppose $\bm{f} = \bm{0}.$ 
Therefore, noting Assumption \ref{assump_3} and then using Theorem 3.3 in \cite{KMT07}, we see that $\bar{u}^{D}_{\delta} \in H^{1,1}(\Omega_{T})$ and further, from the results on boundary regularity of weak solutions in Theorem 3.4 of \cite{KMT07} and the results on interior regularity of weak solutions, it follows that $\bar{u}^{D}_{\varepsilon,\delta} \in H^{1,2}(\Omega_{T}).$ 
Thus the regularity (\ref{reg.sol}) for $u^{D}_{\varepsilon,\delta}$ is verified. 

\subsection{Proof of $P_{s,x}(\sigma_{s}(\Pi) = \infty) = 1$ ~$((s, x) \in D)$} 
\label{subsect.appendix.proof.hitting.time.1}

\noindent
More generally, we show the following. 
{\proposition
\label{hitting.time.1}
Suppose that $\Xi$ is a closed set in $[0, \infty) \times \Gamma$ and that $S(\Xi(t)) = 0$ for almost every $t \geq 0.$ 
Then $P_{s,x}(\sigma_{s}(\Xi) = \infty) = 1$ for $(s, x) \in \left([0, \infty) \times \overline{\Omega}\right) \setminus \Xi.$ 
}
\medskip 

We begin with noting the precise form of the density of the joint distribution of the entrance time and the entrance place of the process $\{X(t)\}$ on the boundary $\Gamma.$ Let $\sigma_{s}(\Gamma) = \inf \{t \geq s: X(t) \in \Gamma\}$ and denote by 
$p_{0}(s, x; t, y)$ the fundamental solution of the equation $\mathcal{P}_{0}u = 0$ on $\Omega$ 
with the homogeneous Dirichlet boundary condition  (see \cite{Ito92} for the fundamental solution).  
Then $P_{s,x}(\sigma_{s}(\Gamma) < \infty) = 1$ for $s \geq 0$ and $x \in \Omega$ (see Corollary 3.3 and its consequence in \cite{Kha12}) and 
\[
P_{s,x}(X(t) \in dy; \sigma_{s}(\Gamma) > t) = p_{0}(s, x; t, y)dy ~~~(0 \leq s < t, x \in \Omega, dy \subset 
\Omega). 
\]
Now we denote the formal adjoint operator of $\mathcal{L}_{0}$ by $\mathcal{L}_{0}^{\ast}$. For the operators, we specify the time variable $s$ or $t$ of those coefficients and the space variable $x$ or $y$ acting on a test function $\varphi(s, x; t, y)$ under the operations such as 
$\mathcal{L}_{0} = \mathcal{L}_{0;x}(s)$ and 
$\mathcal{L}_{0}^{\ast} = \mathcal{L}^{\ast}_{0;y}(t).$ Then we know the following (cf. \cite{Ito92}).

{\lemma
\label{hitting _time_density}
For every $(s, x) \in [0, \infty) \times \Omega,$ it holds 
\[
P_{s,x}((\sigma_{s}(\Gamma), X(\sigma_{s}(\Gamma))) \in (dt, dy)) = \frac{\partial}{\partial\mathcal{N}^{\ast}_{y}(t)}p_{0}(s, x; t, y)dtS(dy),
\]
where $\partial/ \partial\mathcal{N}^{\ast}_{y}(t)$ denotes the conormal derivative associated with $\mathcal{L}^{\ast}_{0;y}(t).$
}
\medskip 

\noindent  
{\bf Proof of Proposition \ref{hitting.time.1}.} 

We first note that for $s \geq 0$ and $x \in \Gamma$ it holds 
\begin{align}
P_{s,x}(X(t) \in \Omega ~\mbox{ for some } t \in (a, a+h) ~\mbox{ with any } a \geq s ~\mbox{ and } h > 0) = 1,
\label{sample_path}
\end{align}
because the sample paths $X(\cdot)$ are continuous, $\Omega$ is an open set  and $P_{s,x}(X(t) \in \Gamma) = 0$ for every $ t > s.$ 
To verify the proposition, we divide the case $x \in \overline{\Omega} \setminus \Xi(s)$ into two cases: the case of $x \in \Omega$ and the case of $x \in \Gamma \setminus \Xi(s).$  
We only treat the first case, since the second case can be treated in the same way.  
Take a sequence $\{\varepsilon_{m}\}_{m=1}^{\infty}$ with $\varepsilon_{m} \downarrow 0$  $(m \to \infty),$ and set 
\[
\Omega_{m}:= \Omega^{-\varepsilon_{m}} = \{y \in \Omega: d(y, \Gamma) \geq \varepsilon_{m}\}.
\]
Without generality, we may suppose that $x \in \Omega_{1} (\subset \Omega_{m}).$ 
For $s \geq 0$ and $m = 1, 2, \ldots,$ define the hitting times $\sigma^{(k)}(m)$ and 
$\tau^{(k)}(m)$ $(k = 0, 1, 2, \ldots)$ as follows: 
\begin{align*}
&\sigma^{(0)}(m):= \sigma_{s}(\Gamma), ~~~\tau^{(0)}(m):= s; \\
&\tau^{(k+1)}(m):= \sigma_{\sigma^{(k)}(m)}(\Omega_{m}), ~~~\sigma^{(k)}(m):= \sigma_{\tau^{(k)}(m)}(\Gamma).
\end{align*}
By using (\ref{sample_path}), we see that 
\begin{align*}
P_{s,x}(\sigma_{s}(\Xi) < \infty) 
\leq \sum_{m=1}^{\infty}\sum_{k=0}^{\infty}P_{s,x}(\sigma_{s}(\Xi)  = \sigma^{(k)}(m) < \infty).
\end{align*}
Since $\Xi$ is closed, using the strong Markov property,
\begin{align*}
&P_{s,x}(\sigma_{s}(\Xi)  = \sigma^{(k)}(m) < \infty) \\ 
&\leq P_{s,x}((\sigma^{(k)}(m), X(\sigma^{(k)}(m)) \in \Xi) \\
&= E_{s,x}\left[\left.P_{t,y}\left((\sigma_{t}(\Gamma), X(\sigma_{t}(\Gamma)) \in \Xi\right)\right|_{t=\tau^{(k)}(m), y=X(\tau^{(k)}(m))}; \tau^{(k)}(m) < \infty\right].
\end{align*}
If $\tau^{(k)}(m) < \infty,$ then $y=X(\tau^{(k)}(m)) \in \Omega.$ 
Hence from the assumption that $S(\Xi(r)) = 0$ for almost every $0 \leq r < \infty$ and Lemma \ref{hitting _time_density}, it follows that 
\[
P_{s,x}((\sigma^{(k)}(m), X(\sigma^{(k)}(m)) \in \Xi) = 0.
\]
Consequently, we have 
$P_{s,x}(\sigma_{s}(\Xi) < \infty) = 0,$ that is, $P_{s,x}(\sigma_{s}(\Xi) = \infty) = 1;$
hence the proposition is proved.
\hspace*{\fill}
$\Box$ 

\subsection{The continuity property of 
$\sigma_{0}(\overline{\Sigma}_{(\infty)})$} 
\label{subsect.continuity.hitting.time}

\noindent

Since $\overline{\Sigma}_{(\infty)} = \overline{\Sigma}_{1(\infty)} \cup \overline{\Sigma}_{2(\infty)},$ we have 
$
\sigma_{0}(\overline{\Sigma}_{(\infty)}) = \sigma_{0}(\overline{\Sigma}_{1(\infty)}) \wedge \sigma_{0}(\overline{\Sigma}_{2(\infty)}).
$
Therefore it is enough to examine the continuity property of $\sigma_{0}(\overline{\Sigma}_{1(\infty)})$ and  $\sigma_{0}(\overline{\Sigma}_{2(\infty)}).$ 

\medskip
\noindent
\textbf{(1) The continuity property of $\sigma_{0}(\overline{\Sigma}_{1(\infty)})$}
\medskip

For a given $s \geq 0$ and $\bm{\omega} \in \bm{U},$ suppose that $L(r, \bm{\omega}) = 0$ $(0 \leq r \leq s),$ and let $\tau_{s}(t) \equiv \tau_{s}(t, \bm{\omega}):= \inf \{u > s: L(u, \bm{\omega}) > t\}$ for $t \geq 0.$
Then we have 
\begin{enumerate}
\item[(1--1)] 
$\displaystyle{
\tau_{s}(t, \bm{\omega}) = \tau_{0}(t, \bm{\omega}). 
}$
\item[(1--2)]
$\tau_{s}(0, \bm{\omega}) = s ~\mbox{ iff } ~ L(u, \bm{\omega}) > s ~\mbox{ for every } u > s$ \\
$ ~~~~~~~~~~~~~~~~~\,  \mbox{ iff } ~ L(u, \bm{\omega}) > s ~\mbox{ for some } u \in (s, s+h] ~\mbox{ with every } h > 0.$
\item[(1--3)]
If~ 
$\displaystyle{
L(u,  \bm{\omega}) = \int_{s}^{u}1_{\Gamma}(X(r,  \bm{\omega}))L(dr,  \bm{\omega})}
$~
for every~$u \geq s,$ 
then the equality $\tau_{s}(0, \bm{\omega}) = s$ implies that $X(u, \bm{\omega}) \in \Gamma$ for some $u \in (s, s+h]$ with every $h > 0.$
\end{enumerate}
For each $s \geq 0$ let 
\begin{align*}
\bm{M}(t) \equiv \bm{M}_{s}(t):= X(t) - X(s) - \int_{s}^{t}\bm{c}(r, X(r))dr - \int_{s}^{t}\bm{\beta}(r, X(r))L(dr).
\end{align*}
Then $\{\bm{M}(t); t \geq s\}$ is an $n$--dimensional square integrable continuous martingale on the filtered probability space $(\bm{U}, \mathcal{U}, P_{s,x}; \mathcal{U}_{t}^{s})$ $( x \in \overline{\Omega})$ with the quadratic variational processes $\langle M_{i}, M_{j}\rangle$ $(i, j = 1, 2, \ldots ,n):$ 
\[
\langle M_{i}, M_{j}\rangle(t) = \int_{s}^{t}2a_{ij}(r, X(r))dr,
\]
where $M_{i}(t)$ is the $i$ th component of $\bm{M}(t)$ and $a_{ij}(r,z)$ is the $(i, j)$ entry of the diffusion matrix $A(r, z).$ 
Moreover let 
\[
\hat{X}(t) \equiv \hat{X}_{s}(t):= X(t) - \int_{s}^{t}\bm{\beta}(r, X(r))L(dr).
\]
Approximating the integral $\displaystyle{\int_{s}^{t}\bm{\beta}(r, X(r, \bm{\omega}))L(dr, \bm{\omega})}$ by a Riemann sum, we get the continuity of the map 
\[
\bm{U} \ni  \bm{\omega} \mapsto \int_{s}^{\cdot}\bm{\beta}(r, X(r, \bm{\omega}))L(dr, \bm{\omega}) \in C([s, \infty) \to \bm{R}^{n}),
\] 
where we suppose that $C([s, \infty) \to \bm{R}^{n})$ is equipped with the locally uniform convergence topology. 
Therefore the map $\bm{U} \ni \bm{\omega} \mapsto \hat{X}(\cdot, \bm{\omega}) \in C([s, \infty) \to \bm{R}^{n})$ is also continuous, and we see the following:
\begin{enumerate}
\item[(1--4)] For $s \geq 0,$ let $\sigma_{s+}(\overline{\Omega}^{\,c}; \hat{X}): = \inf \{t > s: \hat{X}(t) \in \overline{\Omega}^{\,c}\}.$ 
Then, for $x \in \Gamma,$  $P_{s,x}(\sigma_{s+}(\overline{\Omega}^{\,c}; \hat{X}) = s) = 1.$
\end{enumerate}
Indeed, letting $\bm{n} = \bm{n}(x)$ the outward unit normal vector to the boundary $\Gamma = \partial\Omega$ at $x,$ we see that $\{m(t): = \bm{M}(t)\iprod\bm{n}; t \geq s\}$ is a square integrable continuous martingale with the quadratic variational process 
\[
\langle m\rangle(t) = \int_{s}^{t}2A(r, X(r))\bm{n}\iprod\bm{n}\, dr.
\]
By the uniform ellipticity of the diffusion matrix $A,$ we have $\langle m\rangle(t) \approx t - s.$ \\
Therefore $\{B(t):= m(\langle m\rangle^{-1}(t)); t \geq 0\}$ is a one--dimensional Brownian motion 
on $\left(\bm{U}, \mathcal{U}, P_{s,x};\right.$ 
$\left.\mathcal{U}_{\langle m\rangle^{-1}(t)}^{s}\right).$ 
Hence the law of the iterated logarithm for the Brownian motion $\{B(t)\}$ yields the law for 
$\{m(t)\}.$ 
Since 
\[
\hat{X}(t)\iprod\bm{n} = X(s)\iprod\bm{n} + \int_{s}^{t}\bm{c}(r, X(r))\iprod\bm{n}\, dr + m(t),
\]
the behavior of  $\{\hat{X}(t)\iprod\bm{n}\}$ at time near the starting time $s$ is governed by that of 
$\{m(t)\};$ 
hence $\{\hat{X}(t)\iprod\bm{n}\}$ hits $\Omega$ and $\overline{\Omega}^{\,c}$ instantaneously.
Since $\Gamma$ is of class $C^{2,\alpha}$, if necessary, via a local flattening of the boundary, 
we therefore see that the process $\{\hat{X}(t)\}$ has also the same property as $\{\hat{X}(t)\iprod\bm{n}\}$ and (1--4) is verified. 

By (1--2), if $\tau_{s}(0, \bm{\omega}) > s,$ then there exists an $h = h(\bm{\omega}) > 0$ such that 
$L(t, \bm{\omega}) = 0$ for every $t \in [s, s+h].$ 
Therefore, if $P_{s,x}(\tau_{s}(0) > s) > 0,$ then it holds 
\[
P_{s,x}(X(t) = \hat{X}(t) ~\mbox{ for every } ~t \in (s, s+h] ~\mbox{ with some } ~ h > 0) > 0 .
\]
This contradicts the fact (1--4) and hence we get the following:
\begin{enumerate}
\item[(1--5)] For $s \geq 0, ~x \in \Gamma,$ we have $P_{s,x}(\tau_{s}(0) = s) = 1.$
\end{enumerate}
This implies the following fact:
\begin{enumerate}
\item[(1--6)] For $(s, x) \in \, \stackrel{\circ}{\Sigma}_{1(\infty)},$ $P_{s,x}(\sigma_{s+}(\stackrel{\circ}{\Sigma}_{1(\infty)}) = s) = 1.$
\end{enumerate}
Indeed, since $(s, x) \in \, \stackrel{\circ}{\Sigma}_{1(\infty)} \subset [0, \infty) \times \Gamma,$ using the fact (1--5), we have 
\[
P_{s,x}(X(t) \in \Gamma ~\mbox{ with some } ~t \in (s, s+h]  ~\mbox{ for every } ~ h > 0) = 1.
\]
On the other hand, $P_{s,x}((s, X(s)) = (s, x) \in \, \stackrel{\circ}{\Sigma}_{1(\infty)}) = 1$ and the set 
$\stackrel{\circ}{\Sigma}_{1(\infty)}$ is open in $\partial_{L}D_{(\infty)}$ and open in $[0, \infty) \times \Gamma.$ 
Therefore, when $t$ is sufficiently close to $s$ and $X(t, \bm{\omega}) \in \Gamma,$ it holds that $(t, X(t, \bm{\omega})) \in \, \stackrel{\circ}{\Sigma}_{1(\infty)}.$ 
This means $s \leq \sigma_{s+}(\stackrel{\circ}{\Sigma}_{1(\infty)})(\bm{\omega}) \leq t;$ that is, 
the fact (1--6) is verified.

Since  $P_{s,x}\left(\sigma_{s}(\overline{\Sigma}_{1(\infty)}) = \sigma_{s}(\stackrel{\circ}{\Sigma}_{1(\infty)})\right) = 1,$
the strong Markov property yields
\begin{enumerate}
\item[(1--7)] For $s \geq 0, ~x \in D(s) \cup \Gamma^{\prime},$ 
\[
\hspace*{-2cm}P_{s,x}(\tau_{\sigma_{s}(\overline{\Sigma}_{1(\infty)})}(0) = \sigma_{s}(\overline{\Sigma}_{1(\infty)}); \sigma_{s}(\overline{\Sigma}_{1(\infty)}) < \infty) = P_{s,x}(\sigma_{s}(\overline{\Sigma}_{1(\infty)}) < \infty).
\]
\end{enumerate}

Thus we state the continuity property of $\sigma_{0}(\overline{\Sigma}_{1(\infty)})$ in the following.

{\lemma
\label{continuity.1}
Let
\begin{align*}
&\widetilde{\bm{U}}_{1}:= \biggl\{\bm{\omega} \in \bm{U}: \tau_{\sigma_{0}(\overline{\Sigma}_{1(\infty)})(\bm{\omega})}(0,\bm{\omega}) = \sigma_{0}(\overline{\Sigma}_{1(\infty)})(\bm{\omega}), ~\sigma_{0}(\delta\Sigma_{1(\infty)})(\bm{\omega}) = \infty, \\
&\hspace*{5.5cm}L(t, \bm{\omega}) = \int_{0}^{t}1_{\Gamma}(X(r, \bm{\omega}))L(dr, \bm{\omega}) ~\mbox{ for } ~t \geq 0\biggr\},
\end{align*}
provided $\tau_{\infty}(0) = \infty.$
Then, for $s \geq 0$ and $x \in D(s) \cup \Gamma',$ 
$P_{s,x}\bigl(\widetilde{\bm{U}}_{1}\bigr) = 1$ and 
$\left.\sigma_{0}(\overline{\Sigma}_{1(\infty)})\right|_{\widetilde{\bm{U}}_{1}}$ is continuous.
}
\medskip

\noindent
\textbf{Proof.}
The fact $P_{s,x}\bigl(\widetilde{\bm{U}}_{1}\bigr) = 1$ is derived from the facts (\ref{s.to.0}), (1--1), (1--7) and  
\begin{align*}
&P_{s,x}\biggl(\sigma_{s}(\overline{\Sigma}_{1(\infty)}) = \sigma_{0}(\overline{\Sigma}_{1(\infty)}), ~\sigma_{s}(\delta\Sigma_{1(\infty)}) = \sigma_{0}(\delta\Sigma_{1(\infty)}), \\
&\hspace*{3.5cm}L(t) = \int_{0}^{t}1_{\Gamma}(X(r))L(dr) ~\mbox{ for } ~t \geq 0\biggr) = 1.
\end{align*}
Take $\bm{\omega}_{0} \in \widetilde{\bm{U}}_{1}$ with $\sigma_{0}(\overline{\Sigma}_{1(\infty)})(\bm{\omega}_{0}) = \infty.$ 
For $R > 0,$ letting  
\[
G_{[0, R]}(\bm{\omega}_{0}):= \{(t, X(t, \bm{\omega}_{0})): 0 \leq t \leq R\},
\]
we have $d(G_{[0, R]}(\bm{\omega}_{0}), ~\overline{\Sigma}_{1(\infty)}) > 0$ for every $R > 0.$ 
Therefore, when $\bm{\omega} \to \bm{\omega}_{0}$ in $\bm{U},$ we get $\sigma_{0}(\overline{\Sigma}_{1(\infty)})(\bm{\omega}) \to \infty;$ that is,
\[
\lim_{\bm{\omega} \to \bm{\omega}_{0}}\sigma_{0}(\overline{\Sigma}_{1(\infty)})(\bm{\omega}) = \infty = \sigma_{0}(\overline{\Sigma}_{1(\infty)})(\bm{\omega}_{0}).
\]

In the sequel, for an arbitrarily given $\bm{\omega}_{0} \in \widetilde{\bm{U}}_{1}$ with $\sigma_{0}(\overline{\Sigma}_{1(\infty)})(\bm{\omega}_{0}) < \infty,$ we examine the continuity property of $\sigma_{0}(\overline{\Sigma}_{1(\infty)})$ at $\bm{\omega}_{0}.$ 
 In the case where \\ $\sigma:= \liminf_{\bm{\omega} \to \bm{\omega}_{0}}\sigma_{0}(\overline{\Sigma}_{1(\infty)})(\bm{\omega}) < \infty,$ we take a sequence $\left\{\bm{\omega}_{m}\right\}_{m=1}^{\infty} \subset \bm{U}$ such that $\bm{\omega}_{m} \to \bm{\omega}_{0}$ and $\sigma_{0}(\overline{\Sigma}_{1(\infty)})(\bm{\omega}_{m}) \to \sigma$ as $m \to \infty.$ 
Then there exists a decreasing sequence $\left\{\varepsilon_{m}\right\}_{m=1}^{\infty}$ of positive numbers such that $\varepsilon_{m} \downarrow 0$ and 
\begin{align*}
&\left(\sigma_{0}(\overline{\Sigma}_{1(\infty)})(\bm{\omega}_{m}) + \varepsilon_{m}, X(\sigma_{0}(\overline{\Sigma}_{1(\infty)})(\bm{\omega}_{m}) + \varepsilon_{m}, \bm{\omega}_{m})\right) (\in \overline{\Sigma}_{1(\infty)}) \\
&\longrightarrow 
(\sigma, X(\sigma, \bm{\omega}_{0})) (\in \overline{\Sigma}_{1(\infty)}) ~~~\mbox{ as } ~
m \to \infty.
\end{align*}
Hence 
\begin{equation}
\sigma = \liminf_{\bm{\omega} \to \bm{\omega}_{0}}\sigma_{0}(\overline{\Sigma}_{1(\infty)})(\bm{\omega}) \geq \sigma_{0}(\overline{\Sigma}_{1(\infty)})(\bm{\omega}_{0}).
\label{liminf}
\end{equation}
The inequality (\ref{liminf}) is also holds whenever $\sigma = \infty.$ 
This shows 
\[
\widetilde{\bm{U}}_{1} \subset \{\bm{\omega} \in \bm{U}: \sigma_{0}(\overline{\Sigma}_{1(\infty)}) ~\mbox{ is lower semicontinuous at } ~\bm{\omega}\}.
\]
Next let 
\[
(\sigma_{0}, x_{0}): = (\sigma_{0}(\overline{\Sigma}_{1(\infty)})(\bm{\omega}_{0}), X(\sigma_{0}(\overline{\Sigma}_{1(\infty)})(\bm{\omega}_{0}), \bm{\omega}_{0})).
\]
Then 
\[
(\sigma_{0}, x_{0}) = (\sigma_{0}(\stackrel{\circ}{\Sigma}_{1(\infty)})(\bm{\omega}_{0}), X(\sigma_{0}(\stackrel{\circ}{\Sigma}_{1(\infty)})(\bm{\omega}_{0}), \bm{\omega}_{0})) \in \, \stackrel{\circ}{\Sigma}_{1(\infty)}.
\]
Therefore, noting that $\stackrel{\circ}{\Sigma}_{1(\infty)}$ is open in $[0, \infty) \times \Gamma,$                we can provide an interval $[a, b)$ and an open ball $B_{\delta}(x_{0})$ in $\bm{R}^{n}$ such that $a \leq \sigma_{0} < b$ and 
\[
\left([a, b) \times B_{\delta}(x_{0})\right) \cap ([0, \infty) \times \Gamma) \subset \, \stackrel{\circ}{\Sigma}_{1(\infty)};
\]
and for given $\delta$ there exists an $h_{0} > 0$ with $\sigma_{0} + h_{0} \leq b$ and 
\[
X(t, \bm{\omega}_{0}) \in B_{\delta/2}(x_{0}) ~~\mbox{ for every } ~ t \in [\sigma_{0}, \sigma_{0} + h_{0}). 
\]
Since $L(\sigma_{0}, \bm{\omega}_{0}) < L(\sigma_{0} + h, \bm{\omega}_{0})$ for every $h > 0$ and since $L(\sigma_{0}, \bm{\omega})$ and $L(\sigma_{0} + h, \bm{\omega})$ are continuous in 
$\bm{\omega} \in \bm{U},$ it holds for an arbitrarily fixed $h \in (0, h_{0}]$ 
\[
L(\sigma_{0}, \bm{\omega}) < L(\sigma_{0} + h, \bm{\omega}),
\]
provided $\bm{\omega} \in \widetilde{\bm{U}}_{1}$ is sufficiently close to $\bm{\omega}_{0}$ in 
$\bm{U}.$ 
In particular, for any $R > 0,$ we can take $\delta^{\prime} \in (0, \delta/2)$ such that for every $\bm{\omega} \in \widetilde{\bm{U}}_{1}$ with $\|\bm{\omega} - \bm{\omega}_{0}\|_{[0,R]} < \delta^{\prime}$ it holds 
$
L(\sigma_{0}, \bm{\omega}) < L(\sigma_{0} + h, \bm{\omega}) .
$
For such an $\bm{\omega}$, noting 
\[
L(\sigma_{0} + h, \bm{\omega}) - L(\sigma_{0}, \bm{\omega}) = \int_{\sigma_{0}}^{\sigma_{0}+h}
1_{\Gamma}(X(r, \bm{\omega})) L(dr) > 0,
\]
we can choose $t^{\prime} \in [\sigma_{0}, \sigma_{0} + h]$ with $X(t^{\prime}, \bm{\omega}) \in \Gamma.$ 
This implies 
\[
(t^{\prime}, X(t^{\prime}, \bm{\omega})) \in ([a, b) \times B_{\delta}(x_{0})) \cap ([0, \infty) \times \Gamma) \subset \, \stackrel{\circ}{\Sigma}_{1(\infty)} \, \subset \overline{\Sigma}_{1(\infty)}; 
\]
hence  
$
\sigma_{0}(\overline{\Sigma}_{1(\infty)})(\bm{\omega}) \leq t^{\prime} \leq \sigma_{0}+h.
$
Therefore
\[
\limsup_{\bm{\omega} \to \bm{\omega}_{0}, \bm{\omega} \in \widetilde{\bm{U}}_{1}}\sigma_{0}(\overline{\Sigma}_{1(\infty)})(\bm{\omega}) \leq \sigma_{0} = \sigma_{0}(\overline{\Sigma}_{1(\infty)})(\bm{\omega}_{0}). 
\]
Thus we conclude that $\left.\sigma_{0}(\overline{\Sigma}_{1(\infty)})\right|_{\widetilde{\bm{U}}_{1}}$ is continuous.    
\hspace*{\fill}$\Box$

\medskip 
\noindent
\textbf{
(2) The continuity property of $\sigma_{0}(\overline{\Sigma}_{2(\infty)})$}
\medskip 

In the sequel, we assume that the coefficients $A$ and $\bm{c}$ of $\mathcal{L}_{0}$ are extended boundedly onto $[0, \infty) \times \bm{R}^{n}$ with the uniform ellipticity as in Assumption \ref{assump_2}(i) and with the same regularity as in Assumption \ref{assump_3}(iii). 
We begin with preparing necessary facts.
\begin{enumerate}
\item[(2--1)]  For $(s, x) \in \Sigma_{2(\infty)},$  $P_{s,x}\left(\sigma_{s+}(((0, \infty) \times \Omega)\setminus \overline{D_{(\infty)}}) = s\right) =1.$ 
\end{enumerate}
Indeed, for $x \in \Sigma_{2(\infty)}(s) \subset \Omega,$ there exists a neighborhood $G$ of $x$ with $G \subset \Omega.$ 
Because of the uniqueness of solutions to the $\mathcal{P}_{0}$--martingale problem on 
$\bm{R}^{n}$, we get 
\[
\left.P_{s,x}\right|_{\mathcal{U}^{s}_{\sigma_{s}(\partial G)}} = \left.\bar{P}_{s,x}\right|_{\mathcal{U}^{s}_{\sigma_{s}(\partial G)}},
\]
where $\bar{P}_{s,x}$ denotes the solution to the $\mathcal{P}_{0}$--martingale problem on 
$\bm{R}^{n}$ starting at $(s, x).$ 
In particular, 
\begin{align*}
P_{s,x}\left(X(t) \in dy; \sigma_{s}(\partial G) > t\right) = \bar{P}_{s,x}\left(X(t) \in dy; \sigma_{s}(\partial G) > t\right) =:p^{G}(s, x; t, y)dy 
\end{align*}
and hence 
\[
 p(s, x; t, y) \geq  p^{G}(s, x; t, y) ~~~\mbox{ for } ~0 \leq s < t ~\mbox{ and } ~ x, y \in G. 
 \]
Using the upper and lower Gaussian bounds for the density $\bar{p}(s, x; t, y)$ of the transition probability $\bar{P}_{s,x}\left(X(t) \in dy\right)$ in \cite{Aro67} (see also \cite{Aro68}) and the localization argument in \cite{Str08} (see Lemma 5.2.10 or Theorem 5.2.11), we can obtain a lower Gaussian bound for $p^{G}(s, x; t, y)$ in a sufficiently small neighborhood $G'$ of $x$ with $G' \Subset G.$ 
Then we use the argument in the proof of Lemma 2.3 in \cite{CGK06}.  
For $R > 0, \delta > 0$ and $\bar{x} \in \bm{R}^{n},$ set  
\[
\mathcal{T}_{(s,x)}^{R,\delta,\bar{x}}:= \left\{(t, y): s < t < s + \delta, ~ \left|y - x - \bar{x}\sqrt{t - s}\right|^{2} < R^{2}(t - s)\right\}.
\]
It follows from Condition \ref{cond_D_1} for $D$ that $D$ satisfies the following tusk condition: 
for some $R > 0, \delta > 0$ and $\bar{x} \in \bm{R}^{n}$ there exists a $\mathcal{T}_{(s,x)}^{R,\delta,\bar{x}}$ with  
$\overline{\mathcal{T}_{(s,x)}^{R,\delta,\bar{x}}} \cap \overline{D} = \{(s, x)\}.$ 
We take $R, \delta$ so small that 
$\displaystyle{\bigcup_{0 < h < \delta}\,
\mathcal{T}_{(s,x)}^{R,\delta,\bar{x}}(s+h) \subset G'},$
where $\mathcal{T}_{(s,x)}^{R,\delta,\bar{x}}(s+h)$ is the $(s+h)$-section of 
$\mathcal{T}_{(s,x)}^{R,\delta,\bar{x}}.$   
Then 
\begin{align*}
p_{s,x;h}&:= P_{s,x}\left((t, X(t)) \notin \overline{D}  ~\mbox{ for some } t \in (s, s + h]\right) \\
&\geq P_{s,x}\left((s, X(s + h)) \in \mathcal{T}_{(s,x)}^{R,\delta,\bar{x}}\right) \\
&\geq \int_{\mathcal{T}_{(s,x)}^{R,\delta,\bar{x}}(s+h)}p^{G}(s, x; s+h, y) dy \\
&\geq \int_{|z - \bar{x}|^{2} < R}\frac{1}{K}e^{-K|z|^{2}} dz 
\end{align*}
for some constant $K > 0$ independent of $h > 0.$ 
This implies that $\displaystyle{\inf_{0<h<\delta}p_{s,x;h} > 0}$ and hence, by Blumenthal's $0$--$1$ law, we obtain the result (2--1). 
\begin{enumerate}
\item[(2--2)]  For $s \geq 0$ and $x \in D(s) \cup \Gamma^{\prime},$ 
\begin{align*}
P_{s,x}\left(\sigma_{s}(\overline{\Sigma}_{2(\infty)}) = \sigma_{s+}(((0, \infty) \times \Omega)\setminus \overline{D_{(\infty)}})\right) =1.
\end{align*} 
\end{enumerate}
Indeed, since $\sigma_{s}(\Sigma_{2(\infty)}) \leq \sigma_{s+}(((0, \infty) \times \Omega)\setminus \overline{D_{(\infty)}})$ on $\{X(s) = x\}$,   $\sigma_{s}(\overline{\Sigma}_{2(\infty)}) = \sigma_{s}(\Sigma_{2(\infty)})$ on $\{\sigma_{s}(\delta\Sigma_{2(\infty)}) = \infty\}$ and $P_{s,x}\left(\sigma_{s}(\delta\Sigma_{2(\infty)}) = \infty\right) = 1,$ we have 
\begin{align*}
P_{s,x}\left(\sigma_{s}(\overline{\Sigma}_{2(\infty)}) = \sigma_{s+}(((0, \infty) \times \Omega)\setminus \overline{D_{(\infty)}}); \sigma_{s}(\overline{\Sigma}_{2(\infty)}) = \infty\right) = P_{s,x}\left(\sigma_{s}(\overline{\Sigma}_{2(\infty)}) = \infty\right).
\end{align*}
Moreover, it holds 
\begin{align*}
&P_{s,x}\left((\sigma_{s}(\overline{\Sigma}_{2(\infty)}), X(\sigma_{s}(\overline{\Sigma}_{2(\infty)})) \in \Sigma_{2(\infty)}; \sigma_{s}(\overline{\Sigma}_{2(\infty)}) < \infty\right) = P_{s,x}\left(\sigma_{s}(\overline{\Sigma}_{2(\infty)}) < \infty\right), \\
&P_{s,x}\left(\sigma_{s+}(((0, \infty) \times \Omega)\setminus \overline{D_{(\infty)}}) = \sigma_{\sigma_{s}(\overline{\Sigma}_{2(\infty)})+}(((0, \infty) \times \Omega)\setminus \overline{D_{(\infty)}}); \sigma_{s}(\overline{\Sigma}_{2(\infty)}) < \infty\right) \\ 
&\hspace*{5.5cm}= P_{s,x}\left(\sigma_{s}(\overline{\Sigma}_{2(\infty)}) < \infty\right). 
\end{align*} 
Hence, by (2--1) and the strong Markov property, we have (2--2). 

Therefore the continuity property of $\sigma_{0}(\overline{\Sigma}_{2(\infty)})$ is stated as follows: 
{\lemma
\label{continuity.2}
Let 
\begin{align*}
&\widetilde{\bm{U}}_{2}:= \biggl\{\bm{\omega} \in \bm{U}: \sigma_{0}(\overline{\Sigma}_{2(\infty)})(\bm{\omega}) = \sigma_{0+}(((0, \infty) \times \Omega)\setminus \overline{D_{(\infty)}})(\bm{\omega}), \\ 
&\hspace*{8cm}\sigma_{0}(\delta\Sigma_{2(\infty)})(\bm{\omega}) = \infty\biggr\}. 
\end{align*}
Then, for $s \geq 0$ and $x \in D(s) \cup \Gamma^{\prime},$ ~$P_{s,x}(\widetilde{\bm{U}}_{2}) = 1$ and 
\[
\widetilde{\bm{U}}_{2} \subset \{\bm{\omega} \in \bm{U}: \sigma_{0}(\overline{\Sigma}_{2(\infty)}) ~\mbox{ is continuous at}\, ~\bm{\omega}\}. 
\]
}
\medskip

\noindent
\textbf{Proof.} 
From the previous argument, we have $P_{s,x}(\widetilde{\bm{U}}_{2}) = 1,$ and by the same argument as in the proof of Lemma \ref{continuity.1}, we see that 
\[
\widetilde{\bm{U}}_{2} \subset \{\bm{\omega} \in \bm{U}: \sigma_{0}(\overline{\Sigma}_{2(\infty)}) ~\mbox{ is lower semicontinuous at}\, ~\bm{\omega}\}.
\]
On the other hand, since $((0, \infty) \times \Omega)\setminus \overline{D_{(\infty)}}$ is open, it is easy to see that 
\begin{align*}
\widetilde{\bm{U}}_{2} \subset \{\bm{\omega} \in \bm{U}: \sigma_{0+}(((0, \infty) \times \Omega)\setminus \overline{D_{(\infty)}}) ~\mbox{ is upper semicontinuous at}\, ~\bm{\omega}\}; 
\end{align*}
hence the lemma is proved.  
\hspace*{\fill}$\Box$

\medskip 

Consequently, we obtain the continuity property for $\sigma_{0}(\overline{\Sigma}_{(\infty)}):$ 

{\theorem
\label{continuity.final}
Let $\widetilde{\bm{U}} := \widetilde{\bm{U}}_{1} \cap \widetilde{\bm{U}}_{2}.$
Then, for $s \geq 0$ and $x \in D(s) \cup \Gamma^{\prime},$ ~$P_{s,x}(\widetilde{\bm{U}}) = 1$ and $\left.\sigma_{0}(\overline{\Sigma}_{(\infty)})\right|_{\widetilde{\bm{U}}}$ is continuous. 
}

\subsection{Proof of (\ref{cont.Dirichlet_{1}})} 
\label{subsect.continuity.Dirichlet.part_{1}}

We first show the following. 
\begin{enumerate}
\item[(3--1)] 
For arbitrarily given $\epsilon > 0,$ there exists $\delta_{1} > 0$ such that 
\[
P_{s,x}(\sigma_{s}(\Gamma) < s + 2\eta) > 1 - \epsilon 
\]
for every $(s, x) \in [0, \infty) \times \overline{\Omega}$ with $d((s, x), (s_{0}, x_{0})) < \delta_{1}.$ 
\end{enumerate}
Indeed, It follows from the argument in the proof of  the continuity of $\sigma_{0}(\overline{\Sigma}_{1(\infty)})$ (see (1--5) and (1--7) in particular) that 
\[
P_{s,x}\left(\tau_{s}(0) = \tau_{\sigma_{s}(\Gamma)}(0) = \sigma_{s}(\Gamma)\right) = 1
\]
for $(s, x) \in [0, \infty) \times \overline{\Omega}.$
Noting (1--2) in the proof of  the continuity of $\sigma_{0}(\overline{\Sigma}_{1(\infty)}),$ we have 
\begin{align*}
1 &= P_{s_{0}, x_{0}}(\tau_{s_{0}}(0) = s_{0}) \\
&= P_{s_{0}, x_{0}}(L(r) > 0~~\mbox{for some}~ r \in (s_{0}/2, s_{0} + \eta] \cap \bm{Q})
\end{align*}
for any $\eta > 0.$ 
Then we use the same argument as in the proof of Theorem 13.1 of \cite{Wen81} in the following. 
Let $(s_{0}/2, s_{0} + \eta] \cap \bm{Q} = \{r_{1}, r_{2}, \dots\}.$ 
For sufficiently large $m$, it holds that 
\[
P_{s_{0}, x_{0}}(\{L(r_{1}) > 0\} \cup \{L(r_{2}) > 0\} \cup \cdots \cup \{L(r_{m}) > 0\})  \geq 1 - \frac{\epsilon}{3}, 
\]
that is, 
\[
P_{s_{0}, x_{0}}(L(r_{1}) = 0, L(r_{2}) = 0, \ldots, L(r_{m}) = 0)  < \frac{\epsilon}{3}.
\] 
Since $P_{s,x} \longrightarrow P_{s_{0}, x_{0}}$ weakly as $(s, x) \to (s_{0}, x_{0})$ and 
$\{L(r_{1}) = 0, L(r_{2}) = 0, \dots, L(r_{m}) = 0\}$ is closed in $\bm{U}$, we have  
\begin{align*}
&\limsup_{(s, x) \to (s_{0}, x_{0})} P_{s, x}(L(r_{1}) = 0, L(r_{2}) = 0, \dots, L(r_{m}) = 0) \\ 
&\leq P_{s_{0}, x_{0}}(L(r_{1}) = 0, L(r_{2}) = 0, \dots, L(r_{m}) = 0).
\end{align*}
Therefore we find a $\delta_{1} \in (0, \eta)$ for which 
\[
P_{s,x}(L(r_{1}) = 0, \dots, L(r_{m}) = 0) < \epsilon
\]
for every $(s, x) \in [0, \infty) \times \overline{\Omega}$ with $d((s, x), (s_{0}, x_{0})) < \delta_{1},$
here we take $\delta_{1}$ satisfying $0 < \delta_{1} < s_{0}/2$ in the case $s_{0} > 0.$   
Then for such an $(s, x)$, 
\[
P_{s,x}(\tau_{s}(0) \geq s + 2\eta) < \epsilon.
\]
Noting $P_{s,x}\left(\tau_{s}(0) = \sigma_{s}(\Gamma)\right) =1,$  we see that 
\[
P_{s,x}\left(\sigma_{s}(\Gamma) < s + 2\eta\right) > 1 - \epsilon
\]
for every $(s, x) \in [0, \infty) \times\overline{\Omega}$ with $d((s, x), (s_{0}, x_{0})) < \delta_{1}.$
This proves (3--1). \\
Next we show the following.
\begin{enumerate}
\item[(3--2)]
For arbitrarily given $\epsilon > 0,$ there exists $\delta_{2} \in (0, \eta)$ such that 
\[
P_{s,x}(\sigma_{s}(\overline{\Gamma'}) < s + 2\eta) \leq  P_{s_{0},x_{0}}(\sigma_{s_{0}}(\overline{\Gamma'}) < s_{0} + 4\eta) + \epsilon 
\]
for every $(s, x) \in [0, \infty) \times (\Omega \cup \Gamma^{\prime\prime})$ with $d((s, x), (s_{0}, x_{0})) < \delta_{2}.$ 
\end{enumerate}
For given $\eta$, take a piecewise linear function $\lambda_{\eta}(s)$ satisfying 
\[
\lambda_{\eta}(s) = 0 ~~\mbox{for}~ 0 \leq s \leq s_{0} + 3\eta,~~~\lambda_{\eta}(s) = 1 ~~\mbox{for}~ s \geq s_{0} + 4\eta.
\]
Then 
\[
1_{[s_{0}+3\eta, \infty)} \geq \lambda_{\eta}(s) \geq 1_{[s_{0}+4\eta, \infty)}.
\]  
In the same way as in the proof of  
the continuity of $\sigma_{0}(\overline{\Sigma}_{1(\infty)})$,   
we see that $\sigma_{0}(\overline{\Gamma'})$ has an analogous continuity property like as Lemma \ref{continuity.1} with respect to the probability measures $P_{s,x}$ for $(s, x) \in [0, \infty) \times (\Omega \cup \Gamma^{\prime\prime}).$ 
Therefore, using Theorem \ref{mapping.theorem} below, we ensure that $E_{s,x}[\lambda_{\eta}(\sigma_{s}(\overline{\Gamma'}))] = E_{s,x}[\lambda_{\eta}(\sigma_{0}(\overline{\Gamma'}))]$ is continuous in $(s, x) \in [0, \infty) \times (\Omega \cup \Gamma^{\prime\prime}).$ 
Hence for given $\epsilon > 0$ there exists $\delta_{2} \in (0, \eta)$ such that 
\[
E_{s,x}[\lambda_{\eta}(\sigma_{s}(\overline{\Gamma'}))] > E_{s_{0},x_{0}}[\lambda_{\eta}(\sigma_{s_{0}}(\overline{\Gamma'}))] - \epsilon
\] 
for every $(s, x) \in [0, \infty) \times (\Omega \cup \Gamma^{\prime\prime})$ with $d((s, x), (s_{0}, x_{0})) < \delta_{2}.$ 
Accordingly, for such an $(s, x)$, we see that $s + 2\eta < s_{0}+3\eta$ and 
\[
P_{s,x}(\sigma_{s}(\overline{\Gamma'}) \geq s_{0}+3\eta) \geq P_{s_{0},x_{0}}(\sigma_{s_{0}}(\overline{\Gamma'}) \geq s_{0}+4\eta) - \epsilon.
\]
These imply 
\[
P_{s,x}(\sigma_{s}(\overline{\Gamma'}) < s + 2\eta) \leq P_{s_{0},x_{0}}(\sigma_{s_{0}}(\overline{\Gamma'}) < s_{0}+4\eta) + \epsilon; 
\]
hence (3--2) is verified. 
\begin{enumerate}
\item[(3--3)] 
For given $\epsilon > 0,$ put $\delta: = \delta_{1} \wedge \delta_{2}.$ Then 
\[
P_{s,x}(\sigma_{s}(\overline{\Sigma}_{(\infty)}) \geq s + 2\eta) < P_{s_{0},x_{0}}(\sigma_{s_{0}}(\overline{\Gamma'}) < s_{0}+4\eta) + 2\epsilon
\]
for $(s, x) \in D$ with $d((s, x), (s_{0}, x_{0})) < \delta.$ 
\end{enumerate}
Indeed, let 
\[
\Theta:= ([0, \infty) \times \Gamma) \setminus \{\Sigma_{1(\infty)} \cup ([0, \infty) \times \overline{\Gamma'})\}.
\]
By (3--1), we see that for $(s, x) \in D$ with $d((s, x), (s_{0}, x_{0})) < \delta$ 
\begin{align}
1 - \epsilon &< P_{s,x}(\sigma_{s}(\Gamma) < s + 2\eta) \nonumber \\
&= P_{s,x}(\sigma_{s}(\Gamma) < s + 2\eta; \sigma_{s}(\Gamma) = \sigma_{s}(\Theta)) \nonumber \\
&+ P_{s,x}(\sigma_{s}(\Gamma) < s + 2\eta; \sigma_{s}(\Gamma) = \sigma_{s}(\Sigma_{1(\infty)})) \nonumber \\
&+ P_{s,x}(\sigma_{s}(\Gamma) < s + 2\eta; \sigma_{s}(\Gamma) = \sigma_{s}(\overline{\Gamma'})). 
\label{ineq.ent.gamma}
\end{align}
Noting the relationships:   
\begin{align*}
&P_{s,x}(\sigma_{s}(\Theta) \geq \sigma_{s}(\overline{\Sigma}_{(\infty)})) = 1, \\
&\sigma_{s}(\stackrel{\circ}{\Sigma}_{1(\infty)}) \geq\sigma_{s}(\Sigma_{1(\infty)}) \geq \sigma_{s}(\overline{\Sigma}_{(\infty)}), ~ P_{s,x}(\sigma_{s}(\stackrel{\circ}{\Sigma}_{1(\infty)}) = \sigma_{s}(\overline{\Sigma}_{1(\infty)})) = 1,\\
&P_{s,x}(\sigma_{s}(\Gamma) < s + 2\eta; \sigma_{s}(\Gamma) = \sigma_{s}(\overline{\Gamma'})) \leq P_{s_{0},x_{0}}(\sigma_{s_{0}}(\overline{\Gamma'}) < s_{0}+4\eta) + \epsilon
\end{align*}
for $(s, x) \in D$ with $d((s, x), (s_{0}, x_{0})) < \delta$, 
we ensure that the right hand side of (\ref{ineq.ent.gamma}) is dominated by 
\[
P_{s,x}(\sigma_{s}(\overline{\Sigma}_{(\infty)}) < s + 2\eta) + P_{s_{0},x_{0}}(\sigma_{s_{0}}(\overline{\Gamma'}) < s_{0}+4\eta) + \epsilon.
\]
Hence 
\[
P_{s,x}(\sigma_{s}(\overline{\Sigma}_{(\infty)}) \geq s + 2\eta) < P_{s_{0},x_{0}}(\sigma_{s_{0}}(\overline{\Gamma'}) < s_{0}+4\eta) + 2\epsilon;
\]
that is, (3--3) is verified and the proof of (\ref{cont.Dirichlet_{1}}) is complete.


\subsection{Mapping theorem} 
\label{subsect.mapping.theorem}

\noindent
We need the following variation of the usual mapping theorem with respect to weak convergence of probability measures (see \cite{Top67} and also \cite{Bil99, Bil03} for the usual one).
Using Skorohod's theorem (see Theorem 6.7 in \cite{Bil99}) as in the proof of Theorem 25.7 in \cite{Bil03}, we have 
{\theorem
\label{mapping.theorem}
Let $S$ be a polish space with Borel field $\mathcal{S} \equiv \mathcal{B}(S)$ and $S'$ a separable metric space with Borel field $\mathcal{S}' \equiv \mathcal{B}(S').$ 
Consider a family of probability measures $\{\mu_{n}; n= 0, 1, \ldots\}$ on $(S, \mathcal{S})$ with $\mu_{n} \to \mu_{0}$ weakly as $n \to \infty.$ 
Suppose that $h: S \rightarrow S'$ is a measurable map such that there exists a set $\widetilde{S} \in \mathcal{S}$ with the properties: $\mu_{n}(\widetilde{S}) = 1$ for $n= 0, 1, \ldots$ and the restriction $\left.h\right|_{\widetilde{S}}$ is continuous. 
Then $\mu_{n}h^{-1} \to \mu_{0}h^{-1}$ weakly as $n \to \infty.$ 
Here  $\mu_{n}h^{-1}(dy):= \mu_{n}(h^{-1}(dy)).$ 
}

\section{Shape identification inverse problem in a Bayesian framework}  
\label{sect.appendix_B}
(by Hajime Kawakami)

\medskip
\medskip

\noindent
We consider the shape identification inverse problem 
mentioned in \S \ref{sect.appli} in a Bayesian framework based on 
\cite{Das13} and \cite{Stu10}.
In the following, Theorem \ref{continuity.func} plays 
an essentially important role.
Our inverse problem is an  
identification problem of the unknown shape of the Dirichlet boundary $\Sigma$  
from the thermal data $u^D(t,x)$ on $\Gamma^\omega,$ 
where $\Omega$ is fixed and 
$u^D$ is the solution to the parabolic 
problem (\ref{eq_1}) with $\mbox{\boldmath$f$} = \mbox{\boldmath$0$},$ 
$f = 0,$ and $h = 0.$  
The authors proved the uniqueness of this inverse problem 
in \cite{Kaw10}. 

In the Bayesian framework, it is often assumed that 
the observational data is a set of finite dimensional vectors, 
and the posterior distribution of unknown shapes is considered 
under the observational data.   
Although this setting is underdetermined 
as an inverse problem, we can obtain a stability result 
of the posterior distribution as described below.
Let $\mathcal{D}$ be the set of all possible domains 
of our inverse problem. 
Let $G_1$ be a map defined by 
\[
G_1: \mathcal{D} \longrightarrow L^2([0,T] \times \Gamma^\omega),
\hspace{0.5cm}
D \longmapsto \gamma_{\Gamma^\omega}u^D,
\]
where $\gamma_{\Gamma^\omega}$ is the trace operator. 
Let $m$ be a positive integer and 
$G_2$ be a suitable discretization map 
\[
G_2: L^2([0,T] \times \Gamma^\omega) \longrightarrow \mbox{\boldmath$R$}^m, 
\]
where $G_2$ is linear and continuous. 
Then, we define a map $F$ by $F := G_2 \circ G_1$ and 
call it a \textit{forward operator}.
We consider that $\mathcal{D}$ is a metric space 
equipped with Hausdorff metric. 
Then, Theorem \ref{continuity.func} ensures that $F$ is continuous. 
We assume that, for given domain $D \in \mathcal{D},$
we can observe data $y \in \mbox{\boldmath$R$}^m$ 
with noise $\eta \in \mbox{\boldmath$R$}^m,$  
\[
y = F(D) + \eta,  
\]
where $\eta$ is an $m$--dimensional Gaussian random variable. 
Many researches considered the case that $\eta$ 
follows an $m$--dimensional Gaussian distribution with some covariance matrix 
or considered more general situations (e.g. \cite{Das13}). 
For simplicity, we assume that 
the elements $\{\eta_i:\ i = 1, 2, \ldots, m\}$ 
of $\eta = (\eta_i)$ are mutually independent and 
$\eta_i \sim \mathcal{N}(0, \sigma^2),$ 
that is, each $\eta_i$ follows the Gaussian distribution with mean $0$ 
and variance $\sigma^2.$ 

Let $\mathcal{B}$ be the topological $\sigma$-field on $\mathcal{D}.$ 
Then, $F$ is a measurable map 
on $(\mathcal{D}, \mathcal{B}).$ 
Let $\mu_0$ be a probability measure on $(\mathcal{D}, \mathcal{B})$
and denote by $\|\ \cdot\ \|$ the standard Euclidean norm.
Define 
\[
\Psi(D;y) := 
\frac{1}{(2\pi\sigma^2)^{m/2}}
\exp\left(-\frac{1}{2\sigma^2}\|y - F(D)\|^2\right)  
\hspace{1cm}
\left(D \in \mathcal{D}, 
y \in \mbox{\boldmath$R$}^m\right)
\]
and
\[
Z_\Psi(y) 
:= \int_{\mathcal{D}} \Psi(D;y) \mu_0(dD)
\hspace{1cm}
\left(y \in \mbox{\boldmath$R$}^m\right).
\]
Denote by $dy$ the Lebesgue measure on $\mbox{\boldmath$R$}^m$ 
and define a measure $Q_0$ on $\mbox{\boldmath$R$}^m$ by 
\[
Q_0(dy) := Z_\Psi(y) dy. 
\]
Then, $Q_0$ is a probability measure 
from Fubini's theorem. 
Define a measure $\mu$ on $\mathcal{D} \times \mbox{\boldmath$R$}^m$ 
by 
\[
\mu(dD, dy) := \frac{1}{Z_\Psi(y)} 
\Psi(D;y) (\mu_0 \otimes Q_0)(dD, dy) 
= \Psi(D;y) \mu_0(dD) \otimes dy,
\]
where $\otimes$ means the product measure. 
Then, $\mu$ is a probability measure.
For $y \in \mbox{\boldmath$R$}^m$ 
define a measure $\mu^y$ on $\mathcal{D}$ by   
\begin{equation}
\mu^y(dD) 
:= \frac{1}{Z_\Psi(y)} \Psi(D;y) \mu_0(dD). 
\label{eq.sect.problems.3} 
\end{equation}
Then, 
$\mu^y$ satisfies the following from Fubini's theorem:
\begin{itemize}
\item
$\mu^y$ is a probability measure on $\mathcal{D}$ for each 
$y \in \mbox{\boldmath$R$}^m,$ 
\item
for every nonnegative measurable function $f$ on   
$\mathcal{D} \times \mbox{\boldmath$R$}^m,$ the function 
\[
y \longmapsto \int_{\mathcal{D}} f(D,y) \mu^y(dD)
\]
is measurable on $\mbox{\boldmath$R$}^m,$ 
\item
for every nonnegative measurable function $f$ on   
$\mathcal{D} \times \mbox{\boldmath$R$}^m,$ the equation
\[
\int_{\mathcal{D} \times R^m} f(D,y) \mu(dD, dy) 
= \int_{R^m}\left\{\int_{\mathcal{D}} f(D,y) \mu^y(dD)\right\}Q_0(dy)
\]
holds.
\end{itemize}
Therefore, $\mu^y$ is a disintegration of $\mu$  
(see \cite{Cha97} p. 292 for the definition of disintegration)
and $\mu^y$ is a probability measure 
of the conditional random variable $D|y$ 
from the viewpoint of Kolmogorov's approach to conditioning
(see \cite{Cha97} p. 293, \cite{Kaw**}). 
Then, we can consider 
(\ref{eq.sect.problems.3}) as a Bayesian formula,  
where 
\begin{itemize}
\item
$\mu_0$ is the prior,
\item
$\Psi(D;y)$ is the likelihood, 
\item
$\mu^y$ is the posterior.
\end{itemize} 
Furthermore, the strong uniqueness of
disintegrations of $\mu$ also holds (see \cite{Kaw**2}).

For given data $y$ and $y' \in \mbox{\boldmath$R$}^m,$ 
define  
\[
\sigma(y, y') := \sup\left\{
\|y - F(D)\| \vee \|y' - F(D)\|:\ D \in \mathcal{D}\right\}.
\]
Using Theorem \ref{main_result}, we can show that 
there exists $C_F > 0$ such that $\|F(D)\| < C_F$ 
for every $D \in \mathcal{D}.$ 
Therefore, 
we have 
$\sigma(y, y') < \infty.$ 
The Hellinger distance between $\mu^y$ and $\mu^{y'}$ is defined by 
\begin{eqnarray*}
d_{\mbox{\footnotesize{Hell}}}\left(\mu^y, \mu^{y'}\right)
& := & \sqrt{\frac{1}{2}
\int_{\mathcal{D}} 
\left\{
\sqrt{\frac{d\mu^y}{d\mu_0}(D)} 
- 
\sqrt{\frac{d\mu^{y'}}{d\mu_0}(D)} 
\right\}^2 
\mu_0(dD)}
\\
& = & 
\sqrt{\frac{1}{2}
\int_{\mathcal{D}} 
\left\{
\frac{1}{\sqrt{Z_\Psi(y)}} \sqrt{\Psi(D;y)}
-
\frac{1}{\sqrt{Z_\Psi(y')}} \sqrt{\Psi(D;y')}
\right\}^2 
\mu_0(dD)}.
\end{eqnarray*}
Note that 
$\mathcal{D}$ 
is not a function space of the same type as in 
\cite{Das13} and \cite{Stu10}. 
However, 
in the same manner as their theory, 
we can obtain the following Lipschitz 
stability of posterior density: 
\begin{equation}
d_{\mbox{\footnotesize{\rm Hell}}}\left(\mu^y, \mu^{y'}\right)
\leq 
\exp\left(\frac{3\sigma(y,y')^2}{4\sigma^2}\right)
\frac{\sigma(y,y')}{\sigma}
\frac{\|y - y'\|}{\sigma}
\label{eq.theorem.sect.stabilities.1.1}
\end{equation}
for given data $y$ and $y' \in \mbox{\boldmath$R$}^m.$ 
We discuss more details in \cite{Kaw**}.

\end{document}